\documentclass[a4paper,11pt]{article}
\usepackage{amssymb,amsthm,bbm,epic,eepic}
\usepackage{a4wide}
\usepackage{amssymb,amsthm}
\usepackage{xypic}
\newcommand{\la}{\lambda}

\def\Lie{\mathrm{Lie}}
\def\Hom{\mathrm{Hom}}
\def\Ext{\mathrm{Ext}}
\def\End{\mathrm{End}}
\def\Ker{\mathrm{Ker}}
\def\Res{\mathrm{Res}}
\def\Ind{\mathrm{Ind}}
\def\Aut{\mathrm{Aut}}
\def\id{\mathrm{id}}
\def\diag{\mathrm{diag}}
\def\U{\mathsf{U}}

\def\SN{\ensuremath{\mathfrak{S}_n}}
\def\SNN{\ensuremath{\mathfrak{S}_{2,n-2}}}
\def\C{\ensuremath{\mathbbm{C}}}
\def\Q{\mathbbm{Q}}
\def\Z{\mathbbm{Z}}
\def\N{\mathbbm{N}}
\def\R{\mathbbm{R}}
\def\bar{\overline}
\def\gl{\mathfrak{gl}}
\def\sl{\mathfrak{sl}}
\def\so{\mathfrak{so}}
\def\ot{\otimes}
\def\k{\mathbbm{k}}
\def\kt{\mathbbm{k}^{\times}}
\def\Ass{\mathbbm{Assoc}}
\def\AN{\ensuremath{\mathfrak{B}_n}}

\def\om{\omega}

\def\undeuxtrois{(\ref{grouplike}), (\ref{inverse}), (\ref{hexagone})}

\def\eps{\epsilon}
\def\g{\mathfrak{g}}
\def\GM{\mathbbm{G}_m}
\def\AFF{\mathbbm{A}}
\def\A{\mathbbm{A}}
\def\un{\mathbbm{1}}
\def\onto{\twoheadrightarrow}
\def\into{\hookrightarrow}
\def\ii{\mathrm{i}}
\def\Gal{\mathrm{Gal}}
\newtheorem{prop}{Proposition}
\newtheorem{defi}{Definition}
\newtheorem{lemma}{Lemma}
\newtheorem{theo}{Theorem}

\def\d{\mathrm{d}}

\def\s{\sigma}

\def\usection#1{\section*{#1}\markboth{\MakeUppercase{#1}}{\MakeUppercase{#1}}\addcontentsline{toc}{section}{#1}}
\def\llabel#1{\label{#1}}
\def\rref#1{\ref{#1}}

\title{{\bf On the representation theory of braid groups}}
\author{Ivan Marin \\ Institut de Math\'ematiques de Jussieu \\ 175 rue du Chevaleret \\ F-75013 Paris \\ marin@math.jussieu.fr}
\date{May 30, 2006}
\begin{document}

\maketitle

\bigskip

\bigskip

\noindent {\bf Abstract.} This work presents an approach towards the representation theory of the braid groups $B_n$. We focus on finite-dimensional representations over the field of Laurent series which can be obtained from representations of infinitesimal braids, with the help of Drinfeld associators. We set a dictionary between representation-theoretic properties of these two structures, and tools to describe the representations thus obtained. We give an explanation for the frequent apparition of unitary structures on classical representations. We introduce new objects --- varieties of braided extensions, infinitesimal quotients --- which are useful in this setting, and analyse several of their properties. Finally, we review the most classical representations of the braid groups, show how they can be obtained by our methods and how this setting enrich our understanding of them.

\medskip


\vfill

\noindent {\bf Table of contents.} \\
\noindent {\bf 1\ \ \ Introduction \hfill 2} \\
\noindent {\bf 2\ \ \ Preliminaries \hfill 6} \\
\noindent {\bf 3\ \ \ Representation theory over $\k[[h]]$ \hfill 11} \\
\noindent {\bf 4\ \ \ Variety of braided extensions \hfill 18} \\
\noindent {\bf 5\ \ \ Restriction to pure braids \hfill 25} \\
\noindent {\bf 6\ \ \ Basic constructions \hfill 33} \\
\noindent {\bf A\ \ \ Formal and real unitarity \hfill 42} \\
\noindent {\bf B\ \ \ Combinatorial aspects \hfill 45} \\
\noindent {\bf References \hfill 48} \\

\eject
\section{Introduction}

\subsection{Overview}

This paper presents a uniform approach to the representations
of the braid groups $B_n$, in order to get as close as possible to a
``representation theory'' for these finitely generated torsion-free
groups. We focus here on the general setting of this approach, whose
benefits include :
\begin{enumerate}
\item Obtention of new representations of the braid groups (see \cite{KZ,THESE} and part 4).
\item Natural explanations and simple proofs of the appearance
of unitary structures on the classical representations of $B_n$
(see section 3.2.2, appendix A and part 6).
\item Determination of the irreducible components of the derived
representations obtained by tensor products and Schur functors
(see \cite{THESE,IRRED,QUOT,LIETRANSP,KRAMMINF} and part 5).
\item Determination of the algebraic hull of the braid groups
inside classical representations (see \cite{LIETRANSP,KRAMMINF} and
section 5.6).
\item A more simple picture of the different sorts of generic
representations of the braid groups (see part 5 et 6).
\item An arithmetic action of absolute Galois groups on
classical representations of $B_n$ (see \cite{GT}), using parts 3 and 5.
\end{enumerate}
Most of these aspects are developped here, applied to several
representations in \cite{QUOT,LIETRANSP,KRAMMINF}, and some
of them are generalized to other kind of so-called Artin groups
in \cite{HECKINF,IRRED,KRAMMINF,DIEDRAUX}. This work is a slightly
revised and enriched version of \cite{ASSOCVIEUX}.

\subsection{Motivations}
Artin braid group on $n$ strands $B_n$ is one of the finitely generated, infinite and torsion-free group which appears most often in mathematics. Concurrently, and in particular in the past thirty years, lots of linear (finite-dimensional, characteristic-zero) representations
appeared in the most diverse contexts. Although it seems illusory to aim at a complete classification of its representations, there is thus a real need for
some understanding of its representation theory. In particular it seems useful to start unifying as much as possible these approaches on some common
ground --- even if it entails a restriction on the range of representations
which can be considered.

The main difficulty in the investigation of the representations of this kind of groups originates in the fact that they usually belong to families depending on transcendant parameters --- contrary to the finite groups, whose representation theory can be studied over the field of algebraic numbers. Rather than study representations on a pure (algebraically closed) field it is thus useful to
consider on this field additional structures, for instance to assume that it is
(the algebraic closure of) the quotient field of some local ring. Let $\k$ be an arbitrary field of characteristic 0 and $K = \k((h))$ be the field of formal series with coefficients in $\k$. The field $K$ is of infinite transcendance degree over $\k$, and the field of matrix coefficients of a linear representation of $B_n$ over $K$ or its algebraic closure $\bar{K}$ is a finitely generated extension of $\k$. 

We choose to investigate representations of $B_n$ over $\k((h))$. It is well-known that such representations can be obtained from the monodromy of flat vector bundles over $X_n/\SN$, where $X_n = \{ (z_1,\dots,z_n) \in \C^n \ \mid \ z_i \neq z_j \Leftrightarrow i \neq j \}$. Indeed, let $\pi : B_n \to \SN$ be the projection of the braid group onto the symmetric group, $P_n = \Ker \, \pi$ the pure braids group ; one has $P_n = \pi_1(X_n)$, $B_n = \pi_1(X_n/\SN)$. The algebraic variety $X_n$ has good properties with respect to $P_n$ : $X_n$ is an Eilenberg-MacLane space and is the complement
of hypersurfaces in some projective space. Associated to it we then find
a holonomy Lie algebra $\mathcal{T}_n$, and it can be shown (see Kohno \cite{KOHNO}) that the representations sufficiently close to the trivial action on the same vector space are monodromies of representations of the
holonomy Lie algebra. The analogous object for $B_n$ is the Hopf algebra $\AN$ defined as the semidirect product of the group algebra of $\SN$ and
the envelopping algebra $\U \mathcal{T}_n$. The geometric properties of $X_n$ imply that the defining relations of $\mathcal{T}_n$ are \emph{homogeneous} with respect to its natural generators. This allows us to introduce a scalar parameter $h$ inside any representation of $\mathcal{T}_n$ or $\AN$. The computation of the monodromy in terms of K.T. Chen's iterated integrals thus yields representations of $B_n$ over the local ring of formal (and even convergent) power series with \emph{complex}
coefficients.

This linear structure associated to the group structure of $B_n$ is far more tractable than its group algebra, and the monodromy operation commutes with direct sums and tensor products. It follows that one may hope to consider $\AN$ as the Lie algebra of $B_n$, in the sense that the representation theory of connected Lie groups essentially relies on its Lie algebra ; in particular we would like to have at disposal a convenient dictionnary which  avoids the explicit computation of the monodromy, in order to decide whether the representation is irreducible, to decompose tensor products, or decide
invariance with respect to bilinear forms.

Transcendant monodromy however presents some major drawbacks. First of all
it can only lead to representations over (finitely generated extensions of)
the field of complex numbers $\C$. But above all it depends on the choice of a base point. This lack of symmetry raises problems in the study of the dual
representations and the invariance with respect to bilinear forms. In order to solve these problems, we use Drinfeld associators. In his article
\cite{DRIN}, Drinfeld expressed the universal monodromy with respect to
some ``base point at infinity'' in a very simple way, in terms of a formal
power series in two non-commuting variables $\Phi_{KZ}(x,y)$. He then
determined which were the algebraic equations that such a series $\Phi(x,y)$, with coefficients in any field $\k$ of characteristic 0,
had to fulfill so that the same formulas of universal monodromy satisfy
the braid relations. Solutions of these equations are called (Drinfeld) associators, and Drinfeld proved in the same article that there exists an associator with coefficients in the field $\Q$ of rational numbers.

The use of associators enables us to remove the two a priori obstacles to
our purpose --- a contrario, the main obstruction to the generalization of
several aspects of this work to other similar groups, such as generalized Artin groups,
is due to the fact that analogues of Drinfeld associators 
have not been defined and extensively studied yet (see however
\cite{DIEDRAUX} and the references there for Artin groups of type $I_2(m)$
and $B$). Drinfeld associators and its analogues may be seen as non-commutative versions
of Chevalley's formal exponentiation (see \cite{CHEVALLEY}).


\subsection{Outline of the results}

In the first place (part 2) we recall the basic notions on braids, infinitesimal braids and associators, which will be of use to us in the sequel. The third part achieves the program of setting a dictionnary between representations of $\AN$ over $\k$ and representations of $B_n$ over $K$. Beyong linear representations, we give a uniform explanation to the apparition of unitary strutures on ``monodromy'' representations. Let us assume $\k \subset \C$. Once we have obtained representations of $B_n$ over the local ring $\k[[h]]$ of formal power series, representations over a pure field can be deduced by at least two means : either by forgetting the local structure on $\k((h))$ --- whose algebraic closure is isomorphic to $\C$ ---
or by specialization in $h \in \k$. This last method may be used only in case the matrix coefficients of the representation are convergent power series. We show in appendix A that, up to a field automorphism, this
situation can always be assumed. We then get unitary representations by this
method.

In view of understanding representation-theoretic aspects, a first advantage of $\AN$ over $B_n$ is that $\AN$ is a semidirect product, whereas
the short exact sequence $1 \to P_n \to B_n \to \SN \to 1$ is not split. We thus define in part 4, for every representation of $\SN$, the variety
of all representations of $\AN$ of which it is the restriction. This yields
a systematic approach to the search for (irreducible) representations of $B_n$. We analyse in this part how much information on the corresponding representations of $B_n$ is contained in this variety. In particular, we show how its factorization through classical quotients of the braid groups can be detected from the original data. Moreover, a third standard manipulation of the local ring $\k[[h]]$, namely reduction modulo $h$, leads to representations of $B_n$ over $\k$
which are iterated extensions of representations of $\SN$ --- in fact,
of the irreducible components of the original representation.

At the other end of the fundamental exact sequence, another operation is the restriction
to the pure braid group of the representation (part 5). The replacement of $P_n$ by the holonomy Lie algebra $\mathcal{T}_n$ enlights noteworthy phenomenons for a large class of representations --- in particular the ``generic'' irreducibility of tensor products, or the irreducibility of the restriction to various
subgroups of $B_n$. We study the consequences of these properties on the
``infinitesimal quotients'' of $B_n$, i.e. Hopf quotients of $\AN$
whose structure explains in particular the decomposition of tensor products.
We moreover indicate (appendix B) how it is possible to get, in particularly
auspicious situations, an explicit (matrix) description of representations of $B_n$ from their infinitesimal version.

The last part in this work (part 6) reviews the most classical constructions
of representations of the braid group : the Iwahori-Hecke and Birman-Wenzl-Murakami algebras, Yang-Baxter representations, and Long's generalization of Magnus induction. We show how our approach enforces the understanding we have on them. In particular we insist on the unitary
structures which appear on the representations and also indicate the general
properties of part 5 that they satisfy.

\subsection{Beyond}

We underline here a few directions of research which are suggested by the uniform approach settled here. We part them in two types. The first one is concerned with the analysis of standard algebraic or geometric structures, shown here to be relevant for the representation theory of the braid group.

Among them are the infinitesimal quotients defined in section 5.6, and the
reductive Lie algebras associated to them. The work of decomposing the Lie
algebra associated to representations of the Temperley-Lieb algebra and
the Iwahori-Hecke algebra was done
by the author in \cite{QUOT,LIETRANSP}. The next step would be to decompose the Lie algebra associated to the
Birman-Wenzl-Murakami algebra of type A --- the part corresponding
to the Krammer representation is already done in \cite{KRAMMINF}.
Another kind of structures are the varieties of braided extensions
$\mathcal{V}^s(M)$ defined in section 4.1. The analysis of these varieties
and the study of their relationship with the space of extensions is another
interesting task. Apart from the examples studied in part 4, the analysis
was carried out for irreducible $M$ (see \cite{KZ,THESE}). A special
question of particular interest is whether all irreducible representations
of all the generic cyclotomic Hecke algebras considered in section 6.1.1
can be deduced in some way from our approach --- answering this question
is a mainly computational but still delicate matter.
A third kind of structure is given by the connection
of ``essentially pure'' representations with the Deligne-Simpson problem
(section 5.4), namely a
description of which roots of the Kac-Moody algebras defined by
W. Crawley-Boevey corresponds to (irreducible, essentially pure)
representations of the braid group.

The second one is concerned with questions which are consequences of this work. A first one is to determinate the field of matrix coefficients for the representations $\widehat{\Phi}(\rho)$ studied here from the infinitesimal datas, and at least its transcendance degree. A second one is to find a criterium on $\rho$ such that $\widehat{\Phi}(\rho)$ is faithful. For the time being, the approaches to faithfulness questions are based on thorough studies of special representations, such as the Burau representation or the representation of the Birman-Wenzl-Murakami algebra which was intensively studied and shown to be faithful, algebraically by D. Krammer and geometrically by S. Bigelow. Now that the linearity of the braid group is known, we need to find a general criterium for faithfulness. It is our hope that such a criterium could be expressed in terms of the infinitesimal
datas --- weak evidence in this direction can be found in section 4.2. However, while the delicate questions of whether the Hecke algebra representation
and the Burau representation for $n=4$ are faithful remain unsettled, it seems prematurate to state even vague conjectures. The third one is about the ``inverse problem'', that is when does a representation $R$ of $B_n$
``comes from'' some representation $\rho$ of $\AN$. This third question is first not to answer this very vague formulation, but to specify what is meant.
Clearly, some representations of $B_n$ cannot be isomorphic to some $\widehat{\Phi}(\rho)$, because their field of matrix coefficients is an algebraic extension of $\Q$ --- however, it may still be possible to get them by specialization (e.g. the representations of the Iwahori-Hecke algebra at roots of unity). Some other representations cannot be equal to some $\widehat{\Phi}(\rho)$, but are twists of some $\widehat{\Phi}(\rho)$ by
elements of $\Gal(K/\k)$ or other field morphisms (see the 4-dimensional example in section 6.1.1, or take any $\widehat{\Phi}(\rho)$ and replace $h$ by $h^2$). We then wonder how many (semisimple) representations of the braid group it is possible to get from representations of the form $\widehat{\Phi}(\rho)$ by an alternating use of these operations, namely manipulations of the field of matrix coefficients
and representation-theoretic operations such as taking direct sums, sub-modules and quotients.

\bigskip

\noindent {\bf General Notations.} By convention, all rings are commutative with unit, all algebras are associative, all representations are finite dimensional, and all Hopf algebras have antipodes. All the fields occuring in the text are of characteristic 0. Unless
otherwise stated, $n$ (the number of strands) and $N$ (the dimension of the representation) are integers with $n \geq 2$, $N \geq 1$, $\k$ is an arbitrary field of characteristic 0. Whenever $\g$ is a Lie algebra, $\U \g$ is its universal envelopping algebra. Whenever
$A$ is a ring, we shall denote by $A^{\times}$ its set of invertible elements, by $M_N(A)$ the set of $N \times N$ matrices and, if $G$ is a group, by $A\, G$ the group algebra of $G$ over $A$. We let
$\GM$, $\A_r$ and $\mu_r$ for $r \geq 1$ be the algebraic varieties whose $A$-points are
$A^{\times}$, $A^r$ and the set of $r$-th roots of 1 in $A$. We use the notation $\mathrm{diag}(a_1,\dots,a_N)$ to designate diagonal matrices in $M_N(A)$ whose diagonal coefficients are $a_1,\dots,a_N$. We let $\bar{K}$
be the algebraic closure of the field $K$, and denote by $\Gal(L/K)$ the
Galois group of an extension of $L$ over $K$. Finally, if $R$ is a representation, either of a group or of a Hopf algebra, we shall denote by $R^{\vee}$ the dual representation, and by $\un$ the trivial representation.

\vfill
\noindent {\bf Acknowledgments.}
I am indebted to Jo\"el Bella\"\i che for many useful hints in the course of this work. I am also grateful to David Madore and Jo\"el Riou for several fruitful discussions.

\eject
\section{Preliminaries}

\subsection{Braids}

Let $n \geq 2$ be an integer. In the sequel, we shall denote by $\SN$
the symmetric group on $n$ letters, and by $s_i$ the transposition $(i \ i\!+\!1)$ for $1 \leq i < n$. The relations between the $s_i$'s give a presentation, found by Moore in
the end of the 19-th century, of the symmetric group :
$$
< s_1,\dots,s_{n-1} \ \mid \ s_i s_{i+1} s_i = s_{i+1} s_i s_{i+1}, \ 
s_i s_j = s_j s_i \ \ |i-j| \geq 2, \ s_i^2 = 1>.
$$
The braid group on $n$ strands $B_n$ is defined by generators $\s_1,\dots,
\s_{n-1}$ and relations $\s_i \s_{i+1} \s_i = \s_{i+1} \s_i \s_{i+1}$
(so-called braid relation), $\s_i \s_j = \s_j \s_i$ if $|i-j| \geq 2$
(so-called locality relation). We let $\pi$ be the (non-split) projection $B_n \to \SN$ given by $\pi(\s_i) = s_i$. The kernel $P_n$ of $\pi$ is by definition the pure (or
colored) braid group, and is generated by the elements $\xi_{ij}
= \s_{j-1} \s_{j-2} \dots \s_{i+1} \s_i^2 \s_{i+1}^{-1} \dots \s_{j-2}^{-1} \s_{j-1}^{-1}$
for $1 \leq i < j \leq n$. Other fundamental elements of $B_n$ are $\gamma_n = (\s_1 \dots \s_{n-1})^n$, $\delta_n = \s_{n-1} \dots \s_2 \s_1^2 \s_2 \dots \s_{n-1}$. It is readily checked that $\gamma_n$ and $\delta_n$ belong to $P_n$. By convention we let $\mathfrak{S}_1 = P_1 = B_1 = \{ e \}$ (the trivial group), and we mention the obvious isomorphisms $P_2 \simeq \Z \simeq B_2$.

The following algebraic facts are known. See the book by Birman \cite{BIRMAN} and the references there for the most classical ones. The groups $P_n$ and $B_n$ are torsion free. For $n \geq 3$ the center of $B_n$ is the same as the
center of $P_n$, is infinite cyclic and
generated by $\gamma_n$. Each $B_n$ (resp. $P_n$) embeds
into $B_{n+1}$ (resp. $P_{n+1}$) by $\s_i \mapsto \s_i$ (resp. $\xi_{ij}
\mapsto \xi_{ij}$) ; this will be refered to as the ``usual inclusion''
$B_n \subset B_{n+1}$. Hence each $\gamma_r$ and $\delta_r$, for
$2 \leq r \leq n$, will be implicitly considered as an element of $B_n$.
The elements $\xi_{1,n}, \xi_{2,n},\dots, \xi_{n-1,n}$ generate a free
normal subgroup $F_n$ in $P_n$, whose centralizer in $P_n$ is the center of
$P_n$,
and the subgroup generated by
$\s_1^2,\dots,\s_{n-1}^2$ is ``locally free'' : this means that the
relations between these elements are generated by the obvious commutation
relations $\sigma_i^2 \sigma_j^2 = \sigma_j^2 \sigma_i^2$ for $|i-j|
\geq 2$ (see \cite{COLLINS}). The elements $\delta_2,\dots,\delta_n$
generate a free \emph{abelian} subgroup in $P_n$. It is readily checked that $\gamma_n = \gamma_{n-1} \delta_{n-1} = \delta_{n-1} \gamma_{n-1}$, and $\delta_n = \s_{n-1} \delta_{n-1} \s_{n-1}$. The action by conjugation of $B_{n+1}$ on its normal subgroup $P_{n+1}$ restricts to an action of $B_n \subset B_{n+1}$ which leaves the
free subgroup generated by $\xi_{1,n+1},\dots,\xi_{n,n+1}$ invariant.
Hence induces an action (so-called Artin action) on a free group on $n$
generators ; this action is known to be faithful.
By the usual inclusion $P_{n-1} \subset P_n$ we have
$P_n = P_{n-1} \ltimes < \xi_{1,n},\dots,\xi_{n-1,n}>$,
hence $P_n$ is an iterated extension of free groups.

We proceed by recalling some well-known quotients of $B_n$. First note that all $\sigma_i$'s belong to the same conjugacy class in $B_n$. Likewise,
all $\xi_{i,j}$ belong to the same conjugacy class in $B_n$ --- but \emph{not} in $P_n$. The quotient of $B_n$ by its subgroup of commutators $[B_n,B_n]$ is infinite cyclic, generated by the image of any $\s_i$. A canonical isomorphism between $B_n/[B_n,B_n]$ and $\Z$ is given by the
so-called length function $l$ ; one has $l(\sigma_i) = 1$ for $1 \leq i \leq n$. The 
quotient of $B_n$ by the subgroup of commutators of $P_n$ will be denoted
by $\widetilde{\SN}$ and refered to as the \emph{enhanced symmetric group}. This probably non-standard terminology is a translation --- and a specialization to type A --- of Tits notion of {\it groupes de Coxeter \'etendus} (see \cite{TITS}). The kernel of the projection $\widetilde{\SN} \to \SN$ is free abelian of rank $n(n-1)/2$, with free generators given by the images of
the $\xi_{ij}$'s in $\widetilde{\SN}$. Finally, let $Hurw_n$ be Hurwitz' group, also called the braid group of the sphere ; it is defined as the quotient of $B_n$ by its normal subgroup generated by $\delta_n$.

\subsection{Infinitesimal braids}

The Lie algebra of infinitesimal pure braids on $n$ strands $\mathcal{T}_n$ is defined over $\Q$ by generators
$t_{ij}$ for $1 \leq i,j \leq n$ and (homogeneous) relations
$$
\left\lbrace
\begin{array}{lcll}
t_{ij}  =  t_{ji} & & & t_{ii} = 0 \\
{ [ } t_{ij},t_{kl} ] & = & 0 & \mbox{if } \#\{ i,j,k,l \} = 4 \\
{ [ } t_{ij},t_{ik} + t_{kj}  ] & = & 0 & \\
\end{array}
\right.
$$
It is endowed with the natural grading defined by all the generators
$t_{ij}$ being of degree 1. The sum of all these generators is easily seen to be central in $\mathcal{T}_n$,
and even to generate its center. 
A remarkable set of homogeneous elements of degree 1 in $\mathcal{T}_n$ is the following. We define, for every $1 \leq r \leq n$,
$$
Y_r = \sum_{1 \leq i < r} t_{ir} = \sum_{1 \leq u <v \leq r} t_{uv}
- \sum_{1 \leq u<v \leq r-1} t_{uv}. $$
thus $Y_1 = 0$ and $Y_2 = t_{12}$. Considering $\mathcal{T}_r$ for $r < n$ as
embedded in $\mathcal{T}_n$ in the natural way, it is clear that each
$Y_r$ belongs to $\mathcal{T}_r$ and commutes to $\mathcal{T}_{r-1}$.
It follows that the $Y_r$'s commute. The elements $t_{ij}$ and $Y_r$ are to be seen as infinitesimal analogues (or, better, ``residues'') of the elements
$\xi_{ij}$ and $\delta_r$ of $P_n$. Likewise, the Lie subalgebra generated by $t_{1,n}, \dots, t_{n-1,n}$
is shown to be free and an ideal of $\mathcal{T}_n$, which
is a semi-direct product of $\mathcal{T}_{n-1}$ and this Lie algebra. 
It is shown in \cite{COHEN} that the centralizer of $\mathcal{F}_n$ is
the center of $\mathcal{T}_n$ --- in the same way that the centralizer of
$F_n$ in $P_n$ is the center of $P_n$.

There exists a natural action of $\SN$ on $\mathcal{T}_n$ (hence also on
$\U \mathcal{T}_n$) which preserves the grading~: it is defined by
$s. t_{ij} = t_{s(i)s(j)}$ for $s \in \SN$. The semi-direct product $\Q \SN \ltimes \U \mathcal{T}_n$ is a Hopf algebra denoted by $\AN$, and we call it the algebra of infinitesimal braids. The sum of the elements $t_{ij}$ remains central in $\AN$. This algebra is naturally graded, with $\deg t_{ij} = 1$ and $\deg s = 0$ for $ s \in \SN$. We denote by $\widehat{\AN}$
its completion with respect to this grading.

Let $\mathcal{L}$ be the (associative) algebra defined by generators $s,Y,Y',t$ and relations
$$
\left\lbrace
\begin{array}{lll}
s^2 =1 & sYs = Y'-t & sts = t \\
{[}Y,Y'{]} = 0 & [Y + Y',t] =0 \\
\end{array} \right.
$$
The algebra $\mathcal{L}$ is graded, with $Y,Y'$ and $t$ of degree 1, and $s$ of degree 0. Consequences of the defining relations are that $Y + Y'$ is central in $\mathcal{L}$ and that
$
[Y',t] = - [Y,t], \ \ [Y',[Y,t]] = - [Y,[Y,t]]
$.
It is easily checked that, for every $n \geq 3$, there
exists a graded algebra
morphism $\mathcal{L} \to \AN$ given by
$$
\begin{array}{lcllcl}
s & \mapsto &(n \! - \! 1 \ n) & t & \mapsto & t_{n-1,n} \\
Y & \mapsto & Y_{n-1} & Y' & \mapsto & Y_n \\
\end{array}
$$
This morphism naturally extends to the completions with respect to the grading, whence a morphism from the completion $\widehat{\mathcal{L}}$
of $\mathcal{L}$ to $\widehat{\AN}$.

\subsection{Associators}

In its fundamental article \cite{DRIN}, Drinfeld defined, for every commutative $\Q$-algebra $\k$,
and every $\la \in \k$, the set $\Ass_{\la}(\k)$ of all formal series $\Phi$
in two non-commuting variables $A$ and $B$, which satisfy the following relations :

\begin{eqnarray}
\Delta(\Phi) = \Phi \hat{\ot} \Phi \label{grouplike} \\
\Phi(B,A) = \Phi(A,B)^{-1} \label{inverse} \\
e^{\la A} \Phi(C,A) e^{\la C} \Phi(B,C) e^{\la B} \Phi(A,B) = 1 \label{hexagone} \\
\Phi(t_{12},t_{23} + t_{24} )\Phi(t_{13} + t_{23},t_{34}) = 
\Phi(t_{23},t_{34})\Phi(t_{12} + t_{13},t_{24} + t_{34})
\Phi(t_{12} ,t_{23}) \label{pentagone}
\end{eqnarray}
where $C=-A-B$. Our $\la$ is Drinfeld's $\mu/2$, i.e. $\Ass_{\la}(\k)$ is Drinfeld's $M_{\mu}(\k)$ (see \cite{DRIN}, \S 5) with $\mu = 2\la$. The equation (\ref{hexagone}) is called the \emph{hexagonal relation}, the
equation (\ref{pentagone}) the \emph{pentagonal relation} --- these names are motivated by MacLane coherence conditions for monoidal categories. 
In
relation (\ref{grouplike}), the symbol $\hat{\otimes}$ denotes the
completed tensor product associated to the identification of the algebra $\k \ll A,B \gg$ of formal series in non-commuting variables $A$ and $B$
with the completed envelopping bigebra of the free Lie algebra on two generators.
The equation (\ref{pentagone}) lies in $\widehat{\mathsf{U} \mathcal{T}}_4$.
Finally, we note that equation (\ref{grouplike}) is equivalent to $\Phi$
being the exponential of a Lie series, let's say 
$\Phi = \exp \Psi$. In particular, let $A$, $B$ and $Z$ belong to some complete graded $\k$-algebra and have zero constant term ; if $Z$ commutes with $A$ and $B$, then $\Psi(A+Z,B) = \Psi(A,B+Z)
= \Psi(A+Z,B+Z) = \Psi(A,B)$,
and
$$
\Phi(A+Z,B) = \Phi(A,B+Z) = \Phi(A+Z,B+Z) = \Phi(A,B).
$$
It follows that equation (\ref{hexagone}) reads, for
$A+B+C=Z$ \emph{central} i.e. commuting with $A$ and $B$ (hence with $C$),
$$
e^{\la A} \Phi(C,A) e^{\la C} \Phi(B,C) e^{\la B} \Phi(A,B) = e^{\la Z}. $$
Also notice that, if $\Phi \in \Ass_{\la}(\k)$ and $\mu \in \k$,
then $\Phi(\mu A, \mu B) \in \Ass_{\la\mu}(\k)$. Moreover, 
for every commutative $\Q$-algebra $\k$ and every $\la \in \k$, $\Ass_{\la}(\k) = 
\Ass_{-\la}(\k)$. 
The very first terms of an associator, as a formal series, are
easy to find.  In particular, the following result is well-known to specialists in the field. For convenience, if $\mathcal{A} = \bigoplus_{n \geq 0} \mathcal{A}_n$ denotes a graded algebra and $x,y \in \mathcal{A}$, we denote by $x \equiv y$ the relation $x-y \in \bigoplus_{n \geq 4} \mathcal{A}_n$
(``equality up to the order 3'').
\begin{prop} \label{DL} 
For any $\Phi \in \Ass_{\la}(\k)$, there exists $\alpha \in \k$ such that
$
\Phi(A,B) \equiv 1 + \frac{\la^2}{6} [A,B] + \alpha \left( [A,[A,B]] - [B,[B,A]] \right) $
\end{prop} 
\begin{proof}
Because of (\ref{grouplike}), we can write $\Phi(A,B) = exp \Psi(A,B)$,
with $\Psi$ a Lie series in $A$ and $B$ with no constant term. Up to the order 2, $\Psi(A,B)$ takes the form $u A + v B$ with $u,v \in \k$, hence
$\Phi(A,B)$ equals
$1 + u A + v B$ plus higher terms. The pentagonal relation (\ref{pentagone})
implies $u = v = 0$, whence $\Psi(A,B)$ and $\Phi(A,B)$ have no linear term. 
Since
there is only one Lie monomial in degree 2, namely $[A,B]$, there exists $u \in \k$ such that
$\Psi(A,B)$ equals $u [A,B]$ plus higher terms, and
$
\Phi(A,B) = \exp \Psi(A,B)$ equals $1 + u [A,B]$ up to the order 2. Replacing $\Phi(A,B)$ by this value in (\ref{hexagone}) we get
$u = \la^2/6$. 

The space of homogeneous Lie polynomials of degree 3 is spanned by
the two Lie monomials $[A,[A,B]]$ and
$[B,[B,A]]$, thus there exists $\alpha,\beta \in \k$ such that
$
\Phi(A,B) \equiv 1 + \frac{\la^2}{6} [A,B] + \alpha [A,[A,B]] + \beta 
[B,[B,A]].$
Using equation (\ref{inverse}), we get $\beta + \alpha = 0$.
\end{proof}

The state of knowledge about these associators is roughly as follows. Drinfeld defined an explicit associator $\Phi_{KZ} \in \Ass_{\ii \pi}(\C)$ and proved by nonconstructive means, not only that $\Ass_1(\Q) \neq \emptyset$,
but also that $\Ass_1^{\circ} (\Q) \neq \emptyset$, where $\Ass_{\la}^{\circ}(\k) = \{ \Phi \in \Ass_{\la}(\k) \ \mid \ \Phi(-A,-B) = \Phi(A,B) \}$ denotes the set of \emph{even} associators
--- in
particular, one can choose $\alpha = 0$ in the expansion of proposition \ref{DL}. 
The associator $\Phi_{KZ}$ is not even, as show its first terms :
$$
\Phi_{KZ}(A,B) \equiv 1 - \zeta(2) [A,B] + \zeta(3) \left( [[A,B],B] - [A,[A,B]] \right) 
$$
where $\zeta$ is the Riemann zeta function. More generally, Le and Murakami
\cite{LEMURAK} gave an explicit formula for its coefficients, which involves
multiple zeta values.

Deep transcendance conjectures on these values lead many people to think that all the algebraic relations over $\Q$ between the multiple zeta values may be equivalent to those implied by the equations (\ref{grouplike}-\ref{pentagone}). There now exists a rather clear algebraic picture of what these algebraic relations should be ; assuming these relations
to be the only ones, it is then possible to guess explicit --- but still
complicated --- formulas for rational even associators. In particular,
Jean Ecalle announced the existence of a conjectural \emph{canonical}
even rational associator.

In this paper we will not anticipate on the developpments in this field, and we will choose arbitrary associators in accordance with Drinfeld existence
theorems. It should however be kept in mind that, although the
formulas which involve even associators are non-explicit in the actual state of knowledge, it should not be so in a (rather) near future. In addition to explicit formulas, another unknown fact which could be of use to us concerns
the \emph{convergence} of these associators, when $\k$ is a complete topological field. It can be shown that $\Phi_{KZ}(A,B)$ is \emph{not}
universally convergent, since the coefficient of $A^{n-1}B$
(resp. $A B^{n-1}$) is $\zeta(n) \geq 1$. However, one easily shows, from Le and Murakami formulas, that
$\Phi_{KZ}(A,B)$ converges at least for $||A|| \leq \frac{1}{4}
$ and $||B|| \leq \frac{1}{4}$, whenever $A$ and $B$ lies in some Banach
algebra. For general associators, we did not hear of any result concerning convergence. We do not know, in particular,
whether there exists convergent
$\Phi \in \Ass_1(\R)$. One only knows that these cannot be universally
convergent (see \cite{DIEDRAUX}). In order to repair this problem, we prove in appendix  A several approximation results which help us making the series converge. Further developments in the field may (or may not) make these tools
less necessary.

\subsection{Drinfeld isomorphisms}

In the same paper \cite{DRIN}, Drinfeld states that, for any $\Phi \in \Ass_{\la}(\k)$,
there exists an homomorphism
from $\k B_n$ to the completion $\widehat{\AN}$ of $\AN$, given by
$$
\sigma_i \mapsto \Phi(t_{i,i+1}, Y_i) s_i e^{\la t_{i,i+1}}
\Phi(Y_i,t_{i,i+1}).
$$
In particular, $\sigma_1 \mapsto s_1 e^{ \la t_{12}}$ and
$
\sigma_2 \mapsto \Phi(t_{23},t_{12}) s_2 e^{\la t_{23}} \Phi(t_{12},t_{23})
$.
Note that equation (\ref{pentagone}) is not needed for the case $n = 3$.
In order to lighten notations we identify here, whenever
$\Phi$ is fixed, $B_n$ with its image in $\widehat{\AN}$.

We recall that 
$
\delta_n =  \sigma_{n-1} \dots \sigma_2 \sigma_1^2 \sigma_2
\dots \sigma_{n-1}$, and make precise the elusive analogy between $\delta_n$ and $Y_n$.

\begin{prop} \label{toraux}
If $\la \in \k$ and $\Phi$ satisfies \undeuxtrois\  then, 
for all $n \geq 2$,
$
\delta_n = e^{2 \la Y_{n+1}}.
$
\end{prop}
\begin{proof} We use induction on $n$. The assertion is clear for $n=2$, since $\delta_2 = \s_1^2$. We now assume $\delta_n = e^{2\la Y_n}$. Let us set for convenience $t = t_{n,n+1}$, $ s= (n \ n \! + \! 1)$, $Y'_n = Y_{n+1} - t$,
$\Phi_n = \Phi(t,Y_n)$ and $\Psi_n = \Phi(t,Y'_n)$. We have $s Y_n s = Y'_n$
and $s\Phi_n s = \Psi_n$.
By definition, $\delta_{n+1} = \s_n \delta_n \s_n$ equals $\Phi_n s e^{\la t} \Phi_n^{-1} e^{2\la Y_n} \Phi_n s e^{\la t} \Phi_n^{-1}$, hence
$$
\delta_{n+1} = \Phi_n e^{\la t} \Psi_n^{-1} \left( e^{\la Y'_n} \right)^2
\Psi_n e^{\la t} \Phi_n^{-1}.
$$
Because $Y_n + Y'_n + t = Y_n + Y_{n+1}$ commutes with $t$, $Y_n$ et
$Y'_n$, the hexagonal relation reads
$$
\begin{array}{lcl}
\Phi_n e^{\la t} \Psi_n^{-1} e^{\la Y'_n} & = & e^{\la (Y_n + Y_{n+1})}
e^{-\la Y_n} \Phi(Y'_n,Y_n) \\
e^{\la Y'_n} \Psi_n e^{\la t} \Phi_n^{-1} & = & \Phi(Y_n,Y'_n) e^{-\la Y_n}
e^{\la (Y_n + Y_{n+1})}  \\
\end{array}
$$
It follows that $\delta_{n+1} = e^{2 \la (Y_n + Y_{n+1} - Y_n)} = e^{2 \la Y_{n+1}}$, and we conclude by induction.
\end{proof}

Let us notice that the image under Drinfeld isomorphism of $\sigma_{n-1} \in B_n$ is the image in $\widehat{\AN}$ of
$$
\sigma = \Phi(t,Y) s e^{\la t} \Phi(Y,t) \in \widehat{\mathcal{L}}$$
by the morphism $\widehat{\mathcal{L}} \to \widehat{\AN}$ that we already defined.
Using the expression of $\Phi(A,B)$ up to the order 3, with parameters $\la$
and $\alpha$, we now get the expression of $\sigma$ up to the order 3.

\begin{prop} If $\Phi(A,B) \equiv 1 + \frac{\la^2}{6} [A,B] + \alpha ([A,[A,B]] - [B,[B,A]])$, then
$$
\sigma \equiv s \left( e^{\la t} + \frac{\la^2}{3} [Y,t] - \alpha [t,[t,Y]]
+ \frac{\la^3}{6} ( t[Y,t] + [Y,t]t) \right)
$$
\end{prop}

\begin{proof}
One has
$$
\begin{array}{lcl}
\Phi(t,Y) s & \equiv & (1 + \frac{\la^2}{6} [t,Y] + \alpha [t,[t,Y]] - \alpha
[Y,[Y,t]])s \\
 & \equiv & s(1 + \frac{\la^2}{6} [t,Y'-t] + \alpha [t,[t,Y'-t]] - \alpha[Y'-t,[Y'-t,t]]) \\
 & \equiv & s(1+ \frac{\la^2}{6} [t,Y'] + \alpha [t,[t,Y']] - \alpha[Y'-t,[Y',t]]) \\
\end{array}
$$
and $[Y'-t,[Y',t]] = [Y',[Y',t]] - [t,[Y',t]] = - [Y',[Y,t]] + [t,[t,Y']] = [Y,[Y,t]] + [t,[t,Y']]$
hence
$$
\Phi(t,Y)s \equiv s( 1+\frac{\la^2}{6} [Y,t] - \alpha [Y,[Y,t]] )
$$
Since $e^{\la t}$ equals $1 + \la t $ plus terms of order at least 2,
$$
\Phi(t,Y) s e^{\la t} \equiv s( e^{\la t}  + \frac{\la^2}{6} [Y,t] - \alpha [Y,[Y,t]] + \frac{\la^3}{6} [Y,t] t)
$$
Since $\Phi(Y,t) \equiv 1 + \frac{\la^2}{6} [Y,t] + \alpha [Y,[Y,t]] - \alpha [t,[t,Y]]$, we get the conclusion.
\end{proof}
These morphisms commute with the ``addition of a strand'' in the sense that the usual inclusions $\k B_n \into \k B_{n+1}$, $\AN \into \widehat{\mathfrak{B}}_{n+1}$ 
along with the Drinfeld morphism $\k B_r \to \widehat{\mathfrak{B}}_r$ for $r = n,n+1$ associated to the same $\Phi \in \Ass_{\la}(\k)$ form a commutative square. If $n =4$ there is in addition a section of the usual inclusion $B_3 \into B_4$ defined by $\sigma_3 \mapsto \sigma_1$, as well as a section of the usual inclusion $\mathfrak{B}_3 \into \mathfrak{B}_4$. This last one is defined by mapping $t_{i,4}$ to $t_{k,l}$
where $\{i,k,l\} = \{1,2,3\}$ and $\mathfrak{S}_4$ to $\mathfrak{S}_3$ by $s_1,s_3 \mapsto s_1$ and $s_2 \mapsto s_2$. The corresponding square involving Drinfeld morphisms also commutes, as can be easily checked.

\vfill\eject
\section{Representation theory over $\k[[h]]$}
Let $\k$ be a field of characteristic 0, $\k[h]$ the ring of polynomials in one indeterminate $h$ over $\k$, $\k[[h]]$ its completion with respect to the $h$-adic topology, i.e. the ring of formal power series in $h$. We shall denote by $\k(h)$ and $\k((h))$ their field of fractions. The multiplicative group $\GM(\k) = \k^{\times}$ acts on $\k[[h]]$ and $\k((h))$ by continuous automorphisms : $\alpha \in \GM(\k)$ sends $f(h)$ on $f(\alpha h)$. Let $A$ be a subring of $\k[[h]]$ containing $\k[h]$, such that the $h$-adic valuation of $\k[[h]]$ induces on $A$ the structure of a discrete valuation ring with local parameter $h$. We call such an $A$ a discrete valuation subring of $\k[[h]]$. We also assume that $A$ is invariant under $\GM(\k)$.
These very strong conditions are fulfilled of course by $A = \k[[h]]$, but
also by $A = \k \{ h \}$, the ring of convergent power series, in case $\k$ is a complete topological field. Let $K$ be the quotient field of $A$, and consider it as embedded in $\k((h))$. For any $x \in \k[[h]]$ and $v \in \k[[h]]^N$, we shall denote by $\bar{x} \in \k$ and $\bar{v} \in \k^N$ their reduction modulo $h$. By default, tensor products are taken over $\k$.

\def\AA{\mathcal{A}}

We study finite dimensional representations of a finitely generated group $G$ over $K$. We first deal with the general case before specializing to
$G = B_n$ and to representations obtained through associators. In case $\k$ has cardinality at most the continuum, the algebraic closure of $\k((h))$ is isomorphic to $\C$, so these representations can be seen (in a highly non-canonical way) as ordinary representations of $G$. However their field of matrix coefficients, as finitely generated extensions of $\k$, have more structure. In particular their transcendance degree $r$ has important meaning,
namely that such a representation can be considered as a $r$ -dimensional family of representations of $G$ over $\k$.

\subsection{Lifting properties}

Our aim here is to establish properties of a representation of $G$ over
$A$ or $K$ from the study of its coefficients in $h$. In order to distinguish these coefficients from the matrix coefficients of the representation,
and because this appellation has some geometrical meaning, we call them ``infinitesimal datas''. It turns out that the $\k$-algebra generated by these infinitesimal datas already contains an important part of the representation-theoretic information about the representation (see propositions \rref{eqform}, \rref{transfirrgen} and \rref{transfabsirr}). In order to prove this, we first need to draw a parenthesis about idempotents in $GL_N(A)$. Notations : whenever $V$ is a $\k$-vector space, we shall
denote by $V[[h]]$ the $A$-module of formal series with coefficients in $V$. If $V$ is endowed with a $\k$-algebra structure, $V[[h]]$ inherits a $\k[[h]]$-algebra structure. For any $v \in V[[h]]$, we let $\bar{v} \in V \subset V[[h]]$ be its constant term.

\subsubsection{Conjugation of idempotents and Hensel lemma}

Let $S \in GL_N(A)$ be a symmetry, i.e. an element such that $S^2 = 1$. Because $\bar{S} \in GL_N(\k)$ is also a symmetry with the same trace, $S$ is conjugated to $\bar{S}$ in $GL_N(K)$. This conjugation also holds in $GL_N(A)$ because of the following identity
$$
(1 + \bar{S} S) S = S + \bar{S} = \bar{S} (\bar{S} S + 1).
$$
Since $\overline{1 + \bar{S}S} = 2$, $1+ \bar{S} S$ is invertible and its
inverse belongs to $GL_N(A)$. Because symmetries and projectors are linked by linear (rational) relations, the same result holds for idempotents. These two facts are particular instances of the following non-commutative version of Hensel lemma.
\begin{prop} \llabel{henselnoncomm}
Let $\mathcal{A}$ be an associative $\k$-algebra with unit, $Q \in \k[X]$
a polynomial in one indeterminate, and $a \in \mathcal{A}[[h]]$. If
$Q(a) = 0$ and $Q'(\bar{a}) \in \mathcal{A}^{\times}$,
then $a$ is conjugated to $\bar{a}$ in $\mathcal{A}[[h]]$ by a (non-commutative) polynomial in $a$ and $\bar{a}$ with coefficients in $\k$.
\end{prop}
\begin{proof}
For any $n \geq 0$, we let $I_n = \sum_{r=0}^n a^{n-r} \bar{a}^r \in \mathcal{A}[[h]]$. In particular $I_0 = 1$, $I_1 = a + \bar{a}$. We have $I_n \bar{a} - \bar{a}^{n+1} = a I_n - a^{n+1}$. Let us write $Q(X) = \sum_{r=0}^{d} c_r X^r$ with $c_r \in \k$, and define $P = \sum_{r=0}^{d-1} c_{r+1} I_r$. Since
$$
P \bar{a} = \sum_{r=0}^{d-1} c_{r+1} I_r \bar{a} = \sum_{r=0}^{d-1} c_{r+1}
\left( \bar{a}^{r+1} + a I_r - a^{r+1} \right) = aP + Q(\bar{a}) - Q(a)
$$
it follows that $P \bar{a} = aP$, and one easily checks $\bar{P} = Q'(\bar{a})$.
\end{proof}

\subsubsection{General setting}
Let $\AA$ be a $\k$-algebra with unit and $G$ a group. We are interested in triples $(R, \rho,N)$ such that $N$ is a positive integer, $R : G \to GL_N(A)$ and $\rho : \AA \to M_N(\k)$ are representations of $G$ and $\AA$
related by the following conditions :
\begin{itemize}
\item[(i)] $\forall g \in G \ \ \ R(g) \in \rho(\mathcal{A}) \otimes A$
\item[(ii)] $\overline{R(KG) \cap M_N(A)} = \rho(\mathcal{A})$
\end{itemize}
In the sequel, if the integer $N$ is implicit we eventually drop the last index and consider \emph{pairs} $(R,\rho)$ satisfying (i) and (ii). The first condition implies that the coefficients in $h$ of $R$ (its infinitesimal datas) belong to the image of $\rho$ ; it also implies that the images of elements of the group algebra $KG$ which belong to $M_N(A)$ have their reduction modulo $h$ inside $\rho(\AA)$. Condition (ii) forces the reverse inclusion to be true.

We sometimes have to be more specific about this last condition. Let
$\mathfrak{m} = hA$ be the maximal ideal of $A$. Let us assume that we are given a family $\mathcal{S}$ of generators for the $\k$-algebra $\AA$, and a function ${\mathsf{f}} : \mathcal{A} \to KG$. Then consider the following condition on the pair $(R, \rho)$.
$$
\mathrm{(ii)}_{\mathsf{f}} \ \ \ \forall j \in J \ \ R({\mathsf{f}}(a_j)) \in \rho(a_j)
+ M_N(\mathfrak{m})
$$
It is clear that (i) and (ii)$_{\mathsf{f}}$ imply (i) and (ii). Conversely,
if $(R,\rho)$ satisfy (i) and (ii), then $(R,\rho)$ satisfy (i) and (ii)$_{\mathsf{f}}$
for \emph{some} ${\mathsf{f}}$. The specification of the function ${\mathsf{f}}$ is useful in the study of intertwinners. For any triples $(R_1,\rho_1,N_1)$ and $(R_2,\rho_2,N_2)$, we shall denote by $\Hom_{\k}(\rho_1,\rho_2)$ the set of intertwinners of $\rho_1$ and $\rho_2$, and
$$
\Hom_A(R_1,R_2)  =  \Hom_A(A^{N_1},A^{N_2})^G, \ \ 
\Hom_K(R_1,R_2)  =  \Hom_K(K^{N_1},K^{N_2})^G. 
$$
We then have the following result.
\begin{prop} Let $(R_1,\rho_1,N_1)$ and $(R_2,\rho_2,N_2)$ be triples satisfying (i) and (ii)$_{\mathsf{f}}$. Then $\Hom_A(R_1,R_2) = \Hom_{\k}(\rho_1,\rho_2) \otimes A$ and $\Hom_K(R_1,R_2) = \Hom_{\k}(\rho_1,\rho_2) \otimes K$. \llabel{eqform}
\end{prop}
\begin{proof}
We identify $\Hom_A(A^{N_1},A^{N_2})$ with $\Hom_{\k}(\k^{N_1},\k^{N_2}) \otimes A$. Because of condition (i), $\Hom_A(R_1,R_2)$ contains $\Hom_{\k}(\rho_1,\rho_2) \otimes A$. We assume
by contradiction that there exists $e \in \Hom_A(R_1,R_2)$ which do not belong to $\Hom_{\k}(\rho_1,\rho_2) \otimes A$. Since we already proved the reverse inclusion, and since $A$ contains $\k[h]$, we may assume $\bar{e} \not\in \Hom_{\k}(\rho_1,\rho_2)$. Then for every $a \in \AA$ there exists $b \in AG$ and $r \geq 0$ such that $h^{-r} R_i(b) \in \rho_i(a) + M_{N_i}(\mathfrak{m})$ for $ i \in \{ 1,2 \}$. Thus $e R_1(b) = R_2(b) e$ implies $\bar{e} \rho_1(a) = \rho_2(a) \bar{e}$, a contradiction. Finally, the identity
$$
\Hom_{K}(R_1,R_2) = \Hom_A(R_1,R_2)\otimes K = \left( \Hom_{\k} (\rho_1,\rho_2) \otimes_{\k} A \right) \otimes_A K
$$
proves the last assertion.
\end{proof}

\subsubsection{Indecomposability and irreducibility}

Since every couple $(R,\rho)$ satisfying (i) and (ii) satisfies (i) and (ii)$_{\mathsf{f}}$ for some ${\mathsf{f}}$, proposition \rref{eqform} proves in particular that $\End_A(R) = \End_{\k}(\rho) \otimes A$ and $\End_K(R) = \End_{\k}(\rho) \otimes K$. Given such a pair $(R,\rho)$, we would like to lift indecomposability,
i.e. to relate the potential indecomposability of $R$ and $\rho$. Here we have to be careful because, if $\rho$ is indecomposable, then $\End_{\k}(\rho)$ is a local $\k$-algebra and, if $\End_{\k}(\rho)$ contains a non-invertible element $u \neq 0$, then $\End_{K} (R) = \End_{\k}(\rho)\otimes K$ is not local
(indeed, $-u+h$ and $u+h$ are non-invertible elements with invertible sum),
hence $R$ is decomposable over $K$. The convenient definition here is to call $R$ decomposable over $A$ if $A^N$ can be written as a direct sum of two $R(G)$-invariant (free) submodules. We then have
\begin{prop}
Let $(R,\rho)$ be a pair satisfying (i) and (ii). Then $\rho$ is indecomposable iff $R$ is indecomposable over $A$.
\end{prop}
\begin{proof}
If $\rho$ is decomposable, then $R$ is decomposable because of condition (i). Now assume that $R$ is decomposable. This means that there exists a non-trivial idempotent $p \in \End_A(R)$, then conjugated in $GL_N(A)$ to its constant term $\bar{p}$ by proposition \rref{henselnoncomm}. Thus $\bar{p}$ is a non-trivial idempotent in $\End_{\k}(\rho)$ by proposition \rref{eqform}, and $\rho$ is decomposable.
\end{proof}

By definition, a representation $R$ over the local ring $A$ is called (absolutely) irreducible if and only if it is (absolutely) irreducible over its quotient field $K$.

\begin{prop} \llabel{transfirrgen}
Let $(R,\rho)$ be a pair satisfying (i) and (ii). Then $\rho$ is irreducible iff $R$ is irreducible.
\end{prop}
\begin{proof}
Let $U$ be a proper $\AA$-invariant subspace of $\k^N$. Then $U \otimes K \subset \k^N\otimes K = K^N$ is a proper subspace of $K^N$ because $\dim _K U \otimes K = \dim_{\k} U$, and it is a $G$-invariant subspace because of property (i). Thus, if $R$ is irreducible then $\rho$ is irreducible. Conversely, let $E$ be a proper and $G$-invariant subspace of $K^N$, and $E' = E \cap A^N$. Because of condition (i), $A^N$ is a $G$-invariant lattice in $K^N$, and $E'$ is $G$-invariant. Since $E$ is proper, there exists $v \in E \setminus \{ 0 \}$. Since $v$ is non-zero there exists $r \in \Z$ such that $h^r v \in A^N$ and $\bar{h^rv} \neq 0$. It follows that $\bar{E'} \neq \{ 0 \}$. Moreover, if $\bar{E'}$ were equal to $\k^N$, any $\k$-basis of $\k^N$ would lift to a $K$-basis of $K^N$ in $E$, hence $E = K^N$. We thus proved that $\bar{E'}$ is a proper subspace of $\k^N$. Let us choose $a \in \AA$, and prove $\rho(a) \bar{E'} \subset \bar{E'}$. Condition (ii) implies that there exists $b \in KG$ such that $R(b) \in M_N(A)$ and $\bar{R(b)}  = \rho(a)$. Let $v_0 \in \bar{E'}$ and $v$ be a lift of $v_0$ in $(v_0 + \mathfrak{m}^N) \cap E$. We have $R(b) v \in E \cap A^N = E'$, and $\bar{R(b) v} = \bar{R(b)} v_0 = \rho(a) v_0 \in \bar{E'}$. It follows that $\bar{E'}$ is a proper $\AA$-invariant subspace of $\k^N$, and this concludes the proof.
\end{proof}
When $\rho$ is irreducible and $\k$ is algebraically closed, in which case $\rho$ is absolutely irreducible, it is not clear a priori that $R$ is absolutely irreducible, because $\k((h))$ is not algebraically closed.
However, this is true :

\begin{prop} \llabel{transfabsirr}
Let $(R,\rho)$ be a pair satisfying (i) and (ii). Then $\rho$ is absolutely irreducible iff $R$ is absolutely irreducible.
\end{prop}
\begin{proof}
Let us recall that $\rho$ is absolutely irreducible if and only if
it is surjective, and that $R$ is absolutely irreducible if and only if
$R : KG \to M_N(K)$ is surjective. If $R : KG \to M_N(K)$ is surjective,
then for every $m \in M_N(\k) \subset M_N(K)$, there exists $b \in KG$
such that $R(b) = m$. Because of condition (i), this implies that $m$ belongs to the
image of $\rho$. Conversely, let us assume that $\rho$ is surjective. 
There exists then $a_1,\dots,a_{N^2} \in \mathcal{A}$ such that
$\rho(a_1),\dots,\rho(a_{N^2})$ form a $K$-basis of the $K$-vector space $M_N(K) \simeq K^{N^2}$. Let $b_1,
\dots b_{N^2}$ be the corresponding elements in $KG$ given by (ii)
such that $R(b_i) \in M_N(A)$ and $\bar{R(b_i)} = \rho(a_i)$ for $1 \leq i \leq N^2$. 
Because of condition (ii), the determinant of the family
$R(b_1),\dots,R(b_{N^2})$ is an element of $A$ with invertible constant term, hence is invertible in $K$. This family then forms
a $K$-basis of $M_N(K)$, thereby proving that $R$ is absolutely irreducible.
\end{proof}

\subsection{Representations of $B_n$}

Let $\Phi \in \Ass_{\la}(\k)$ with $\la \in \kt$. This associator $\Phi$ defines a morphism from the group algebra $\k B_n$ to $\widehat{\AN}$, which is a Hopf algebra morphism for the canonical structures on $\k B_n$ and $\widehat{\AN}$. On the other hand, if $\rho : \AN \to M_N(\k)$ is a representation of $\AN$ over $\k$, we may extend $\rho$ to a representation
$\rho'$ of $\AN$ over $\k[h]$ by $\rho'(t_{ij}) = h \rho(t_{ij})$, hence to a representation $\rho''$ of $\widehat{\AN}$ over $\k[[h]]$. Let $\widehat{\Phi}(\rho)$ be the induced representation of $B_n$ over $\k[[h]]$. Let $A$ be a discrete valuation subring of $\k[[h]]$ such that
for all $N>0$ and all $\rho : \AN \to M_N(A)$, we have $\widehat{\Phi}(\rho)(B_n) \subset GL_N(A)$. For instance, if $\k = \R$ and $\Phi \in \Ass_{\la}(\R)$ is convergent, $A = \k \{ h \}$ may be chosen. We study general properties of these functors $\widehat{\Phi}$, from the category of finite-dimensional $\AN$-modules to the category of representations of $G$ over $A$.

\subsubsection{Linear representations}

Since the Drinfeld morphism $\k B_N \to \widehat{\AN}$ preserves coproducts and antipodes, the functor $\widehat{\Phi}$ commutes with tensor products and duals : for any representations $\rho_1, \rho_2$ of $\widehat{\AN}$, 
$$
\widehat{\Phi}(\rho_1 \otimes \rho_2) = \widehat{\Phi}(\rho_1) \otimes
\widehat{\Phi}(\rho_2), \ \ \widehat{\Phi}(\rho^{\vee}) = \widehat{\Phi}(\rho)^{\vee}.
$$
Let $\rho : \AN \to M_N(\k)$ be a representation of $\AN$, $R = \widehat{\Phi}(\rho)$. Then the pair $(R,\rho)$ satisfies the properties (i) and (ii) of the previous section, for $\mathcal{A} = \AN$ and $G = B_n$. Moreover, it satifies (ii)$_{\mathsf{f}}$ for a universal function ${\mathsf{f}}$ independant of the given pair. Indeed, let us choose as generators of $\AN$ the set $\mathcal{S} = \{ s_i \ | \ 1 \leq i \leq n-1 \} \cup \{ t_{12} \}$. Then
$R(\sigma_i) \in \rho(s_i) + h M_N(A)$ and $R(\sigma_1^2) \in 1 + 2\la h \rho(t_{12}) + h^2 M_N(A)$, so a convenient ${\mathsf{f}}$ is defined by ${\mathsf{f}}(s_i) = \sigma_i$ and ${\mathsf{f}}(t_{12}) = (\sigma_1^2-1)/2\la h$. Consequences of the above results can be condensed in the following theorem.

\begin{theo} \llabel{theolift}
Let $\Phi \in \Ass_{\la}(\k)$ with $\la \in \kt$. Then $\widehat{\Phi}$ preserves indecomposability, irreducibility and absolute irreducibility. It commutes with direct sums, tensor products and duals. Moreover, for any given representations $\rho_1 : \AN \to M_{N_1}(\k)$ and $\rho_2 : \AN \to M_{N_2}(\k)$, we have
$$
\Hom_A(\widehat{\Phi}(\rho_1),\widehat{\Phi}(\rho_2)) = \Hom_{\k}(\rho_1,\rho_2) \otimes A, \ \ \Hom_K(\widehat{\Phi}(\rho_1),
\widehat{\Phi}(\rho_2)) = \Hom_{\k}(\rho_1,\rho_2) \otimes K.$$
\end{theo}

In addition, note that $\widehat{\Phi}$ is exact, since the functor
$\underline{\ \ } \otimes A$ from the category of finite dimensional vector spaces over $\k$ to the category of free $A$-modules is exact. Finally, $\widehat{\Phi}$ obviously commutes with the natural restrictions corresponding to the removal of the last strand --- with a slight imprecision in the notations, we may write :
$$
\Res_{B_{n-1}} \widehat{\Phi}(\rho) = \widehat{\Phi} \left( \Res_{\mathfrak{B}_{n-1}} \rho \right) .
$$
The determination of the field of matrix coefficients of $\widehat{\Phi}(\rho)$ is a far more delicate matter. A coarse lower bound for its transcendance degree over $\k$ is given by the trivial case $n=2$ : since $B_2$ is isomorphic to $\Z$, we only have to determine the field $L$ of matrix coefficients of $\exp( hX)$ when $X \in M_N(\k)$ is known. At least if $\k$ is algebraically closed, a straightforward use of Jordan canonical form shows that, if the spectrum $Sp(X)$ of $X$ equals $\{ a_1,\dots,a_r \}$, then $L = \k ( e^{h a_1},\dots,e^{h a_r} )$
if $X$ is semisimple and $L = \k (h,e^{ha_1},\dots,e^{h a_r })$
otherwise. In particular, writing $\dim_{\Q} S$ for the dimension of the $\Q$-vector space spanned by $S \subset \k$, its transcendance degree over $\k$ is equal to $\dim_{\Q} Sp(X)$ if $X$ is semisimple and $\dim_{\Q} Sp(X) + 1$ otherwise.

This results extends to general $n$ if $\rho(\mathcal{T}_n)$ is commutative --- then $\Phi(\rho(t_{i,i+1}),\rho(Y_i)) = 1$ for $1 \leq i \leq n-1$, $\widehat{\Phi}(\rho)$ do not depend on $\Phi$ and factorizes through $\widetilde{\SN}$. In the general case, it seems difficult to get a non-trivial upper bound for the transcendance degree ; a thiner lower bound is given by taking into account the spectrum of the $\rho(Y_i)'s$, since $\widehat{\Phi}(\rho)(\delta_i) = \exp(2 \lambda h \rho(Y_i))$. Unfortunately, these eigenvalues are usually rational linear combinations of $Sp(\rho(t_{12}))$
--- see appendix B --- so this may not improve the lower bound at all.

\medskip

A last feature of this construction, and another consequence of the lifting properties is the following. Let $\mathcal{A}$ be one of the envelopping algebras $\U \mathcal{T}_n$, $\U \mathcal{F}_n$, $\U \mathcal{I}_n$, where
we denote by $\mathcal{F}_n$ (resp. $\mathcal{I}_n$) the Lie subalgebra of $\mathcal{T}_n$ generated by the $t_{i,n}$ (resp. $t_{i,i+1}$) for $1 \leq i \leq n-1$, and $G$ be the corresponding group $P_n$, $F_n$ or $I_n$, where we denote by $F_n$ (resp. $I_n$) the subgroup of $P_n$ generated by $\xi_{i,n}$ (resp. $\sigma_i^2$) for $1 \leq i \leq n-1$. Let $\rho : \AN \to M_N(\k)$ be a representation and consider $R = \widehat{\Phi}(\rho)$. It
can be easily checked that
$$
\widehat{\Phi}(\rho)(\xi_{i,j}) \in 1 + 2 \la h \rho(t_{ij}) + h^2 M_N(A)
$$
whence an immediate consequence of the lifting properties is the following.

\begin{prop}
Let $\la \in \kt$, $\Phi \in \Ass_{\la}(\k)$ and $\rho : \AN \to M_N(\k)$. Let $(H,\g)$ be one of the following pairs : $(P_n,\mathcal{T}_n)$,
$(F_n,\mathcal{F}_n)$, $(I_n,\mathcal{I}_n)$. Then $\Res_H \widehat{\Phi}(\rho)$ is indecomposable (resp. irreducible, absolutely irreducible) if and only if $\Res_{\g} \rho$ is so.
\end{prop}

\subsubsection{Orthogonal, symplectic and unitary representations}

An immediate corollary of theorem \rref{theolift} is that the functors $\widehat{\Phi}$ for $\Phi \in \Ass_{\la}(\k)$, $\la \in \kt$, preserve
semisimplicity. It is thus a natural question to ask whether $B_n$-invariant bilinear or sesquilinear form can be detected at the infinitesimal level. This is all the more challenging that it appeared in the past twenty years, starting from the Squier form on the Burau representation, that the most classical representations of the braid group are naturally endowed with a ``unitary'' structure. In order to answer this question, we first need
some notations.

Let $L$ be a field of characteristic 0, $N$ a positive integer and $\beta$ a non-degenerate bilinear form on $L^N$. To any $x \in M_N(L)$ we associate its transpose $x^{\dagger}$ with respect to $\beta$ : $\beta(xu,v) = \beta(u,xv)$. As usual, $x \in M_N(L)$ is called symmetric if $x^{\dagger} = x$, antisymmetric if $x^{\dagger} = -x$ and isometric if $x^{\dagger} x = 1$. If $\beta$ is symmetric, we
denote by $O_N(L,\beta)$ the group of isometries in $GL_N(L)$. If $\beta$ is skew-symmetric, we denote it by $SP_N(L,\beta)$. If $L$ is embedded in some larger field $L'$, then $\beta$ is trivially
extended to a non-degenerate bilinear form $\beta'$ on $(L')^N$. For the sake of simplicity, we write $O_N(L',\beta) = O_N(L',\beta')$ and $SP_N(L',\beta) = SP_N(L',\beta')$. Now assume that $L$ admits an involutive
non-trivial field automorphism $\eps$. Given a non-degenerate symmetric bilinear form $\beta$ on $L^N$, we make $\eps$ act on $M_N(L)$ coefficientwise and define
$$
U_N^{\eps}(L,\beta) = \{ x \in GL_N(L) \ \mid \ \eps x^{\dagger} = x^{-1} \}.
$$
If $\beta$ is implicit or is the standard bilinear form on $L^N$ we simply write $U_N^{\eps}(L)$. In particular, the ordinary unitary group $U_N$ is defined as $U_N^{z \mapsto \bar{z}}(\C)$. Recall that the Squier form of the Burau representation (see \cite{SQUIER}) sends $B_n$ to $U_n^{\eps}(\R(q))$ where $\eps \in \Gal(\R(q)/\R)$ is defined
by $\eps(q) = q^{-1}$, and ordinary unitary representations are obtained
by specializing at $q \in \C$ of modulus 1. On $K$, we will consider the
field automorphism $\eps$ defined by $f(h) \mapsto f(-h)$, i.e. the only
automorphism in $\Gal(K/\k)$ continuous for the $h$-adic topology such that
$\eps(h) = -h$.

\medskip

The conditions on the infinitesimal datas in order for our construction to provide orthogonal, symplectic or unitary representations of $B_n$ are the following ones.

\begin{defi} \llabel{defsym}
A representation $\rho : \AN \to M_N(\k)$ is said to be orthogonal (resp. symplectic) if there exists a symmetric (resp. skew-symmetric) non-degenerate bilinear form $\beta$ on $\k^N$ such that $\rho(t_{12})^{\dagger} = - \rho(t_{12})$ and $\rho(\SN)$ is contained in $O_N(K,\beta)$ (resp. $SP_N(K,\beta)$). In case $\beta$ is symmetric, $\rho$ is said to be unitary
if $\rho(t_{12})^{\dagger} = \rho(t_{12})$ and $\rho(\SN) \subset O_N(\k,\beta)$.
\end{defi}

Note that the condition $\rho(t_{12})^{\dagger} = \pm \rho(t_{12})$ implies that $\rho(t_{ij})^{\dagger} = \pm \rho(t_{ij})$ for all $1 \leq i,j \leq n$, since $\rho(s)$ is isometric for all $s \in \SN$. Under these conditions, we have the following properties.

\begin{prop}
Let $\rho : \AN \to M_N(\k)$ be orthogonal (resp. symplectic) with respect to $\beta$, and $\Phi \in \Ass_{\la}(\k)$ with $\la \in \k$. Then $\widehat{\Phi}(\rho)$ factorizes through $O_N(K,\beta)$ (resp. $SP_N(K,\beta)$).
\end{prop}
\begin{proof}
It is sufficient to show that $\widehat{\Phi}(\rho)(\sigma_i)$ is isometric
for all $1 \leq i < n$. For this, we use its explicit expression. Since $O_N(\k,\beta) \subset O_N(K,\beta)$ and
$SP_N(\k,\beta) \subset SP_N(K,\beta)$, $\rho(s_i)$ is isometric. The same assertion for $\exp(\la h \rho(t_{ij}))$ is standard under these hypothesis. Now we have to show that
$\Phi(hx,hy)$ for $x,y \in M_N(\k)$ is isometric as soon as $x^{\dagger} = -x$, $y^{\dagger} = -y$. This follows
from the fact that $\Phi$ is the exponential of a Lie series $\Psi$,
and that the set of skew-symmetric elements forms a Lie subalgebra of
$M_N(\k) = \gl_n(\k)$.
Hence $\Psi(hx,hy)^{\dagger} = - \Psi(hx,hy)$ and $\Phi(hx,hy)^{\dagger} = \Phi(hx,hy)^{-1}$.
\end{proof}

\begin{prop} \llabel{liftunitaire}
Let $\rho : \AN \to M_N(\k)$ be unitary with respect to $\beta$ and $\Phi \in \Ass_{\la}(\k)$ with $\la \in \k$. Then $\widehat{\Phi}(\rho)$ factorizes through $U_N^{\eps}(K,\beta)$, with $\eps : f(h) \mapsto f(-h)$.
\end{prop}
\begin{proof}
It is sufficient to show that each $\widehat{\Phi}(\rho)(\sigma_i)$
belong to $U_N^{\eps}(K,\beta)$. We have $\rho(s_i) \in O_N(\k,\beta)
\subset U_N^{\eps}(K,\beta)$ and, since $\rho(t_{ij})^{\dagger} = \rho(t_{ij})$,
then $\exp(\la h \rho(t_{i,i+1})) \in U_N^{\eps}(K,\beta)$. It remains to show that $\Phi(hx,hy)$ belong to $U_N^{\eps}(K,\beta)$ as soon as $x$ and $y$ are selfadjoint. We know that $\Phi$ is the exponential of a Lie series $\Psi$. Let us introduce the outer Lie algebra automorphism $\tau : x \mapsto -x^{\dagger}$ of $\gl_N(\k)$, trivially extended to $\gl_N(K) = \gl_N(\k) \otimes K$. One has $\tau(x) = -x$, $\tau(y) = -y$, and
$$
-\Psi(hx,hy)^{\dagger} = \tau \Psi(hx,hy) = \Psi(h\tau x,h \tau y) = \Psi(-hx,-hy)
$$
hence $\Psi(hx,hy)^{\dagger} = -\Psi(-hx,-hy)$ and 
$
\Phi(hx,hy)^{\dagger} = \Phi(-hx,-hy)^{-1} = \eps \Phi(hx,hy)
$
so that $\Phi(hx,hy) \in U_N^{\eps}(K,\beta)$.
\end{proof}

In case $\k$ is a subfield of $\R$, the field of coefficients of $\widehat{\Phi}(\rho)$ is a finitely generated extension of $\k$ contained in $\R((h))$, not necessarily in the field of convergent power series
$\R(\{h \})$. We show in the appendix how this representation may be twisted by field isomorphisms in order to get convergent coefficients. It is then possible to get ordinary unitary representations of $B_n$ from representations of $B_n$ into $U_N^{\eps}(K)$, by specialization in $h$ small and real --- we refer to appendix A for more details on this topic.

\vfill\eject
\section{Variety of braided extensions}
One of the major discomforts in the study of braid group representations is that the extension $1 \to P_n \to B_n \to \SN \to 1$ of the symmetric group is not split. In particular
we can not use the well-known representation theory of the symmetric group 
and study a restriction to the symmetric group of a given representation of $B_n$.

Things change drastically with the infinitesimal version $\AN$, which is a semi-direct product. It then appears an important algebraic variety associated to any $\k \SN$-module $M$. Its set of $\k$-points is formed by all representations
of $\AN$ whose restriction to $\k \SN$ is $M$. For reasons that we shall
make clear later, we call this variety the \emph{variety of braided extensions} of the module $M$.

Notations : for $\underline{\alpha} = (\alpha_1,\dots,\alpha_r) \in \Z_{>0}^r$ of size
$\alpha_1+\dots+\alpha_r = n$, we let $\mathfrak{S}_{\underline{\alpha}}$
be the corresponding Young subgroup of $\SN$. In particular,
$\mathfrak{S}_{2,n-2} \simeq \mathfrak{S}_2 \times \mathfrak{S}_{n-2}$.

\subsection{Definition and elementary properties}

Let $M$ be a finite-dimensional $\k \SN$-module. We define
$$
\mathcal{V}(M) = \{ \rho : \AN \to \End_{\k}(M) \ \mid \ \Res_{\SN} \rho = M \}.
$$
In order to have at disposal a more tractable definition of the same object, we introduce the following algebraic functions on $\End_{\k}(M)$ :
$C_1(x) = [x,(1 \ 3). x + (2 \ 3) . x]$, $C_2(x) = [x, (1 \ 3)(2 \ 4) . x]$, where we denote by $x \mapsto s . x$ the action of $s \in \SN$
on $\End_{\k}(M) = M \otimes M^{\vee}$. We let $C_2$ (resp. $C_1$) be zero
if $n < 4$ (resp. $n < 3$).

\begin{lemma}
The map $\rho \mapsto \rho(t_{12})$ defines a bijection from $\mathcal{V}(M)$ to the set $\{ \tau \in \End_{\SNN}(M) \ \mid \ C_1(\tau) = C_2(\tau) = 0 \}$.
\end{lemma}
\begin{proof}
Let $\mathcal{S}$ be the set introduced in the statement. If $\rho \in \mathcal{V}(M)$
and $\tau = \rho(t_{12})$, we have $C_1(\tau) = \rho([t_{12},t_{13}+t_{23}])$,
$C_2(\tau) = \rho([t_{12},t_{34}])$ so $\tau \in \mathcal{S}$ and the map is well-defined. Since $t_{12}$ and $\SN$ generate $\AN$ it is injective. Surjectivity
is a consequence of Frobenius reciprocity law
$$
\End_{\SNN}(M) = \Hom_{\SNN} \left(\un, \Res_{\SNN} \End_{\k}(M) \right)
\simeq \Hom_{\SN} \left(\Ind_{\SNN}^{\SN} \un, \End_{\k}(M) \right).
$$
Indeed, let $e_{\{i,j\}}$ for $1 \leq i \neq j \leq n$ be a basis of the induction of the trivial representation of $\SNN$, such that the action of $\SN$ is given by $s. e_{\{i,j\}} = e_{\{s(i),s(j)\}}$. It is readily checked that the canonical isomorphism given by Frobenius reciprocity law is such that $\tau \in \End_{\SNN}(M)$ is sent to a morphism $\varphi_{\tau}$ satisfying $\varphi_{\tau}(e_{\{1,2\}}) = \tau$,
whence $\rho(t_{ij}) = \varphi_{\tau}(e_{\{i,j\}})$ defines a preimage of $\tau$ in $\mathcal{V}(M)$ --- all the defining relations of $\mathcal{T}_n$ being consequences of $[t_{12},t_{13}+t_{23}] = [t_{12},t_{34}] = 0$
by $\SN$-invariance. 
\end{proof}

This correspondance thus endows $\mathcal{V}(M)$ with the structure of an affine algebraic variety over the field of definition of $M$. Since
all ordinary representations of the symmetric group can be defined
over $\Q$, we usually consider $\mathcal{V}(M)$ as an affine algebraic
variety over $\Q$. As a set, that is when this algebraic structure is not considered, we simply denote by $\mathcal{V}(M)$  the set of $\k$-points of this algebraic variety.
Moreover, we identify both sets under this correspondance, and denote an
element of $\mathcal{V}(M)$ indifferently by $\rho : \AN \to \End_{\k}(M)$
or $\tau \in \End_{\k}(M)$. We let $N = \dim M$, $A = \k[[h]]$, $K = \k((h))$. Any $\Phi \in \Ass_{\la}(\k)$ sends $\rho \in \mathcal{V}(M)$ to a $K$-point of the variety of $N$-dimensional representations of $B_n$.
The images of $\mathcal{V}(M_1)$ and $\mathcal{V}(M_2)$ under this map are disjoint provided that $M_1$ is not isomorphic to $M_2$. Indeed, let $\chi$ be the character of $\SN$ associated to the $\SN$-module $M$. Then the trace of $\widehat{\Phi}(\rho)(\sigma)$ for $\sigma \in B_n$ belongs to $A$ and its reduction modulo $h$ is $\chi \circ \pi(\sigma)$, hence $\chi$ is determined
by $\widehat{\Phi}(\rho) \in \Hom(B_n,GL_N(K))$
and so is $M$ up to isomorphism. In accordance with the general conventions
about variety of representations, we consider \emph{absolutely} irreducible
representations instead of irreducible ones. In case $\la \in \kt$, the elements
in $\mathcal{V}(M)$ which correspond to these are precisely the
absolutely irreducible representations of $\AN$, i.e. elements
$\rho \in \mathcal{V}(M)$ which are surjective. We denote this subset (resp. open subvariety) of $\mathcal{V}(M)$ by $\mathcal{V}^s(M)$. There exists a scheme morphism $\mathcal{V}(M) \to \AFF_1$ given by the trace. All its fibers are isomorphic and for all $l \in \k$ and $\tau \in \mathcal{V}(M)$
we have $\tau + l \in \mathcal{V}(M)$. These fibers as well as the open
subvariety $\mathcal{V}^s(M)$ are invariant under the action of the
reductive group $G(M) = \Aut_{\SN}(M)$ and the schematic quotient $\mathcal{V}^s(M)//G(M)$ is the natural object of study for classification
purposes.

\medskip

The tensor product of representations of $B_n$ admits as infinitesimal counterparts the morphisms $\mathcal{V}(M_1) \times \mathcal{V}(M_2) \to \mathcal{V}(M_1 \otimes M_2)$ given by $(\tau_1,\tau_2) \mapsto \tau_1 \otimes 1 + 1 \otimes \tau_2$, which restrict to maps
$\mathcal{V}_0(M_1) \times \mathcal{V}_0(M_2) \hookrightarrow \mathcal{V}_0(M_1 \otimes M_2)$, where we let $\mathcal{V}_0(M)$ be the set of
all $\tau \in \mathcal{V}(M)$ with zero trace. These maps are injective on the $\k$-points, and closed immersions at the algebraic level. The dual corresponds to the involution $\tau \mapsto - \tau$ on $\mathcal{V}(M)$.
Also note that $\tau \in \mathcal{V}(M)$ has trace $\alpha \in \k$ iff the corresponding representation of $B_n$ deduced from some $\Phi \in \Ass_{\la}(\k)$ is such that (one whence) every $\sigma_i$ act with determinant $\exp(\la h \alpha)$. In view of the transcendance degree of the field of matrix coefficients, another important information encoded in $\mathcal{V}(M)$ is the dimension of the $\Q$-vector space spanned by the eigenvalues of $\tau \in \mathcal{V}(M)$ --- i.e. the roots of its characteristic polynomial. 
As an illustration, let us look for monic polynomials $P \in \Q[X]$ of degree $N>1$, such that the $\Q$-linear span
of its roots has dimension $N$. First $P$ cannot admit any rational root,
for the coefficient of $X^{N-1}$ in $P$ would provide a rational linear relation. For the same reason, this coefficient cannot be zero. Finally, $P$ of course has to be prime with its derivative $P'$. These conditions
are not yet sufficient (see $P= (X-1)^4-2$). They nevertheless imply that this
$\Q$-linear span is at least 2-dimensional. Let $\alpha_1,\dots,\alpha_N
\in \bar{\Q}$ be the roots of $P$ and $L = \Q(\alpha_1,\dots,\alpha_N)$.
A sufficient additional condition is that $G = \Gal(L/\Q)$ acts
2-transitively on the roots. It is indeed a classical fact that the permutation $\Q G$-module with basis $\{ e_{\alpha_i} \ \mid \ 1 \leq i \leq N \}$ and action $g.e_{\alpha_i} = e_{g(\alpha_i)}$ then splits into
two components, the invariant vector $e_{\alpha_1}+\dots+e_{\alpha_N}$
and a $(N-1)$-dimensional \emph{irreducible} supplement. The kernel of its natural
projection onto the $\Q$-linear span of the roots do not contain the
invariant vector (the coefficient of $X^{N-1}$ is non-zero) and do not
contain its supplement because the image has to be of dimension at least 2.
Then this morphism is bijective and the $\Q$-linear span of the roots
has dimension $N$. 
Other standard properties that can be read on $\mathcal{V}(M)$ are dealt with in the next section.

\medskip

The study of $\mathcal{V}(M) = \mathcal{V}^s(M)$ for irreducible $M$ was carried out
in \cite{THESE} and the result is exposed in \cite{KZ}. For example, if $M$
is the irreducible 2-dimensional representation of $\mathfrak{S}_3$,
there exists matrix models of $M$ such that $(1 \ 2)$ acts as the diagonal matrix $\diag(1,-1)$, and the $\SNN$-invariance condition implies
that $\tau = \diag(a,b)$ for $a,b \in \k$. We then have $G(M) \simeq \GM$ and $\mathcal{V}(M) = \mathcal{V}^s(M) = \mathcal{V}^s(M)//G(M) \simeq \AFF_2$. We recall from this earlier work that, even when $M$ is irreducible, the variety $\mathcal{V}^s(M)$
may admit several irreducible components.
As a concrete example for non-irreducible $M$, we now study the case where $M$ is the sum of the 2-dimensional irreducible representation of $\mathfrak{S}_3$ and of the trivial one. Assuming that $\k$ contains square
roots of 3, we use Young orthogonal models to write down matrices for $s_1,s_2$ and for the general form of $\tau \in \End_{\mathfrak{S}_{2,1}}(M)$.
$$
s_1 = \left( \begin{array}{ccc} 1 & 0 & 0 \\ 0 & 1 & 0 \\ 0 & 0 & -1 \\
\end{array} \right) \ \
s_2 = \left( \begin{array}{ccc} 1 & 0 & 0 \\ 0 & \frac{-1}{2} & \frac{\sqrt{3}}{2} \\ 0 & \frac{\sqrt{3}}{2} & \frac{1}{2} \\
\end{array} \right) \ \
\tau = \left( \begin{array}{ccc} a & b & 0 \\ c & d & 0 \\ 0 & 0 & l \\ \end{array} \right)
$$
The representation is absolutely irreducible iff $b \neq 0$ and $c \neq 0$. Note that two such $\tau$'s define isomorphic representations iff the products $bc$ are the same. An easy calculation shows that, assuming $bc \neq 0$, $\tau \in \mathcal{V}(M)$ iff $2a-d-l =0$. The
spectrum of $\tau$ then equals $\{2a-d,\frac{a+d}{2} \pm \frac{1}{2} \sqrt{(a-d)^2 + 4bc} \}$. In view of the field of matrix coefficients
for the associated representations of $B_n$ we note that, provided $\k$ is large enough, the parameters can be chosen so that these three eigenvalues are linearly independant over $\Q$. For instance, if $\k \supset \Q(\sqrt{2},\sqrt{3})$ we may choose $a=1$, $d=\sqrt{2}$, $bc = 1/\sqrt{2}$.

If the original $\k \SN$-module $M$ is endowed with a $\SN$-invariant
non-degenerate symmetric bilinear form $\beta$, we let $M_{\beta}$ be the corresponding structure and $\mathcal{V}(M_{\beta}) = \{ \tau
\in \mathcal{V}(M) \ \mid \ \tau^{\dagger} = \tau \}$. Elements of $\mathcal{V}(M_{\beta})$ give rise to unitary representation of $B_n$ because of proposition \rref{liftunitaire}. In the example above, with $\beta$ being the standard scalar product on $\k^3$, $\tau \in \mathcal{V}(M_{\beta})$ iff $b=c$. If $\k = \R$ note in addition that, for any $\tau \in \mathcal{V}(M)$ the set
$G(M). \tau \cap \mathcal{V}(M_{\beta})$ is non-empty as soon as $bc > 0$.
By comparison, the elements $\tau \in \mathcal{V}(M)$ giving rise to orthogonal representations, i.e. such that $\tau^{\dagger} = - \tau$, form
a family of smaller dimension : they have to fulfill the conditions $a=d=l=0$, $b=-c$ and their eigenvalues can never be linearly independant over $\Q$.

\medskip

We end this section by a rough geometric description of $\mathcal{V}(M)$ and the related
varieties in case $M$ is a $\Q \SN$-module with $\dim_{\Q} M = 4$
and $n =3$. Our motivation is to illustrate how these varieties may depend on the original $\Q \SN$-module.
We identify partitions of 3 and representations of $\mathfrak{S}_3$ : in particular we respectively denote by $[3]$ and $[1^3]$ the trivial and sign representation.

First consider $M = 2 [2,1]$. Then $\End_{\SNN}(M) \simeq \AFF_8$. Let $p : \AFF_8 \to \AFF_6$ be the map $(y_1,\dots,y_8) \mapsto
(y_1-y_2,y_3-y_4,y_5,\dots,y_8)$ and $S \subset \AFF_6$ be the vanishing
variety of the polynomials $x_3 x_6 - x_4 x_5$, $x_5 x_1 - x_3 x_2$,
$x_4 x_2 - x_1 x_6$. Then $\mathcal{V}(M) \simeq p^{-1}(S)$ and $\mathcal{V}^s(M) = \emptyset$. Let now $M = [2,1] + [3] + [1^3] = [2,1]^{\otimes 2}$. Again
$\End_{\SNN}(M) \simeq \AFF_8$. Let $p : \AFF_8 \to \AFF_6$ be the map
$(y_1,\dots,y_8) \mapsto (y_1+y_2-2y_3,y_1+y_2 - 2y_4, y_5,\dots,y_8)$ and
$S \subset \AFF_6$ be the vanishing variety of the polynomials $y_1 y_6,y_2y_6,y_2y_8, y_1y_7$. Then $\mathcal{V}(M) \simeq p^{-1}(S)$ and
$\mathcal{V}^s(M) \simeq p^{-1}(S')$ where $S' \subset S \subset \AFF_6$ is
$\{ (0,0) \} \times \GM^4$, hence $\mathcal{V}^s(M) \simeq \AFF_2 \times \GM^4$. We have $G(M) \simeq \GM^2$ and $\mathcal{V}^s(M) //G(M) \simeq \AFF_2 \times \GM^2$. Let $(x,y,u,v)$ be a $\k$-point in $\AFF_2 \times \GM^2$. Then the spectrum of $4\tau$ under this isomorphism is $\{ 3x + y \pm \sqrt{(x-y)^2 + u}, 3y+x \pm \sqrt{(x-y)^2+v} \}$. If $\k$ is large enough (e.g. $\k \supset \Q(\sqrt{2},\sqrt{3},\sqrt{5})$)
it is possible to get four eigenvalues linearly independant over $\Q$
(take $x=1$, $y= 1+ \sqrt{3}$, $u=1$, $v=2$). Moreover, $2t_{12} + 2t_{13} + 2t_{23}$ acts by the scalar $3(x+y)$, which is not an integer linear combination of the eigenvalues of $\tau$. It follows that the field of matrix coefficients of $R = \widehat{\Phi}(\rho)$ for $\Phi \in \Ass_{\la}(\k)$, $\la \in \kt$ is larger than the field generated by the eigenvalues of $R(\sigma_1)$, since $R(\gamma_3)$ is then the exponential of $3\la h (x+y)$.

Let us now
consider the morphism $\psi : \mathcal{V}([2,1])^2 \to \mathcal{V}(M)$. The
intersection of the image of $\psi$ with $\mathcal{V}^s(M)$ is projected in $\mathcal{V}^s(M)//G(M) \simeq \AFF_2 \times \GM^2$ on quadruples $(x,y,u,v)$
such that $x=y$, whence the corresponding eigenvalues for $\tau$ cannot be
linearly independant over $\Q$. Likewise, for $M = [2,1] + 2[3]$ it
is not possible to get linearly independant eigenvalues over $\Q$. More precisely,
if $\tau \in \mathcal{V}^s(M)$ and $R$ is a representation of $B_3$
corresponding to it, the spectrum $\{ a,b,c,d \}$ of $R(\sigma_1)$
satisfy the algebraic relation $cd = a^4 b^6$. A similar result of course
holds for $M = [2,1] + 2[1^3]$. Finally, in case $M$ is a sum of $1$-dimensional irreducible representations, one easily gets $\mathcal{V}^s(M) = \emptyset$. It follows that the only mean to get
by these methods a 4-dimensional linear irreducible representation of $B_3$ such that $\sigma_1$ has algebraically independant eigenvalues
is to start with $M = [2,1] + [3]+[1^3]$. We shall use this result later (see section 6.1.1).

\subsection{Quotients of the braid group}

Let $M$ be a $\k \SN$-module, $\tau \in \mathcal{V}^s(M)$ and $\rho : \AN \to \End_{\k}(M)$ be the corresponding representation with $\tau = \rho(t_{12})$. We choose once and for all $\Phi \in \Ass_{\la}(\k)$ with $\la \in \kt$ and let $R = \widehat{\Phi}(\rho)$. We look for conditions on $\tau$ in order to characterize factorizations of $R$ through the usual
quotients $Hurw_n$, $\Z = B_n/[B_n,B_n]$ and $\widetilde{\SN}$.
We also consider the map $\pi \times l$ from $B_n$ to $\SN \times \Z$. We let $T$ be the sum of the $\rho(t_{ij})$ for $1 \leq i < j \leq n$ and recall that $\bar{R}$ designates the reduction modulo $h$ of $R$. We assume that $n \geq 3$ and state the results as a list of lemmas.
\begin{lemma} We have $tr(\tau) = 0$ iff $T = 0$ iff $R$ factorizes through
$B_n/Z(B_n)$.\llabel{quotbncentre}
\end{lemma}
\begin{proof}
Indeed, this element $T$ commutes to $\rho(\AN)$, hence is scalar because $\rho$ is absolutely irreducible. It is then zero if and only if its trace is zero. But since all the $\rho(t_{ij})$'s are conjugated in $\AN$,
it follows that $tr(T) = 0$ iff $tr(\tau) = 0$. The last assertion
is a direct consequence of $R(\gamma_n) = \exp(2 \la hT)$.
\end{proof} 
\begin{lemma} $R$ factorizes through $Hurw_n$ iff $\rho(Y_n) = 0$
iff $\rho(Y_n)$ commutes with $\rho(\SN)$ and $tr(\tau) = 0$. \label{quothurwitz}
\end{lemma}
\begin{proof}
The first part comes from $R(\delta_n) = \exp(2 \la h \rho(Y_n))$. If
$\rho(Y_n) = 0$, it of course commutes with $\rho(\SN)$. Conversely, for $n \geq 3$,
if $\rho(Y_n)$ commutes to $\SN$, it also commutes to $\rho(\AN)$
since it commutes to $\rho(t_{12})$. It is then a scalar, and the equality
$tr(\rho(Y_n))
= (n-1) tr(\tau)$ concludes the proof.
\end{proof}
\begin{lemma} \llabel{quotsnparz}
$\tau$ commutes to $\rho(\SN)$, i.e. $\rho(t_{ij}) = \tau$ for all $1 \leq i < j \leq n$, iff for all $\sigma \in B_n$ we have $R(\sigma) = \bar{R}(\sigma) z^{l(\sigma)}$ for some $z \in GL_N(A)$.
\end{lemma}
\begin{proof}
If $t$ commutes to $\rho(\SN)$, $R(\sigma_i) = \rho(s_i) \exp(\la h t)$
and the conclusion follows with $z = \exp(\la h t)$. Conversely,
if $R(\sigma) = \bar{R}(\sigma) z^{l(\sigma)}$ for all $\sigma \in B_n$,
$R(\sigma_1) = \rho(s_1) z$ implies $z = \exp(\la h t)$. Moreover,
$R(\sigma_i) \in \bar{R}(\sigma_i) ( 1 + \la h \rho(t_{i,i+1}))+ h^2M_N(A)$
and $\bar{R}(\sigma_i) z \in \bar{R}(\sigma_i)(1 + \la h\tau) + h^2 M_N(A)$,
hence $R(\sigma_i) = \bar{R}(\sigma_i) z$ implies $\rho(t_{i,i+1}) = \tau$
for $1 \leq i \leq n$. In particular $\rho(t_{23}) = \tau$, i.e.
$\tau$ commutes with $\rho((1 \ 3))$. Since it already commutes to $\rho(\SNN)$, it also commutes to all $\rho(\SN)$.
\end{proof}
\begin{lemma} \llabel{quotextsym}
The following conditions are equivalent : (i) The image $\rho(\mathcal{T}_n)$ of $\mathcal{T}_n$ is commutative ; (ii) $\rho([t_{12},t_{23}])=0$ ; (iii) $R$ factorizes through $\widetilde{\SN}$.
\end{lemma}
\begin{proof}
Of course $i) \Rightarrow ii)$. Condition $i)$ means that $\rho([t_{ij},t_{kl}]) = 0$
for all subsets $\{i,j,k,l\}$ of $\{1,\dots,n \}$, such that $i \neq j$, $k \neq l$. If $\# \{i,j,k,l\} = 4$, this is part of the definition of $\mathcal{T}_n$ ; if $\# \{i,j,k,l\} = 2$ this only means $\rho([t_{ij},t_{ij}]) = 0$, and if $\# \{i,j,k,l\} = 3$, this expression
is conjugated to $\rho([t_{12},t_{23}])$ by some $\rho(s)$ for $s \in \SN$. Hence $ii) \Rightarrow i)$. Now $i) \Rightarrow iii)$ because in
that case $R(\sigma_i) = \rho(s_i) \exp(\la h \rho(t_{i,i+1}))$. Conversely, since
$R(\sigma_1^2 \sigma_2^2) - R(\sigma_2^2  \sigma_1^2)$
belongs to $ 4 \la^2 h^2 \rho([t_{12},t_{23}]) + h^3M_N(A)
$,
it follows that $iii) \Rightarrow ii)$.
\end{proof}

Note that the hypothesis that $\rho$ is absolutely irreducible was useful only for lemma \rref{quothurwitz}. A characterization of the absolutely irreducible
representations $R$ which factorize through $\Z$ is very easy,
since $M$ must then be 1-dimensional. Lemma \rref{quotsnparz} gives a characterization of a more subtle situation, namely when $R$ factorizes through the non-surjective
map $l \times \pi : B_n \to \Z \times \SN$. Note that such an $R$ in particular factorizes through $\widetilde{\SN}$.  

\subsection{Braided extensions}

\subsubsection{General construction}

Let $A$ be a discrete valuation subring of $\k[[h]]$ and $K$ be its quotient field. Let $R : B_n \to GL_N(A)$ be a representation of the braid group.
As usual, we let $\bar{R}$ be its reduction modulo $h$, and we assume
that $\bar{R}(\sigma_1^2) = 1$, i.e. that $\bar{R}$ factorizes through $\SN$. This condition is of course satisfied if $R = \widehat{\Phi}(\rho)$
for some $\rho : \AN \to M_N(\k)$. More precisely, if $\rho \in \mathcal{V}^s(M)$, $\bar{R}$ corresponds to the $\k \SN$-module $M$. By conjugation in $GL_N(K)$ we can get other representations $S : B_n \to GL_N(A)$. We denote this relationship between $R$ and $S$ by
$R \sim S$. It means that $R$ and $S$ are isomorphic as representations over $K$, but correspond to different $B_n$-invariant lattices in $K^N$. It is a classical fact that the simple modules appearing in the Jordan-H\"older decomposition of $\bar{R}$ and $\bar{S}$ are the same.
\begin{prop} If $R$ is irreducible and $\bar{R}$ is not irreducible,
then there exists $S \sim R$ such that $\bar{S}$ do not factorize through
the symmetric group.
\end{prop}
\begin{proof} Under these hypothesis it is known (see \cite{JOEL,JOELGRAF}) that there
exists $S \sim R$ such that $\bar{S}$ is not semisimple. Since
all representations of the symmetric group are semisimple, the conclusion follows.
\end{proof}
More generally, let $M$ be a $\k \SN$-module, and let $(M_i)_{i \in I}$
be its simple submodules. Let $\tau \in \mathcal{V}^s(M)$, $\Phi \in \Ass_{\la}(\k)$ with $\la \in \kt$ such that $\widehat{\Phi}(\rho)(B_n) \subset GL_N(A)$. The relationship between the reductions modulo $h$ of $\widehat{\Phi}(\rho)$ given by different choices of lattices in $K^N$
on the one hand and $B_n$-extensions of the modules $M_i$ on the other hand
enter a general setting extensively studied by Bella\"\i che and Graftieaux in \cite{JOEL,JOELGRAF}. This is the justification for
the name we have chosen of ``variety of braided extensions'' of the
given $\k \SN$-module $M$. A case of particular interest is when $M = M_1 \oplus M_2$, with $M_1$ and $M_2$ irreducible. Then any $\tau \in \mathcal{V}^s(M)$ allows us to define, following the original idea of Ribet \cite{RIBET} in
the framework of algebraic number theory, elements in $\Ext_{B_n}(M_1,M_2)$
and $\Ext_{B_n}(M_2,M_1)$. We now explain this in more detail.

Assume $M = M_1 \oplus M_2$ with $M_1$ and $M_2$ being simple $\k \SN$-modules.
Let $\tau \in \mathcal{V}^s(M)$ and $\rho : \AN \to \End_{\k}(M)$ be the
associated representation of $\AN$. We choose $\Phi \in \Ass_{\la}(\k)$
with $\la \in \kt$. 
For simplicity we assume that $A = \k[[h]]$.
If, according to the decomposition $M= M_1 \oplus M_2$, we write $\rho(s_i)$
and $\rho(t_{i,i+1})$ for $1 \leq i \leq n$ in matrix form as follows
$$
\rho(s_i) = \left( \begin{array}{cc} s_i^{(1)} & 0 \\ 0 & s_i^{(2)} \\ \end{array} \right) \ \ \ \rho(t_{i,i+1}) = \left( \begin{array}{cc} \tau_{i,i+1}^{(11)} & \tau_{i,i+1}^{(12)} \\ \tau_{i,i+1}^{(21)} & \tau_{i,i+1}^{(22)} \\ \end{array} \right)
$$
then, on the lattices $hM_1 \oplus M_2$ and $M_1 \oplus h M_2$, the reduction of $\widehat{\Phi}(\rho)(\sigma_i)$ modulo $h$ reads
$$
\left( \begin{array}{cc} s_i^{(1)} & \la s_i^1 \tau_{i,i+1}^{(12)} \\
0 & s_i^{(2)} \\ \end{array} \right) \ \mbox{ and }
\ \left( \begin{array}{cc} s_i^{(1)} & 0 \\ \la s_i^{(2)} \tau_{i,i+1}^{(21)} & s_i^{(2)} \\ \end{array} \right)
$$
respectively. Note that these expressions do \emph{not} depend on the choice of $\Phi \in \Ass_{\la}(\k)$. Moreover, since $\AN$ is generated by the $s_i$'s and $t_{i,i+1}$'s, the fact that $\tau$ belongs to $\mathcal{V}^s(M)$
implies that there exists $i$ and $j$ such that $\tau_{i,i+1}^{(12)}$
and $\tau_{j,j+1}^{(21)}$ are non zero, hence that the elements obtained in $\Ext_{B_n}(M_1,M_2)$ and $\Ext_{B_n}(M_2,M_1)$ are non trivial. The most
elementary case is the following fact :
\begin{prop} \label{vsvide}
Let $B$ and $C$ be two non-trivial $\k \SN$-module. Then
$$ \Hom_{\SNN}(C,B) = 0 \Rightarrow
 \mathcal{V}^s(B \oplus C) = \emptyset, \ \ 
\Ext_{B_n}(B,C) = \Ext_{B_n}(C,B) = 0.
$$
\end{prop}
\begin{proof}
We assume $\Hom_{\SNN}(C,B) = 0$. Since all $\k \SNN$-modules are semisimple, we also have $\Hom_{\SNN}(B,C) = 0$. The $\SNN$-endomorphisms of $B \oplus C$ decompose as follows : 
$\End_{\SNN}(B \oplus C)$ equals
$$ \End_{\SNN}(B) \oplus
\Hom_{\SNN}(B,C) \oplus
\Hom_{\SNN}(C,B) \oplus
\End_{\SNN}(C,C).
$$
Then every $\tau \in \mathcal{V}(B \oplus C)$ leaves $B$ and $C$ invariant, hence can
not belong to $\mathcal{V}^s(B \oplus C)$. Let us now consider
an element of $\Ext_{B_n}(C,B)$, i.e. a representation of $B_n$
on the vector space $B \oplus C$ which can be written in matrix form
$$
\sigma \mapsto \left( \begin{array}{cc} \sigma^B & \sigma^{CB} \\
0 & \sigma^C \\ \end{array} \right)
$$
Since $B$ and $C$ are representations of $\SN$, we have $(\sigma_1^B)^2 = 1$ and $(\sigma_1^C)^2 = 1$. Moreover,
$\sigma_1^2$ commutes to $B_{2} \times B_{n-2} \subset B_n$. Since the image
of $\sigma_1^2$ has the form $\left( \begin{array}{cc} 1 & x \\
0 & 1 \\ \end{array} \right)$, this means that $\sigma^B x = x \sigma^C$
for every $\sigma \in B_2 \times B_{n-2}$, hence that $x \in
\Hom_{\SNN}(C,B)$. By hypothesis we then have $x = 0$, hence $\sigma_1^2 = 1$. Since the braids $\sigma_i^2$ for $1 \leq i \leq n$ are conjugated in $B_n$,
it follows that $\sigma_i^2 = 1$ for all $1 \leq i \leq n$ and
the representation factorizes through $\SN$, hence is semi-simple and
$\Ext_{B_n}(C,B) = 0$. Similarly, $\Ext_{B_n}(B,C) = 0$.
\end{proof}

In general, elements in $\Ext_{B_n}(C,B)$ can then be constructed from elements in $\mathcal{V}^s(B \oplus C)$. Note that $\mathcal{V}^s(B \oplus C)$ is not a $\k$-vector space in general but only a cone, contrary
to $\Ext_{B_n}(B , C)$. A natural question, which we leave open, is whether these elements
generate $\Ext_{B_n}(B,C)$ as a $\k$-vector space.

\subsubsection{Example : extensions between hooks}

Let $\un$ be the trivial representation of $\SN$, $E$ its permutation representation over $\k$. Then $E = \un \oplus U$, where $U$ is an irreducible representation of $\SN$ corresponding to the partition $[n-1,1]$ of $n$. Let $(e_1,\dots,e_n)$
be a basis for $E$ such that $s. e_i = e_{s(i)}$ for $s \in \SN$. The variety $\mathcal{V}^s(E)$ was studied in \cite{THESE}. It is shown there that $\mathcal{V}^s(E)$ is of the form $\kt \tau$ for some non-zero $\tau \in \mathcal{V}^s(E)$, and that the corresponding irreducible representation of $B_n$ deduced from any associator $\Phi \in \Ass_{\la}(\k)$
is given by
$$
\left\lbrace \begin{array}{lcl}
\sigma_k . e_k & = & q e_{k+1} \\
\sigma_k . e_{k+1} & = & q e_k \\
\sigma_k. e_r & = & e_r \mbox{ if } r \not\in \{k,k+1\} \\
\end{array} \right.
$$
with $q = e^{\alpha h}$, $\alpha$ depending on $\la$ and on the choice of $\tau \in \mathcal{V}^s(E)$. By the general procedure described above, we get the
following non-split indecomposable representations of $B_n$ over $\k$~:
$$
\left\lbrace \begin{array}{lcl}
\sigma_k.f_k & = & f_{k+1} + \alpha \frac{n-2}{n^2} v \\
\sigma_k.f_{k+1} & = & f_k + \alpha \frac{n-2}{n^2} v \\
\sigma_k . f_r & = & f_r - \frac{2\alpha}{n^2} v  \\
\end{array} \right.
\left\lbrace \begin{array}{lcl}
\sigma_k.g_k & = & g_{k+1} - \frac{2\alpha}{n^2} w + \frac{\alpha}{n} (g_k+g_{k+1})  \\
\sigma_k.g_{k+1} & = & g_k - \frac{2\alpha}{n^2} w + \frac{\alpha}{n} (g_k+g_{k+1}) \\
\sigma_k . g_r & = & g_r - \frac{2\alpha}{n^2} w + \frac{\alpha}{n} (g_k+g_{k+1})  \\
\end{array} \right.
$$
with the following notations : $r \not\in  \{k,k+1 \}$, $(f_1,\dots,f_n)$ and $(g_1,\dots,g_n)$ are some basis of $E$, $v = f_1 + \dots + f_n$, $w = g_1+ \dots + g_n$. The $p$-th exterior power $\Lambda^p U$ of $U$ is irreducible under $\SN$ and corresponds to the hook
diagram $[n-p,1^p]$. One has $\Lambda^p V = \Lambda^{p-1} U \oplus \Lambda^p U$. From these non-trivial elements in $\Ext_{B_n}(\un,U)$ and $\Ext_{B_n}(U,\un)$ we then get non-trivial elements in $\Ext_{B_n}(\Lambda^{p-1}U,\Lambda^p U)$ and $\Ext_{B_n}(\Lambda^p U, \Lambda^{p-1}U)$. For the sake of concision, we introduce the Koszul differential on the exterior algebra $\Lambda^{\bullet}E$, defined by
$$
\d(e_1 \wedge \dots \wedge e_m) = \sum_{j=1}^m (-1)^{j+1} e_1 \wedge \dots
\wedge \widehat{e_j} \wedge \dots \wedge e_m.
$$
If $f$ is of the form $f_{i_1} \wedge \dots \wedge f_{i_m}$ for some $m$
such that $\{ i_1,\dots,i_m \} \cap \{ k,k+1 \} = \emptyset$,
then the action of $\sigma_k$ in the first case is given by
$$
\left\lbrace
\begin{array}{lcl}
\sigma_k.f & = & f - \frac{2\alpha}{n^2} v \wedge \d f\\
\sigma_k.(f_k \wedge f) & = & f_{k+1} \wedge f - \frac{2\alpha}{n^2} v \wedge \d f + \alpha \frac{n-2}{n^2} v \wedge f \\
\sigma_k.(f_{k+1} \wedge f) & = & f_k \wedge f - \frac{2 \alpha}{n^2} v \wedge \d f + \alpha \frac{n-2}{n^2} v \wedge f \\
\sigma_k. (f_k \wedge f_{k+1} \wedge f) & = &
f_{k+1} \wedge f_k \wedge f + \alpha \frac{n-2}{n^2} v \wedge (f_k - f_{k+1}) - \frac{2 \alpha}{n^2}(f_{k+1} - f_k) \wedge v \wedge \d f \\
\end{array}
\right.
$$
In order to describe the action of $\sigma_k$ in the second case, we introduce
a permutation action of $\SN$ on the vector space $\Lambda^{\bullet} E$, defined by
$s.(g_{i_1}\wedge \dots \wedge g_{i_m}) = g_{s(i_1)} \wedge \dots \wedge g_{s(i_m)}$ for $s \in \SN$.
We then have
$$
\sigma_k.g = s_k.g + \left( \frac{\alpha}{n}(g_k + g_{k+1})  - \frac{2\alpha}{n^2} w \right) \wedge \d g
$$
for all $g \in \Lambda^{\bullet} E$.
Note that, since the restriction to $\SNN$ of $\Lambda^{p} U$ and $\Lambda^q U$ have no isomorphic irreducible component for $|p-q| \geq 2$, it follows from proposition \ref{vsvide} that
$\Ext_{B_n}(\Lambda^p U, \Lambda^q U) = 0$ in these cases. The cohomology groups of $B_n$ in its representation $U$ were computed
by Vassiliev in \cite{VASS}. In particular, he there showed that
$\Ext_{B_n}(\un,U) = H^1(B_n,U) \simeq \k$ for $n \geq 3$.

\vfill\eject
\section{Restriction to pure braids}

Keeping in mind the short exact sequence $1 \to P_n \to B_n \to \SN \to 1$,
we now focus on the restriction of representations of $B_n$ to its
normal subgroup of finite index $P_n$. Being considered as an iterated
extension of free groups, $P_n$ do not have a well-understood
representation theory to build on. On the contrary, it seems hopeless
to classify its representations --- whatever this could mean. Nevertheless, as it is a normal subgroup
of finite index, most irreducible representations of $B_n$ are
already irreducible under the action of $P_n$. Moreover, its infinitesimal
structure is a graded Lie algebra, that is a more elementary structure than the
infinitesimal braid algebra $\AN$. It follows that studying the restriction to $\mathcal{T}_n$ helps in understanding the decomposition of
semisimple representations into irreducible components. Since $\mathcal{T}_n$ is generated by its
elements of degree 1, it is first a natural task to investigate the
image of these under a representation $\rho$ of $\AN$, and their role in
irreducibility properties.

We show that the restriction of $\rho$ to this subspace of
homogeneous degree 1 elements of $\mathcal{T}_n$ is injective, unless the associated $B_n$-representation factorizes
through one of the classical quotients of the braid group. We then focus
on representation-theoretic properties which are \emph{tannakian}, i.e.
such that the collection of representations satisfying these forms
a tannakian subcategory of the category of representations.
For instance, semisimplicity is a tannakian property. In order to enlarge
this setting, we define a notion of generically tannakian properties, making use of the natural $\GM(\k)$-action on $\mathcal{T}_n$. Two generically tannakian properties naturally appear. One of them is essential purity. Essentially pure irreducible representations of $\AN$  
restrict to irreducible representations of $\mathcal{T}_n$. As a matter
of fact, they are already irreducible under the action of several
Lie subalgebras of $\mathcal{T}_n$. This means that the associated
representations of $B_n$ are irreducible under the action of
several subgroups. The other one is agregation. This notion was first defined in \cite{IRRED}, and is a generalization of the fact that tensor products of irreducible representations of a free group are ``usually'' irreducible.

Another interesting property of a representation $\rho$ is when the
eigenvalues of $\rho(t_{12})$ form a simplex over $\Q$.
We conclude this part by studying the incidence of all these properties on
the Hopf
quotients of $\AN$ induced by $\rho$,
and mention the connection with the algebraic hull of the
pure braid group in the corresponding representations. For convenience, an associator $\Phi \in \Ass_{\la}(\k)$, $\la \in \kt$, integers $n \geq 2$, $N \geq 1$ are
chosen once and for all.

\subsection{Linear independance}

Let $\rho : \AN \to M_N(\k)$ and $R = \widehat{\Phi}(\rho)$,
for some $\Phi \in \Ass_{\la}(\k)$ and $\la \in \kt$. We showed in the last part that it is possible to characterize the potential factorization of
$R$ through $B_n/Z(B_n)$, $\SN \times \Z$, $Hurw_n$ or
$\widetilde{\SN}$ in terms of $\rho$. We show here that,
if $R$ do not factorize through one of these groups, then the $\k$-vector
space spanned by $\rho(t_{ij})$ for $1 \leq i < j \leq n$ is $n(n-1)/2$-dimensional. In particular, since $B_2 = \widetilde{\mathfrak{S}_2}$, we assume $n \geq 3$.

We first need to recall several elementary facts about representations
of the symmetric group. Le $E$ be the $n(n-1)/2$-dimensional $\k \SN$-module with basis $e_{ij} = e_{\{i,j\}}$ and action $s. e_{ij}
= e_{s(i),s(j)}$. It is the induced representation of the trivial representation of the Young subgroup $\SNN$. Littlewood-Richardson rule
shows that provided $n \geq 4$ it is the sum of three irreducible non-isomorphic
components respectively corresponding to the partitions $[n]$, $[n-1,1]$
and $[n-2,2]$ whose dimensions are $1$, $n-1$ and $n(n-3)/2$. If $n =3$ there are only two irreducible components, corresponding to $[n]$ and $[n-1,1]$. In both cases, the first
component is generated by the sum $y$ of all the $e_{ij}$'s for $1 \leq i < j \leq n$, and the
direct sum of the first two components admits as basis the $y_i$'s defined
for $1 \leq i \leq n$ as $y_i = \sum_{j \neq i} e_{ij}$.

\begin{prop}
If $R$ do not factorize through $B_n/Z(B_n)$, $Hurw_n$ or
$\widetilde{\SN}$, then the $\rho(t_{ij})$'s are linearly independant.
\end{prop}
\begin{proof}
We let $L$ be the sub-$\SN$-module of $\End(M)$ linearly 
generated by the $\rho(t_{ij})$'s. 
We use the notations above and let $u_i = \varphi(y_i)$. There is a natural surjective $\SN$-morphism
$\varphi : E \to L$ defined by $\varphi(e_{ij}) = t_{ij}$. 
Consequences of the characterizations in section 4.2 are the following ones.
Lemma \rref{quotbncentre} implies $\varphi(y) \neq 0$. Lemma \rref{quothurwitz} implies that $u_n$ cannot commute to every $s\in \SN$,
hence there exists $i,j$ such that $u_i -u_j \neq 0$. Lemma \rref{quotsnparz} implies that $L$ has dimension at least 2. It follows
that $\Ker \, \varphi$, as a $\SN$-subrepresentation of $E$, cannot contain
the invariant vector of $E$ and may contain at most one of the remaining
irreducible components. One of them is generated by the elements $y_i - y_j$, hence the kernel cannot contain this one. It follows that, if $\varphi$ were not surjective $L$ would be of dimension $n$, generated by $u_1,\dots,u_n$, and $n \geq 4$. We assume this and show that it would contradict lemma \rref{quotextsym} --- that is, the family $u_1,\dots,u_n$
would be commutative.

Let $u_1+ \dots + u_n = \mu \neq 0$. Since $\rho(t_{12})$ is
a linear $\SNN$-invariant linear combination of the $u_i$'s, there exists
$b,c\in\k$ such that $\rho(t_{12}) = b(u_1+u_2) + c \mu$, hence $\rho(t_{ij}) = b(u_i+u_j) + c \mu$ for all $1 \leq i< j \leq n$. It follows
that $u_i = (n-2)bu_i + ((n-1)c + b) \mu$ whence $b = 1/(n-2) \neq 0$
and $u_i + u_j$ is a linear combination of $\rho(t_{ij})$ and $\mu$.
The infinitesimal braid relations and $\SN$-invariance then imply that $[u_i+u_j,u_k] = 0$ whenever $\# \{i,j,k \} = 3$. Choosing $i,j$ such that
$\# \{ 1,2,i,j\} = 4$ one has
$
[u_i,u_j] = [u_i+u_1,u_j] - [u_1,u_j] = -[u_1,u_j] = [u_1,u_2]$ and similarly $[u_j,u_i] = [u_1,u_2]$. By $\SN$-invariance it follows that $[u_i,u_j] = 0$ for all $i,j$, hence a contradiction.
\end{proof}

\subsection{Tannakian and generically tannakian properties}

Let $\mathcal{A}$ be a graded Hopf algebra over $\k$. We shall study here
special properties that representations of $\mathcal{A}$ may have,
in order to apply this setting to the case where $\mathcal{A} = \U \mathcal{T}_n$ or $\mathcal{A} = \AN$. For a given property $(P)$, we study in particular
their stability by the standard operations in representation theory,
namely taking sub-modules and quotients, dual and tensor products, and
extension of the base field. 
\begin{defi}
A property (P) is said to be \emph{tannakian} if it is stable by 
taking sub-modules and quotients, dual and tensor products, and
extension of the base field. In other words, if \\
\indent 1) \ If $\rho$ satifies (P), then so does $\rho^{\vee}$, and all its
quotients and sub-modules. \\
\indent 2) \ If $\rho^1,\dots,\rho^r$ satifies (P), then so does $\rho^1
\otimes \dots \otimes \rho^r$. \\
\indent 3) \ If $\rho$ satifies (P), then $\rho \otimes L$ satifies (P) over any extension $L$
of $\k$.
\end{defi}

If a property $(P)$ is tannakian, the full sub-category
of the (left) $\mathcal{A}$-modules that satisfy (P) inherits 
the tannakian structure, hence corresponds to a pro-algebraic
group scheme over $\bar{\k}$. In particular, the property of factorizing
through some Hopf quotient of $\mathcal{A}$ is the basic example of a tannakian property.

We also define a larger class of properties, by relaxing condition 2).
Since $\mathcal{A}$ is graded, it is naturally endowed with
an action of $\GM(\k)$ such that $\alpha \in \GM(\k)$ sends
an homogeneous element $x \in \mathcal{A}$ of degree $r$ to $\alpha^r x$.
If $\rho$ is a representation of $\mathcal{A}$ and $\alpha \in \GM(\k)$, we let $\rho_{\alpha}$ be the representation twisted
by this action : $\rho_{\alpha}(x) = \rho(\alpha.x)$. In case $\mathcal{A} = \AN$ and $R = \widehat{\Phi}(\rho)$,
$\widehat{\Phi}(\rho_{\alpha})$ is deduced from $R$
by the field automorphism $\alpha \in \GM(\k) \subset \Gal(\k((h))/\k)$,
this inclusion being defined by $\alpha.f(h) = f(\alpha h)$ for
$f \in \k((h))$.  

\begin{defi} A property (P) is said to be \emph{generically tannakian} if
1), 2') and 3) are satified, where
\begin{itemize}\item[2')] If $\rho^1,\dots,\rho^r$ satifies (P), then so does $\rho^1_{\alpha_1}
\otimes \dots \otimes \rho^r_{\alpha_r}$ for generic
$(\alpha_1,\dots,\alpha_r) \in \GM(\k)^r$.
\end{itemize}
\end{defi}
\noindent In this definition, the genericity has to be understood with respect to the Zariski
topology on $\GM(\k) = \mathrm{Spec} \, \k[t^{-1},t]$.

We end this section by studying properties of $\U \mathcal{T}_n$ and
$\AN$ which are tannakian. All the proofs are standard. The first property under consideration is semisimplicity. We first need a lemma.
\begin{lemma} \label{restrss}
Let $G$ be a finite group acting on a Lie algebra $\g$ by automorphisms. A representation of $\U_G \g = \k G \ltimes \U \g$ is semi-simple if and only if its
restriction to $\U \g$ is semi-simple.
\end{lemma}
\begin{proof}
Let $\rho : \U_G \g \to M_N(\k)$ be a representation of $\U_G \g$,
and $\rho'$ be its restriction to $\U \g$. We first assume that $\rho$ is irreducible. Let us choose an irreducible subrepresentation of $\rho'$ and let $U$ be the underlying subspace of $\k^N$. Then any vector space $\rho(s) U$ for $s \in G$ is $\rho'$-invariant and irreducible as a representation of $\U \g$,
and the sum of all of these is a non-zero $\rho$-invariant subspace of $\k^N$. Since $\rho$ is irreducible, this subspace is $\k^N$ and $\rho'$ is semisimple. If $\rho$ is semisimple, the same conclusion of course holds.
Conversely, we assume that $\rho'$ is semisimple. Let $\psi$ be a sub-representation of $\rho$ and $U$ be the underlying subspace. Let $i: U \hookrightarrow \k^N$ be the canonical inclusion. Since $\rho'$ is semisimple, $i$ admits a left inverse $j$ as a $\U \g$-module. Let $J = (\sum
\psi(s) \circ j \circ \rho(s)^{-1})/(\#G)$, with the sum being taken over all elements of the finite group $G$. One has $J \circ i = \id_U$ and $J$
is a morphism of $\U_G \g$-modules. Then $\rho$ is semisimple
and the lemma is proved.
\end{proof}

\begin{prop}
The property of being semisimple is tannakian for $\mathcal{A} = \U \mathcal{T}_n$ and $\mathcal{A} = \AN$.
\end{prop}
\begin{proof}
For every algebra $\mathcal{A}$, every sub-quotient of a semisimple module is semisimple. Condition 3) is
moreover satified because $\k$ is perfect, being of characteristic 0. In case $\mathcal{A} = \U \mathcal{T}_n$, the dual is semisimple and condition 2) is satified because $\mathcal{A}$ is an envelopping algebra. The general case follows from lemma \rref{restrss}.
\end{proof}

We call a representation $\rho$ of $\AN$ unitary if is unitary with respect to some non-degenerate symmetric bilinear form.

\begin{prop}
The property of being unitary is tannakian.
\end{prop}
\begin{proof}
Since any unitary representation is semisimple, condition 1) is true.
Condition 2) holds because, if $\rho,\rho'$ are unitary with respect to
$\beta,\beta'$, then $\rho \otimes \rho'$ is unitary with respect to $\beta \otimes \beta'$, which is symmetric and non degenerate. Condition 3) is clear.
\end{proof}

\subsection{Agregating representations}

A less standard property is the property of agregation.
This is our first example of a generically tannakian property. We recall from \cite{IRRED} its definition and its main application.

\begin{defi} An element in $M_N(\k)$ is called regular
it is diagonalizable with distinct eigenvalues. A representation $\rho$ of $\mathcal{T}_n$ is said to be \emph{agregating} if $\rho(V)$ contains
a regular element, where $V$ is the vector space generated by the $\rho(t_{ij})$. A representation of $\AN$ is called agregating if its
restriction to $\U \mathcal{T}_n$ is agregating.
\end{defi}
\begin{prop}
The property of being agregating is generically tannakian.
\end{prop}
\noindent It is sufficient to show this for $\U \mathcal{T}_n$. See \cite{IRRED}
for a proof. A consequence of this property is the following~:
\begin{prop} (see \cite{IRRED}) \llabel{agregtens}
Let $\rho$ be a representation of $U \mathcal{T}_n$. If $\rho$ is agregating and irreducible, then it is absolutely irreducible. If $\rho^1,\dots,\rho^r$
are
agregating 
representations
of $\U \mathcal{T}_n$, then $\rho^1_{\alpha_1} \otimes \dots \otimes \rho^r_{\alpha_r}$ is irreducible (resp. indecomposable) for
generic values of $\alpha_1,\dots,\alpha_r \in \GM(\k)$ if and only if all the $\rho^i$'s are so.
\end{prop}

\subsection{Essential purity}

In this section, we let $s = (1 \ 2)$, $t = t_{12}$. As an associative algebra with unit, $\mathfrak{B}_2$ is generated by $s$ and $t$, with relations $s^2=1$,
$st = ts$.

\begin{defi}
A representation $\rho$ of $\mathfrak{B}_2$ is said to be \emph{essentially pure} if $\rho(s)$ is a polynomial in $\rho(t)$. A representation
of $\AN$ is called essentially pure if its restriction to $\mathfrak{B}_2$ is so.
\end{defi}

If $\rho(t)$ is diagonalizable, it means that $\rho(s)$ acts by a scalar on
each sub-vector space $\Ker(\rho(t) - \alpha)$ for $\alpha \in \k$. Note that a representation $\rho$ of $\AN$ is essentially pure if and only if
\emph{any} $\rho((i \ j))$ is a polynomial in $\rho(t_{ij})$.

This property has special incidences on the restrictions
of a representation of $B_n$ or $\AN$. The notations $F_n$, $I_n$,
$\mathcal{F}_n$ and $\mathcal{I}_n$ were defined in section 3.2.1.

\begin{prop}
Let $\rho$ be a representation of $\AN$, and $\Phi \in \Ass_{\la}(\k)$
with $\la \in \kt$.
If $\rho$ is essentially pure, then the irreducible components
of $\rho$ (resp. $\widehat{\Phi}(\rho)$) are irreducible under the action
of $\mathcal{T}_n$, $\mathcal{F}_n$ and $\mathcal{I}_n$
(resp. $P_n$, $F_n$ and $I_n$).
\end{prop}
\begin{proof}
If $\rho$ is essentially pure, the image $\rho(( i \ j))$ of
each transposition is a polynomial in $\rho(t_{ij})$. We can suppose
that $\rho$ is irreducible. Since the transpositions generate $\SN$
, its restriction to $\mathcal{T}_n$ is irreducible. Since
the transpositions of the special form $(i \ n)$ (resp. $(i \ i \! + \! 1)$)
generate $\SN$, its restriction to $\mathcal{F}_n$ (resp. $\mathcal{I}_n$) is also irreducible.
The consequences for $\widehat{\Phi}(\rho)$ follow from proposition \ref{transfirrgen}.
\end{proof}

Since $F_n$ is a free group, it follows from this that every essentially pure representation of the braid group $B_n$ yields a solution of the Deligne-Simpson problem on $n$-tuples, in the particular case where the first $n-1$ conjugacy classes are the same. Recall that this problem concerns the classification of $n$-tuples $(C_1,\dots,C_n)$ of conjugacy classes in $GL_N(K)$ such that there exists a $n$-tuple $(X_1,\dots,X_n) \in GL_N(K)^n$ satisfying $X_i \in C_i$ for $1 \leq i \leq n$, $K^N$ is irreducible under $X_1,\dots,X_n$, and $X_1 \dots X_n = 1$. Here $X_i = R(\xi_{i,n})$ for $1 \leq i \leq n-1$ and $X_n = R((\xi_{1,n} \dots \xi_{n-1,n})^{-1})$. Conversely, a simple solution of this problem in this particular case could be of great benefit for the representation theory of the braid group. An infinitesimal variant of this problem was considered and solved by specialists in the field, namely the classification of $n$-tuples $(c_1,\dots,c_n)$ of conjugacy classes in $M_N(\k)$ such that there exists a $n$-tuple $(x_1,\dots,x_n) \in M_N(\k)^n$ satisfying $x_i \in c_i$ for $1 \leq i \leq n$, $\k^N$ is irreducible under $x_1,\dots,x_n$ and $x_1+\dots+x_n = 0$. Here $x_i = \rho(t_{i,n})$ for $1 \leq i \leq n-1$ and $x_n = -\rho(t_{1n} + \dots + t_{n-1,n})$. Let us mention the solution recently given by Crawley-Boevey \cite{CRAWLEY} in terms of representations
of quivers and the associated Kac-Moody algebra. The solutions
given by essentially pure irreducible representations of $B_n$ correspond to a special kind of quivers, and it would be most interesting to decide which roots of the
associated Kac-Moody algebra are concerned by this construction.
Let us moreover notice that
the so-called Katz rigidity
index is very easy to compute from a given essentially pure $\AN$-representation, and gives an interesting invariant for such
representations --- for instance, the Burau representation leads to a rigid tuple in this
sense, and it is the only one doing so among representations of the
Iwahori-Hecke algebra of type A (see \cite{BURAURIGIDE}).

\begin{lemma} The representation $\rho : \AN \to M_N(\k)$ is essentially pure iff it is essentially pure over $\bar{\k}$. \llabel{esspurebar}
\end{lemma}
\begin{proof} If $\rho$ is essentially pure, it is in particular essentially pure over $\bar{\k}$. Conversely, let $P \in \bar{\k}[X]$ be such that
$\rho(s) = P(\rho(t))$, and let $L$ be smallest normal extension of $\k$ containing the coefficients of $P$. Since $L$ is Galois and finite
one has $\rho(s) = Q(\rho(t))$ with $Q = (\sum \sigma P)/(\# \Gal(L/\k)) \in \k[X]$, the sum being taken over all $\sigma \in \Gal(L/\k)$.
\end{proof}

\begin{prop}
If $\rho^1,\dots,\rho^r$ is a family of essentially pure
representations then, for generic $\alpha_1,\dots,\alpha_r \in \GM(\k)$, $\rho^{\alpha}  = \rho^1_{\alpha_1} \otimes \dots \otimes \rho^r_{\alpha_r}$ is essentially pure.
\end{prop}
\begin{proof}
Because of lemma \rref{esspurebar}, we may assume that $\k$ is algebraically closed. Let $s^i = \rho^i(s)$ and  $t^i = \rho^i(t)$. Let $d^i$ be the semisimple component of the endomorphism $t^i$. Then $t^{\alpha} = \rho^{\alpha}(t)$ is the sum of the elements $\alpha_i 1 \otimes \dots \otimes t^i \otimes
\dots 1$ for $1 \leq i \leq r$, and the sum $d^{\alpha}$ of the elements $\alpha_i 1 \otimes \dots \otimes d^i \otimes
\dots 1$ for $1 \leq i \leq r$ is its semisimple component. 
In particular $Sp(t^{\alpha}) = Sp(d^{\alpha})$ and there exists a polynomial $Q_{\alpha} \in \k[X]$ depending on $\alpha_1,\dots,\alpha_r$ 
such that
$d^{\alpha} = Q_{\alpha}(t^{\alpha})$. 
Since the representations $\rho^i$ are essentially pure,
there exists for all $1 \leq i \leq r$ polynomials $P_i \in \k[X]$ such that $s^i = P_i(t^i)$ --- hence $s^i = P_i(d^i)$ since $s^i$ is semisimple.
It follows that there exists set-theoretic maps ${\mathsf{f}}_i : Sp(t^i)=Sp(d^i) \to Sp(s^i)$
such that $ d^i v = \mu v$ implies $s^{i}v = {\mathsf{f}}_i (\mu) v$.

Let $v_i$, $1 \leq i \leq r$, be a non-zero vector in the underlying vector space of $\rho^i$, and $v = v_1 \otimes \dots \otimes v_r$. There exists a set-theoretic map ${\mathsf{f}} = \prod {\mathsf{f}}_i : Sp(t^1)\times \dots
\times Sp(t^r) \to Sp(s^{\alpha})$ such that $\forall i \ d^i v_i = \mu_i v_i$ implies $s^{\alpha}v = {\mathsf{f}}(\mu_1,\dots,\mu_r) v$.
For generic $(\alpha_1,\dots,\alpha_r) \in \GM(\k)$ we may assume that the
condition $d^i v_i = \mu_i v_i$ for all $1 \leq i \leq r$ is equivalent to the condition $d^{\alpha} v = ( \sum \alpha_i \mu_i) v$. Hence there exists a map ${\mathsf{f}}' : Sp(d^{\alpha}) \to Sp(s^{\alpha})$ such that $d^{\alpha} v = \mu v \Rightarrow s^{\alpha} v = {\mathsf{f}}'(\mu)$. Since
$d^{\alpha}$ is semisimple and $\k$ is algebraically closed, it follows that $s^{\alpha} = P(d^{\alpha})$ for some $P \in \k[X]$, hence $s^{\alpha} = (P \circ Q_{\alpha})(t^{\alpha})$ and $\rho^{\alpha}$ is essentially pure.
\end{proof}

It easily follows :

\begin{prop}
The property of being essentially pure is generically tannakian.
\end{prop}

\subsection{Simplicial representations}

We study here $N$-dimensional representations of $B_n$ such that 
the image of $\sigma_1$ or, equivalently, of $\sigma_1^2$,
is semisimple and its field of matrix coefficients is a purely
transcendantal extension of $\k$.  As we saw before, if $R = \widehat{\Phi}(\rho)$ for some $\Phi \in \Ass_{\la}(\k)$ with $\la \in \kt$ and $\k$ is algebraically closed, this means that the eigenvalues of $\rho(t_{12})$ have no $\Q$-linear relations between them --- in particular, this condition only depends
on the restriction to $\mathcal{T}_n$, hence fits our preoccupations in this part. Since a rescaling of $R$ by multiplication of some scalar (equivalently, a rescaling of $\rho$ by addition of some scalar) do not change the representation in a significant way, the convenient setting at the infinitesimal level is not linear algebra but affine geometry over $\Q$. In this framework, a simplex will refer to the structure of $\Q$-affine vector space of $\k$.

\begin{defi}
A representation $\rho : \AN \to M_N(\k)$ is called \emph{simplicial} if
$\rho(t_{12})$ is diagonalizable over $\k$ and its spectrum forms a
simplex.
\end{defi}

We introduce in appendix B a representation-theoretic criterium ensuring that the interior of this simplex contains the eigenvalues of $2T/n(n-1)$, with $T$ being the sum of the $\rho(t_{ij})$'s --- thus illustrating the relevance
of this geometric definition.


For any Hopf algebra $\mathcal{A}$ and representation $\rho$ of $\mathcal{A}$, the derived representations of $\rho$ are the representations
of the form $\rho^{\otimes p} \otimes (\rho^{\vee})^{\otimes q}$ for $p,q \geq 0$, their subrepresentations and quotients. Besides its influence 
on the field of matrix coefficient, the most interesting
aspect of simpliciality is that it makes essential purity extend to derived representations.

\begin{prop}
Let $\rho$ be simplicial and 
essentially pure. Then all its derived 
representations are essentially pure. \llabel{simplpure}
\end{prop}
\begin{proof}
Let $\mathcal{S} \subset \k$ be the spectrum of $\rho(t_{12})$.
Essential purity means that there exists a set-theoretic map
$\varepsilon : \mathcal{S} \to \{ -1,1\}$ such that, for all $\alpha \in S, \ v \in \k^N$, $\rho(t_{12})v = \alpha 
v \Rightarrow \rho(s) =
\varepsilon(\alpha)v$.
Let $r \geq 0$. We first show that $\rho^{\otimes r}$ is essentially pure. We let $t = t_{12}$, $s = (1 \ 2)$ as before. Let
$\underline{v} \in (\k^N)^{\otimes r}$ be a pure tensor formed by tensoring several eigenvectors for $t$ : $n_1$ 
with eigenvalue $\alpha_1$, \dots, $n_r$ with eigenvalue $\alpha_r$. One 
has,
$$
\rho^{\otimes r}(t) \underline{v} = \left( \sum_{i=1}^{\# S} n_i \alpha_i \right) 
\underline{v} \ \ \ \rho^{\otimes r} (s) \underline{v} = \prod_{i=1}^r \varepsilon(\alpha_i)^{n_i}
\underline{v}.
$$
Since $\sum n_i = r$ and the $\alpha_i's$ are affine independant, the value
of $\rho^{\otimes r}(t_{12})$ on $\underline{v}$ determines the $n_i$'s as barycentric
coordinates, hence the value of $\rho^{\otimes r}(s)$ on $\underline{v}$. Since the tensors like $\underline{v}$ form a basis of $(\k^N)^{\otimes r}$
which are eigenvectors for $\rho^{\otimes r}(t_{12})$ it follows that
$\rho^{\otimes r}(s)$ is a polynomial of $\rho^{\otimes r}(t_{12})$
and
$\rho^{\otimes r}$ is essentially pure.
The same more 
generally holds for $\hat{\rho} = \rho^{\otimes p} \otimes 
(\rho^{*})^{\otimes q}$ for $p,q \geq 0$ : the eigenvalues of $\hat{\rho}(t)$ are then of the 
form $\sum (n_i - m_i) \alpha_i$, with $\sum n_i = p$, $\sum m_i = q$,
and uniquely determine the $(n_i-m_i)$'s, hence the parity of the $(n_i + 
m_i)$'s
and the corresponding action $\prod \varepsilon(\alpha_i)^{n_i + m_i}$
of $\hat{\rho}(s)$. Since essential purity is generically tannakian and
in particular satisfies condition 1) the conclusion follows.
\end{proof}

\subsection{Infinitesimal quotients and algebraic hulls}

Let $\rho: \AN \to M_N(\k)$ be a representation and $\mathcal{H}_{\rho} = \rho(\mathcal{T}_n)$ the Lie subalgebra of $\gl_N(\k)$ generated by the $\rho(t_{ij})$'s. The sum $T$  of the $\rho(t_{ij})$'s belong to the center of
$\mathcal{H}_{\rho}$. The symmetric group acts on $\mathcal{H}_{\rho}$ by automorphisms and the Lie algebra morphism $\mathcal{T}_n \to \mathcal{H}_{\rho}$ extends to a surjective Hopf algebras morphism
$\AN \to \mathcal{Q}_{\rho}$ where $\mathcal{Q}_{\rho} = \k \SN \ltimes
\U \mathcal{H}_{\rho}$ is called the \emph{infinitesimal quotient} of $B_n$
associated to $\rho$.
Note that, if $\rho$ can be decomposed as $\rho_1 \oplus \dots \oplus \rho_r$, then $\mathcal{H}_{\rho}$ naturally embeds into $\mathcal{H}_{\rho_1} \oplus \dots \oplus \mathcal{H}_{\rho_r}$ by an $\SN$-equivariant Lie algebra morphism, hence $\mathcal{Q}_{\rho}$
embeds into $\mathcal{Q}_{\rho_1} \oplus \dots \oplus \mathcal{Q}_{\rho_r}$
as a Hopf algebra.

The representation-theoretic conditions introduced above, when satisfied
by $\rho$, have the following consequences on $\mathcal{H}_{\rho}$.

\begin{prop} Let $\rho$ be essentially pure and semisimple.
Then $\mathcal{H}_{\rho}$ is a reductive Lie algebra
and $Z(\mathcal{H}_{\rho}) = \rho(Z(\mathcal{T}_n)) = \k T$. In particular,
if $\rho(t_{12}) \in \gl_N(\k)$ has non-zero trace,
then $\mathcal{H}_{\rho}$ has 1-dimensional center.
Moreover, if $\rho$ is agregating then $\mathcal{H}_{\rho}$ is split over
$\k$.
\llabel{Hrhoreduct}
\end{prop}
\begin{proof}
The representation $\rho$ induces a faithful representation $\rho'$ of
$\mathcal{H}_{\rho}$ on $\k^N$. Since $\rho$ is essentially pure, every simple component of $\rho$ is irreducible under the action of $\mathcal{T}_n$, hence $\rho'$ is faithful and semisimple, and $\mathcal{H}_{\rho}$ is reductive.

One always has $\rho(Z(\mathcal{T}_n)) = \k T$. Since
$\mathcal{H}_{\rho}$ is reductive, one has
$\mathcal{H}_{\rho} = Z(\mathcal{H}_{\rho})
\times [\mathcal{H}_{\rho},\mathcal{H}_{\rho}]$
and $\rho([\mathcal{T}_n,\mathcal{T}_n]) = [\mathcal{H}_{\rho},
\mathcal{H}_{\rho}]$. If $\mathcal{T}_n^1$ denotes the subspace generated
by the $t_{ij}'s$ inside $\mathcal{T}_n$ one has
$\mathcal{T}_{n} =\mathcal{T}_{n}^1 \times [\mathcal{T}_{n},
\mathcal{T}_{n}]$ as vector spaces and $\SN$-modules,
and $\mathcal{H}_{\rho} = Z(\mathcal{H}_{\rho})
\times [\mathcal{H}_{\rho},\mathcal{H}_{\rho}]$ is also a decomposition
as $\SN$-modules because $\rho$ is essentially pure
hence $\SN$ acts trivially on $Z(\mathcal{H}_{\rho})$. We denote by $p$ the corresponding $\SN$-equivariant
projection $\mathcal{H}_{\rho}
 \onto Z(\mathcal{H}_{\rho})$.
It follows that
$Z(\mathcal{H}_{\rho}) \subset p \circ \rho(\mathcal{T}_{n}^1)^{\rho(\SN)}$
hence it is a quotient of the $\SN$-module $\mathcal{T}_n^1$.
Because of the decomposition of this latter semisimple $\SN$-module
one has $p \circ \rho(\mathcal{T}_{n}^1)^{\rho(\SN)} = p\circ \rho((\mathcal{T}_n^1)
^{\SN}) = p(\k T)$. Hence $\k T \subset Z(\mathcal{H}_{\rho} \subset
p(\k T)$ with $\dim \k T \geq \dim(p(\k T))$ and the conclusion follows.
Since the trace of $\rho(T)$ is $\frac{n(n-1)}{2}$-times the trace
of $\rho(t_{12})$, if the latter is non-zero then $\k T$ has dimension 1.

If $\rho$ is agregating, there exists $x \in \mathcal{T}_n$ such that $u = \rho(x)$ is a regular element in $\gl_N(\k)$. It is then a regular element
of the Lie algebra $\mathcal{H}_{\rho}$ in the traditional sense and a
Cartan subalgebra of $\mathcal{H}_{\rho}$ is given by
$\mathcal{C}_{\rho} = \{ y \in \mathcal{H}_{\rho} \ \mid \ \exists m \
(\mathrm{ad} \ u)^m y =0 \}$. But since
$\mathcal{H}_{\rho}$ is reductive, then $\mathcal{C}_{\rho}$ is
commutative. It follows that $\mathcal{C}_{\rho} \subset M_N(\k)$ is a set
of polynomials in $u$, hence all elements in $\mathcal{C}_{\rho}$ are
diagonalizable and $\mathcal{H}_{\rho}$ is split over $\k$.
\end{proof}

As before in this part, we are interested here in whether the irreducible components of a representation of $\mathcal{Q}_{\rho}$ are already irreducible under the action of $\mathcal{H}_{\rho}$. Let $\rho = \rho_1 \oplus \dots \oplus \rho_r$
be a semisimple, essentially pure and \emph{simplicial} representation
of $\AN$. It easily follows from propositions \rref{simplpure}
and \rref{Hrhoreduct} that the decomposition of any tensor product of
$\rho_1,\dots,\rho_r$ and their duals only depends on the decomposition of the reductive Lie algebra $\mathcal{H}_{\rho}$ into simple factors. For applications of this situation to the most classical representations of the
braid group, see \cite{QUOT,LIETRANSP,KRAMMINF}.

These aspects are closely connected to the algebraic hull of $R(P_n)$,
where $R = \widehat{\Phi}(\rho)$, inside $GL_N(K)$. We first need a lemma.

\begin{lemma} Let $\Gamma$ be a Zariski-closed subgroup
of $GL_N(K)$ and $\Lie \Gamma$ its Lie algebra over $K$. For all
$x \in M_N(A)$, if $\exp(hx) \in \Gamma$ then $x \in \Lie \Gamma$.
\end{lemma}
\begin{proof}
Let $a_1,\dots,a_r$ be polynomial functions on $M_N(K)$ with coefficients
in $A$ such that $\Gamma = \{ m \in GL_N(K) \ | \ \forall i \in [1,r] \ \ 
a_i(m)=0 \}$. Since $\Gamma$ is a subgroup of $GL_N(K)$, one has
$a_i(\exp(nhx))=0$ for all $i \in [1,r]$ and $n \in \Z$. Let $Q_i
= a_i(\exp(\alpha h x)) \in (\k[\alpha])[[h]]$. Since $\k$ has caracteristic
0, $\Z$ is Zariski-dense in $\k$ hence $Q_i=0$ and $\exp(\alpha h x)$
is a $\k[[\alpha]]$-point of $\Gamma$. It follows from \cite{CHEVALLEY}
t. 2 \S 8 prop. 4 that $h x \in \Lie \Gamma$ hence $x \in \Lie \Gamma$.
\end{proof}

Recall that $F_n$ is a normal subgroup of $P_n$ such that $P_n/F_n \simeq
P_{n-1}$ and the centralizer of $F_n$ in $P_n$ is $Z(P_n)$. Similarly,
$\mathfrak{F}_n$ is an ideal of $\mathcal{T}_n$, $\mathcal{T}_n/\mathfrak{F}_n
\simeq \mathcal{T}_{n-1}$ and its centralizer in $\mathcal{T}_n$ is the
center of $\mathcal{T}_n$. We let $\check{\mathcal{H}}_{\rho} =
\rho(\mathfrak{F}_n)$.

\begin{prop} Let $H$ and $\check{H}$ be the Zariski closure of $R(P_n)$
and $R(F_n)$ inside $GL_N(K)$, and $\mathfrak{h}$, $\check{\mathfrak{h}}$
be their Lie algebras. Then $\mathcal{H}_{\rho} \otimes K
\subset \mathfrak{h}$ and $\check{\mathcal{H}}_{\rho} \otimes K
\subset \check{\mathfrak{h}}$.
\end{prop}
\begin{proof} 
It easily follows from the definition of $\widetilde{\Phi}$ that
$\widetilde{\Phi}(\xi_{ij}) = \Phi_{ij} e^{2 \la t_{ij}} \Phi_{ij}^{-1}$
with $\Phi_{ij} \in \exp \widehat{\mathcal{T}}_n$, hence
$\widetilde{\Phi}(\xi_{ij}) = \exp \widetilde{\psi}_{ij}$ where
$\widetilde{\psi}_{ij} \in \mathcal{T}_n$ and
$\widetilde{\psi}_{ij} = 2 \la t_{ij}$ plus higher terms. In particular
$R(\xi_{ij}) = \exp h \psi_{ij}$ with $\psi_{ij} \in \rho(\mathcal{T}_n)
\otimes A$ and $\overline{\psi_{ij}} = 2 \la \rho(t_{ij})$. It follows from the
lemma that $\psi_{ij} \in \mathfrak{h}$. Now the $K$-Lie algebra
generated by the $\psi_{ij}$'s is contained in $\mathcal{H}_{\rho} \otimes
K$, and has dimension at least $\dim \mathcal{H}_{\rho}$ since the
$\overline{\psi_{ij}}$'s generate $\mathcal{H}_{\rho}$. It follows
$\mathcal{H}_{\rho} \otimes K \subset \mathfrak{h}$.

Because $\mathcal{F}_n$ is an ideal of $\mathcal{T}_n$ one has
$\widetilde{\psi}_{ij} \in \widehat{\mathcal{F}_n}$ and one gets
$\check{\mathcal{H}}_{\rho} \otimes K
\subset \check{\mathfrak{h}}$ by the same argument.
\end{proof}

In fact we proved more, namely that the same result holds when
$P_n$ and $F_n$ are replaced by the subgroups generated by the
$\xi_{ij}^k$'s, and by the $\xi_{in}^k$'s, for some fixed integer $k$.
As a consequence of this, we get
\begin{prop} Let $G$ be a finite index subgroup of $B_n$, $\Gamma$
the Zariski closure of $R(G)$ and $\mathfrak{h}$ the
Lie algebra of $\Gamma$. Then $\mathcal{H}_{\rho} \otimes K \subset \mathfrak{h}$.
The same holds for any subgroup of $B_n$ containing $\xi_{ij}^k$,
$1 \leq i< j \leq n$, for some fixed non-zero integer $k$.
\end{prop}
\begin{proof}
If $G$ has finite index, there exists $k \in \Z \setminus \{ 0 \}$ such
that G contains the subgroup $P_n^{(k)}$ of $B_n$ generated
by the $\xi_{ij}^k$, $1 \leq i< j \leq n$. Then $R(G) \supset R(P_n^{(k)})$
and the Lie algebra of the Zariski closure of $R(P_n^{(k)})$
contains $\mathcal{H}_{\rho} \otimes K$ by the above remark
hence $\mathcal{H}_{\rho} \otimes K \subset \mathfrak{h}$.
\end{proof}

The last subgoups of interest are the right-angled group
$C_n^{(k)}$ generated by the elements $\sigma_1^{2k},\dots,
\sigma_{n-1}^{2k}$, for some non-zero integer $k$.
\begin{prop} For all $k \in \Z \setminus \{ 0 \}$, let $\Gamma_k$ be the
Zariski closure of $R(C_n^k)$ and $\mathfrak{h}_k$ its Lie algebra.
If $\mathcal{H}_{\rho}$ is generated by the elements $\rho(t_{i,i+1})$
for $1 \leq i \leq n-1$ then
$\mathcal{H}_{\rho} \otimes K \subset \mathfrak{h}_k$
\end{prop}
\begin{proof}
We have $\widetilde{\Phi}(\sigma_i^{2k}) = \exp \widetilde{\psi}_i$
with $\widetilde{\psi}_i \in \widehat{\mathcal{T}_n}$ and
$\widetilde{\psi}_i = 2k \la t_{i,i+1}$
plus higher terms. Hence $R(\sigma_i^{2k}) = \exp h \psi_{i}$ with
$\overline{\psi_i } = 2 k \la \rho(t_{i,i+1})$ and $\psi_{i} \in
\mathfrak{h}_k$. Let $\mathcal{J}$ be the Lie algebra generated
by the $\psi_i$. One has $\mathcal{J} \subset \mathfrak{h}_k$,
$\mathcal{J} \subset \mathcal{H}_{\rho} \otimes K$ and
$\mathcal{J}$ has dimension $\dim \mathcal{H}_{\rho} \otimes K$ since
the elements $\overline{\psi_i}$ generate $\mathcal{H}_{\rho}$.
\end{proof}

\vfill\eject
\section{Basic constructions}

In this last part we review the basic constructions of representations
of the braid groups which appeared in the last decades. We show how
these can be dealt with in our terms, and how our approach -- in particular our uniform approach to unitary representations -- may simplify
or improve either the constructions or the proof of their main properties. 
Incidentally, we shall get three different proofs of the unitarizability of the Burau representation : as a representation of the Iwahori-Hecke algebra (see 6.1.1), of the Birman-Wenzl-Murakami algebra (see 6.2.2) and as a product of Long's induction (see 6.3). For the sake of simplicity, we let $A = \k[[h]]$, $K = \k((h))$, and $\Phi \in \Ass_{\la}(\k)$, $\la \in \kt$ is chosen once and for all.

\subsection{Small representations}

Let $R : B_n \to GL_N(\k)$ be an irreducible representation of $B_n$. To any such representation are attached several integer invariants. Besides the dimension $N$, the most classical ones are the cardinality of the spectrum
and the rank of $R(\sigma_1)-1$. A systematic step-by-step search
for irreducible representations of $B_n$ can thus be carried out by imposing small values on these parameters. First note that, letting $R = \widehat{\Phi}(\rho)$, $\tau = \rho(t_{12})$, $R(\sigma_1) \in GL_N(K)$ is semisimple if
and only if $\tau \in M_N(\k)$ is semisimple.
We now review what happens when one of these
parameters is ``small''.

\subsubsection{Small spectrum}

The first parameter that we consider is the cardinality of the spectrum of $R(\sigma_1)$. We assume that $R(\sigma_1)$ is diagonalizable and let $\mathcal{S} =
\{ a_1,\dots,a_r\} \subset K^{\times}$ be this spectrum. Then $R$ factorizes through
the $K$-algebra $H_n(a_1,\dots,a_r)$ defined as the quotient of $K B_n$ by the relation
$
(\sigma_1 -a_1) (\sigma_1 - a_2) \dots (\sigma_1 - a_r) = 0.
$ If $\mathcal{S}$ has cardinality 1, then $H_n(a_1) = K$ and $R$ factorizes
through a morphism $B_n \to \Z \to GL_1(K)$. If $R = \widehat{\Phi}(\rho)$
the irreducible representation $\rho$
must factorize through a morphism $\AN \to \k \mathfrak{S}_2 \ltimes \k$.
We may thus assume $r \geq 2$.

\paragraph{Iwahori-Hecke algebras.}
If $\mathcal{S} = \{ a_1 ,a_2 \}$ has cardinality 2, the algebra $H_n(a_1,a_2)$ is known as the Iwahori-Hecke
algebra of type $A_{n-1}$ with parameters $a_1,a_2$. It is well-known
that this algebra is isomorphic to $K \SN$ as soon as $a_1/a_2 \in K$ is not
a root of unity and is a square. In particular, if $a_1$ and $a_2$ are algebraically
independant over $\Q$ and $a_1/a_2$ is a square, then $H_n(a_1,a_2) \simeq K \SN$ and,
in case $R = \widehat{\Phi}(\rho)$, the spectrum of $\tau = \rho(t_{12})$ has to be of cardinality 2. We assume that $Sp(\tau) = \{ u,v \} \subset \k$ with $u \neq v$. Then $s = \rho(s_1)$ acts as a scalar on each eigenspace of $\tau$, otherwise $R(\sigma_1)$ would admit three distinct eigenvalues ; 
these scalars are distinct otherwise $\rho$ would be reducible. We can thus assume that
$\tau(x) = ux$ (resp. $\tau(x) = vx$)
iff $s(x) = x$ (resp. $s(x) = -x$), up to tensoring by the sign representation of $\SN$.
It follows that there exists $\alpha,\beta \in \k$
such that $\tau = \alpha + \beta s$.

Conversely, 
let 
$M$ be
an irreducible representation of $\SN$.
For every $\alpha,\beta \in \k$ an element $\tau \in \mathcal{V}(M)$
can be defined by $\tau = \alpha + \beta s$, with $s$ denoting the action of $s_1$ on $M$. If we let $\rho$
be the corresponding representation of $\AN$, then $\widehat{\Phi}(\rho)(\sigma_1)$
is diagonalizable over $K$ with two eigenvalues 
$\{ e^{\la h (\alpha+\beta)} ,- e^{\la h (\alpha - \beta)} \}$,
which are algebraically independant as soon as $\beta \neq 0$,
and it follows that $R$ factorizes through
the Hecke algebra. The cases when $M$ corresponds to the partition $[n-1,1]$ or $[2,1^{n-2}]$ lead to the reduced Burau representation. Note that the existence of a unitary structure
on the Hecke algebra representations (for parameters of modulus 1
and not roots of unity) can be seen in the following way. If $\k \subset \R$, there is a canonical euclidian
structure on $M$ for which the $\SN$-action is orthogonal. In particular,
each transposition is orthogonal and involutive, hence selfadjoint,
and $\tau$ is selfadjoint. This gives then rise to a unitary representation
of $\AN$. Moreover, these representations are simplicial, essentially pure if $\beta \neq 0$, and the elements $Y_n$ are sent to so-called Jucys-Murphy elements of the
symmetric group. It follows that these representations are also agregating.
We determined the algebraic hull of $B_n$ in these representations in \cite{LIETRANSP}.

\paragraph{Cubic Hecke algebras.}

For generic values of $a,b,c$ -- in particular when $a,b,c \in K$
are algebraically independant over $\k$ -- the algebras $H_n(a,b,c)$
for $n \leq 4$ and, to the best of our knowledge conjecturally, $n=5$, are
finite-dimensional and semi-simple. This is related to the fact that they appear
as deformations of group algebras of exceptional finite complex reflections groups.
We first study in some detail the case $n=3$.

The finite complex reflection group called $G_4$ in Shephard and Todd classification (see \cite{SHEPTODD}) admits the presentation
$
<S,T \ \mid \ S^3 = 1, T^3 =1, STS = TST >
$.
It is obviously a quotient of $B_3$ through the map $\sigma_1 \mapsto S$, $\sigma_2
\mapsto T$, and has order 24. Let $j$ be a 3-root of 1. This finite group admits three 1-dimensional irreducible representations $X_{\om}$ for $\om \in \{1,j,j^2\}$ defined by
$S,T \mapsto \omega$, three 2-dimensional irreducible representations $U_{\om}$ characterized by the fact that $S$ and $T$ act with eigenvalues $\{1,j,j^2 \} \setminus \{ \omega \}$, and one 3-dimensional irreducible representation $V$. On the representation $V$, $S$ and $T$ act with eigenvalues $\{1,j,j^2 \}$.
It follows from \cite{BROUEMALLE} that $K G_4$ is isomorphic to its ``cyclotomic Hecke algebra'' $H_3(a,b,c)$ for generic $a,b,c$ and in particular
for $a,b,c$ algebraically independant over $\k$. The irreducible representations of $K G_4$ uniquely deform into irreducible representations $X_{\om}$, $U_{\om}$ and $V$ of $H_3(a,b,c)$ with
$\om \in \{a,b,c\}$. More precisely, a computation of the discriminant shows that $H_3(a,b,c)$ is absolutely semi-simple if and only if
$$
a.b.c.(c-b)(a-c).(b-a).(c^2-cb+b^2).(b^2-ba+a^2).(a^2+bc).(b^2+ac).(c^2+ab) \neq 0. 
$$
A first description of $V$ was obtained in \cite{BROUEMALLE}. We obtained a symmetric model in \cite{THESE} : $\sigma_1$ acts as the diagonal matrix with entries $(a,b,c)$, and
$\sigma_2$ acts as the following matrix~:
$$
\left(
{\begin{array}{ccc}
{\displaystyle \frac {(b + c)\,b\,c}{(a - b)\,(a - c)}}  & 
{\displaystyle \frac {c\,(a\,c + b^2)}{(a - b)\,(a - c)}
}  & {\displaystyle \frac {b\,(a\,b + c^{2})}{(a - b)\,(a - c)
}}  \\ [2ex]
  {\displaystyle \frac {c \,(b\,c + a^{2})}{(b - a)\,(b - c
)}}  &  {\displaystyle \frac {(a + c)\,a\,c}{(b - a)\,(b - c)}
}  &   {\displaystyle \frac {a \,(a\,b + c^{2})}{(b - a)\,(
b - c)}}  \\ [2ex]
{\displaystyle \frac {b\,(b\,c + a^{2})}{(c - a)\,(c - b)}
}  & {\displaystyle \frac {a\,(a\,c + b^2)}{(c - a)\,(c
 - b)}}  & {\displaystyle \frac {(a + b)\,a\,b}{(c - a)\,(c - b)}
} 
\end{array}}
 \right)
$$
The representations $X_{\om}$, $U_{\om}$ arise from Hecke algebra representations, namely the trivial and reduced Burau representations. It follows from \cite{BMR} that all these representations, including $V$, can be obtained by monodromy over the configuration space associated to $G_4$
from the corresponding $G_4$-representations.
A natural question is to ask whether these representations arise
as $\widehat{\Phi}(\rho)$ for some $\rho$. For all irreducible representations except $V$, this is a consequence of
the study of Hecke algebra representations. It is also true for $V$,
due to our work on $\mathcal{V}^s([2,1] \oplus [3])$ in section 4.1.
It follows that every finite-dimensional representation can be obtained
as a direct sum of representations of the form $\widehat{\Phi}(\rho)$,
each twisted by some element in $\Gal(K/\k)$. Note however that some non-irreducible representations cannot be obtained directly as $\widehat{\Phi}(\rho)$ for some $\rho$. This
is in particular the case for the representation $U_a \oplus U_b \oplus U_c$~: the determinant of $\sigma_1$ should be equal to $(abc)^2$ and have $-1$ as constant term at the same time.

For every $n \geq 2$, there exists a finite-dimensional quotient of $H_n(a,b,c)$
known as the Birman-Wenzl-Murakami algebra. From the knowledge of the representation
theory of $G_4$, it is easy to determine from the infinitesimal datas
whether the corresponding monodromy representation factorize through
this algebra or not, and all representations of this algebra
can be obtained by our method -- we refer to \cite{QUOT} for details on these two points, and only mention that the corresponding representations of $\AN$ are simplicial, essentially pure and agregating. 

Important representations of the Birman-Wenzl-Murakami algebra
are the Krammer representations --- essentially the only
known faithful representation of the braid groups. The algebraic
hull of the braid groups inside them was computed in \cite{KRAMMINF}.

\paragraph{Quartic relations and beyond.}

In the same vein, $H_3(a,b,c,d)$ is a deformation of the group algebra of
the finite complex reflection group called $G_8$. It is finite-dimensional,
and isomorphic to $K G_8$ as soon as $a,b,c,d  \in K$ are algebraically independant over $\k$ and $abcd$ admits a square root in $K$. For a detailed study of $H_3(1,u,v,w)$ we refer to the work of Brou\'e and Malle \cite{BROUEMALLE}. The irreducible representations of $G_8$
of dimension less than 4 are not of interest for
us, since they can be obtained by specialization from the cyclotomic Hecke algebra of $G_4$.  The rest is formed by two 4-dimensional irreducible representations. The corresponding two irreducible representations of $H_3(a,b,c,d)$ can be distinguished by the scalar action of $(\sigma_1 \sigma_2)^3$, and can be deduced from each other through the action of
$\Gal(L/\k(a,b,c,d))$, where $L \subset K$ is the minimal (quadratic) extension of $\k(a,b,c,d)$ containing $\sqrt{abcd}$. 
These representations can be obtained by our method. Indeed, we proved in section 4.1 that there exists $\tau \in \mathcal{V}^s([2,1]+[3]+[1^3])$
such that the eigenvalues of $\tau$ can be chosen linearly independant
over $\Q$, hence the corresponding eigenvalues of $\sigma_1$ are
algebraically independent. Moreover, the determinant of $\sigma_1$ has all its square roots inside $\k[[h]]$ --- hence its companion representation is deduced by using some element in $\Gal(K/\k)$. This example shows however that we definitely have to take into account these Galois actions : there is no way to get this other representation directly as $\widehat{\Phi}(\rho)$ for some $\rho$. The algebra $H_3(a,b,c,d,e)$ is again a deformation of
the group algebra of a finite complex reflections group called $G_{16}$. However,
this process does not go further : the quotient of $B_3$ by the relation
$\sigma_1^r = 1$ for $r \geq 6$ is infinite.

\subsubsection{Small dimension}

Small-dimensional irreducible representations of $B_n$ were studied by Formanek and Sysoeva (see \cite{FORMANEK,SYSOEVA1,SYSOEVA2}). In the infinitesimal setting,  
we denote by $M$ a $\k \SN$-module, and ask how small can the dimension of $M$
be with $\mathcal{V}^s(M) \neq \emptyset$. If $M$ decomposes as a
sum of 1-dimensional representations, then $M$ should itself be of dimension 1 under this condition (see proposition \ref{vsvide}). This happens in particular if $M$ is of dimension
at most $n-2$. If $M$ is of dimension at most $n-1$, then
the action of $\SN$ is the standard $(n-1)$-dimensional representation, or its tensorization
by the sign representation and elements in $\mathcal{V}^s(M)$
correspond to the reduced Burau representation. If $M$ has dimension $n$
then $M$ is, up to tensorization by the sign representation, the
sum of the standard representation plus a 1-dimensional representation.
This 1-dimensional representation must be the trivial one by proposition
\ref{vsvide}. The analysis
of $\mathcal{V}^s([n-1,1] \oplus \un )$ carried out in \cite{THESE} shows
that we get the irreducible representation of $\widetilde{\SN}$
described in section 4.3.2. These observations are infinitesimal
analogues of the results of Formanek and Sysoeva. The other representations obtained in their work do not fits our approach, since the eigenvalues of $R(\sigma_1)$ are then algebraic over $\k$ : thay can be obtained only after specialization of the field of matrix coefficients --- for instance, one of them is a $(n-2)$-dimensional representation of the Hecke algebra with parameters roots of unity.

\subsubsection{Small rank}
A central ingredient in Formanek's and Sysoeva's approach is the analysis
of the irreducible representations $R$ such that $R(\sigma_1) - 1$
has rank 1. 
Then $R(\sigma_1)-1$ is either semisimple or nilpotent. 
Let $R = \widehat{\Phi}(\rho)$, $\tau = \rho(t_{12})$ for
$\rho : \AN \to M_N(\k)$.
In case $R(\sigma_1)-1$ is semisimple it factorizes through the Hecke algebra of type $A$, and this general situation was already studied. In this particular case, since $\tau$ is then a linear combination of $\rho(s_1)$ and 1, one easily shows that $R(\sigma_1)-1$ has rank one iff $R$ is the
reduced Burau representation. If it is nilpotent then $R(\sigma_1)-1 \in \rho(s_1)-1 + h M_N(A)$ implies that $\rho$ is the trivial action of $\SN$
on $\k^N$. It thus contradicts our assumption that $R$ is irreducible
for $N \geq 2$. Besides this, it also implies that $R(\sigma_1)-1 \in
\la \rho(s_1) \tau h  + h^2 M_N(A)$, hence that $\tau$ is nilpotent and
of rank 1. A richer infinitesimal analogue of this analysis is then to ask what happens when $\tau$ is a nilpotent transvection. The answer is that any such representation $R$ will factorize through $\widetilde{\SN}$. Since
it is equivalent to the fact that $\rho([t_{12},t_{23}])=0$, this
is an immediate consequence of the following result.

\begin{lemma}
Let $\rho : \mathfrak{B}_3 \to M_N(\k)$ be a representation such that
$\rho(t_{12})$ is a nilpotent transvection. Then $\rho([t_{12},t_{23}]) = 0$.
\end{lemma}
\begin{proof}
For all $1 \leq i < j \leq n$, there exists $v_{ij} \in \k^N \setminus \{ 0 \}$ and a non-zero linear form $\varphi_{i,j}$ on $\k^N$ such that $\rho(t_{ij})(x) = \varphi_{ij}(x) v_{ij}$
for all $x \in \k^N$ and $\varphi_{ij}(v_{ij})=0$. Let us denote by $U \subset \k^N$ the $\k$-vector space spanned
by the $v_{ij}$'s. It is invariant under the action of $\mathcal{T}_3$. Moreover,
for all $x \in \k^N$ and $s \in \mathfrak{S}_3$, 
$$\rho(st_{ij}s^{-1}) (x) = \varphi_{ij} \left( \rho(s^{-1})(x) \right) \rho(s)(v_{ij}) = \rho(t_{s(i) s(j)})(x) = \varphi_{s(i)s(j)}(x) v_{s(i)s(j)}
$$
and it follows that $U$ is invariant under $\mathfrak{B}_3$. Now $U$ is of dimension at most 3. This dimension cannot be $0$ because $v_{12} \neq 0$. If it has
dimension 1, then $v_{12}$ and $v_{23}$ are proportional,
$\varphi_{12}(v_{12}) = 0$ implies $\varphi_{12}(v_{23}) = 0$
and similarly $\varphi_{23}(v_{12}) = 0$, then $\rho([t_{12},t_{23}]) = 0$.
If $U$ has dimension 3, i.e. $v_{12},v_{13}$ and $v_{23}$ are linearly independant, then $\rho([t_{12},t_{13} + t_{23}]) = 0$
implies $\varphi_{12}(x) \varphi_{23}(v_{12}) = 0$ for all $x \in \k^N$,
hence $\varphi_{23}(v_{12}) = 0$. Similarly, $\rho([t_{23},t_{12} + t_{13}]) =0$ implies $\varphi_{12}(v_{23}) = 0$, hence $\rho([t_{12},t_{23}]) = 0$.
The only remaining case is $\dim U = 2$. If $\mathfrak{S}_3$
acts on $U$ by the trivial or sign action then $U$ would be of dimension 1.
If not, $\rho(s_1)$ would admit two distinct eigenvalues and $\rho(t_{12})$, commuting to $\rho(s_1)$, would be semisimple hence zero on $U$. Then $ 0 =\rho(t_{12}) v_{23} = \varphi_{12}(v_{23}) v_{12}$ implies
$\varphi_{12}(v_{23}) = 0$ and similarly $\varphi_{23}(v_{12}) = 0$. This concludes
the proof.
\end{proof}

\subsection{Yang-Baxter representations}

\subsubsection{Preliminaries}
All the material here is standard. We recall that $\k$ has characteristic zero. Let $\g$ be a semisimple
Lie algebra over $\k$, and $K_{\g}$ its Killing form. Since $\g$ is semisimple, this bilinear symmetric invariant form is nondegenerate. Let $(e_{\la})_{\la \in \Lambda}$ be a basis for $\g$, and $(e^{\la})_{\la
\in \Lambda}$ the dual basis of $\g$ with respect to $K_{\g}$ : $K_{\g}(e_{\la},
e^{\mu} ) = \delta_{\la,\mu}$. The Casimir element $C = \sum e_{\la} e^{\la} \in \U \g$ do not depend on the choice of the basis $(e_{\la})$ and is central in $\U \g$. We let $\Delta$ be the coproduct of the Hopf algebra $\U \g$, and introduce $c = \Delta(C) - C \otimes 1 - 1 \otimes C \in (\U \g)^{\otimes 2}$, $\tau = c \otimes 1 \in ( \U \g)^{\otimes n}$. With respect to the action of $\SN$ on $(\U \g)^{\otimes n}$ by permutation of the factors,
$\tau$ is $\mathfrak{S}_2 \times \mathfrak{S}_{n-2}$-invariant (because
$\U\g$ is cocommutative), and we define $\tau_{ij}=\tau_{ji}$ for $1 \leq i < j \leq n$ by the characteristic property
$\tau_{s(1),s(2)} = s. \tau$ for $s \in \SN$. 
It is well known that $c$ commutes to $\Delta(\g) \subset (\U \g)^{\otimes 2}$ and that there exists a $\SN$-equivariant algebra morphism
$\U \mathcal{T}_n \to (\U \g)^{\otimes n}$ defined by $t_{ij} \mapsto \tau_{ij}$. 

As a consequence of this, there is an algebra morphism $\AN \to \k \SN \ltimes (\U \g)^{\otimes n}$. Since $\g$ is semisimple, $\g^n$ is semisimple
and every representation of $(\U \g)^{\otimes n} = \U (\g^n)$ is
completely reducible. Moreover, every irreducible representation of $\g^n$ has the form $V_1 \otimes \dots \otimes V_n$ where the $V_i$'s are irreducible $\g$-modules. By lemma \ref{restrss}, every representation
of $\k \SN \ltimes ( \U \g)^{\otimes n}$ is completely reducible. A special kind of irreducible representation of $\k \SN \ltimes (\U \g)^{\otimes n}$ are of the type $V \otimes \dots \otimes V$ for $V$ an irreducible representation of $\g$. Since $c$ commutes to $\Delta(\g)$, the diagonal
action of $\g$ on any $(\U \g)^{\otimes n}$-module (resp. $\k \SN \ltimes (\U \g)^{\otimes n}$-module) commutes with the action of $\mathcal{T}_n$
(resp. $\AN$).

\subsubsection{Semisimplicity and unitarity}

Let us consider $n$ \emph{self-dual} irreducible $\g$-modules $V_1,\dots,V_n$, for which we fix an isomorphism $V_i \simeq V_i^{\vee}$. Let $\beta_i$ for $1 \leq i \leq n$ be the non-degenerate $\g$-invariant bilinear form on $V_i$ defined by the composite $V_i \otimes V_i \simeq V_i \otimes V_i^{\vee} \to \un$, where $\un$ designates the trivial representation of $\g$ and the last map is the canonical contraction. The form $\beta_i$ is either symmetric or skew-symmetric,
depending on whether the $\g$-module $V_i$ is orthogonal ($ \un \hookrightarrow S^2 V_i$) or symplectic ($ \un \hookrightarrow \Lambda^2 V_i$). We define a non-degenerate bilinear form $\beta$ on $V_1 \otimes \dots \otimes V_n$ as the product of the forms $\beta_i$. It is easily seen that $g_1 \otimes \dots \otimes g_n \in (\U \g)^{\otimes n}$ is selfadjoint with respect to $\beta$ as soon as there are an even number of $g_i \in \g$ while the other ones belong to $\k$. In particular, each $\tau_{ij}$ is selfadjoint with respect to $\beta$ and, if $V_1 = \dots = V_n$, the action of $\SN$ is obviously isometric. This implies that the action of $\mathcal{T}_n$ (resp. $\AN$) is semisimple.

If the form $\beta$ is symmetric and $V_1 = \dots = V_n = V$, then the action of $\AN$ is unitary. This happens for instance if $V$ is the adjoint representation of $\g$, in which case each $\beta_i$ is a multiple of the Killing form and is symmetric,
or if $n$ is even. However the most interesting situation is when
$\k \subset \R$, but then $\beta$ is not in general positive definite.

In one case however we can get unitary representations of $\AN$ on $V^{\otimes n}$ for any field of characteristic zero. Let $\alpha$ be
the canonical symmetric nondegenerate bilinear form on $\k^m$ and $\g=
\so_m(\k)$, $V = \k^m$. The forms $\beta_i$ coincide with $\alpha$, and
$\beta$ is the canonical bilinear form on $V^{\otimes n} = (\k^{m})^{\otimes n} = \k^{mn}$. It follows that the action of $\AN$ and the corresponding action of $B_n$ are unitary. This explains and proves the existence of a unitary structure on the representations of the Birman-Wenzl-Murakami algebra, since they appear in this way. More generally, if each representation $V_i$ of $\g$ is faithful and $\k = \R$, this situation happens exactly when the
Lie algebra $\g$ defined over $\R$ is compact ; indeed, $\beta_i$ can be chosen to be positive definite iff the image of $\g$ in $V_i$ is a
compact Lie algebra.

\subsubsection{Other properties}

By its very definition, the eigenvalues of $\tau$ on $V^{\otimes n}$ are closely related to the eigenvalues of the Casimir operator
on $V \otimes V$. In particular, if the decomposition of $V \otimes V$ as 
a $\g$-module is multiplicity-free and the Casimir operator acts with distinct eigenvalues on each irreducible component, then the action of $\AN$ on $V^{\otimes n}$ is essentially pure. More generally, the eigenvalues of the action of each
$Y_r$ are related to the eigenvalues of the action of the Casimir
operator on $V^{\otimes r}$. We note the following.
\begin{lemma} If $[\tau_{12},\tau_{23}] = 0$ then $\g$ is commutative.
\end{lemma}
\begin{proof} By extension of the base field, we can assume $\bar{\k} = \k$,
hence there exists a basis $(e_{\la})_{\la \in \Lambda}$ of $\g$
which is self-dual with respect for the Killing form and
$ c = 2\sum e_{\la} \otimes e_{\la}$. Then
$
[\tau_{12},\tau_{23}]$ is the sum over all $\la,\mu \in \Lambda$ of the elements $e_{\la} \otimes [e_{\la} , e_{\mu}]
\otimes e_{\mu}$ and the conclusion follows.
\end{proof}
In particular such a representation will in general not factorize through the enhanced symmetric group. In case $\g = \sl_2(\k)$, every representation
is self-dual, so the actions of $\mathcal{T}_n$ will be semisimple. We recall from \cite{IRRED} the following.
\begin{prop}
If $V_1,\dots,V_n$ are irreducible representations of $\sl_2(\k)$, then
the $\k$-vector space of highest weight vectors of $V_1 \otimes \dots \otimes  V_n$ is an agregating representation of $\mathcal{T}_n$.
\end{prop}

\subsection{Artin action and Long induction}

\subsubsection{Artin action}
We denote by $f_1,\dots,f_n$ given free generators of the free group $L_n$.
The well-known faithful Artin action of $B_n$ on $L_n$ is given by making $\sigma_i$ for $1 \leq i < n$ act on $L_n$ by
$f_i \mapsto f_{i+1}$, $f_{i+1} \mapsto f_{i+1}^{-1} f_i f_{i+1}$
and $f_j \mapsto f_j$ for $j \not\in \{ i, i+1 \}$. A theorem of Artin says that this identifies $B_n$ with the subgroup
of automorphisms of $L_n$ which preserve the set of conjugacy classes of $\{ f_1, \dots, f_n \}$ and the product $f_1\dots f_n$. The geometric meaning of this is that $B_n$ can be identified with a set of diffeomorphisms of the disc minus $n$ punctures, whose fundamental group is a free group on $n$ generators. This action enables one
to form the semi-direct product $B_n \ltimes L_n$, which embeds in
$B_{n+1}$ through the inclusion $B_n \subset B_{n+1}$ and the map
$f_i \mapsto \xi_{i,n+1}$, thus identifying $L_n$ and $F_{n+1}$. In particular, every representation
of $B_{n+1}$ restricts to a representation of $B_n \ltimes L_n \simeq B_n \ltimes F_{n+1} \subset B_{n+1}$.
Also
note that $L_n \rtimes B_n$ can be mapped onto $\Z \times B_n$ by sending
each $f_i$ to 1. In particular, every (irreducible) representation
of $B_n$ and any non-zero scalar $\alpha$ yield a (irreducible)
representation of $B_n \ltimes L_n$. In general every representation
$R$ of $B_n \ltimes L_n$ can be rescaled by any non-zero scalar $\alpha$
in order to get a new representation $R^{\alpha}$, letting $R^{\alpha}(\sigma) = R(\sigma)$ for $\sigma \in B_n$, $R^{\alpha}(f_r) = \alpha R(f_r)$. We define here an infinitesimal analogue of Artin action.

Let $\mathcal{L}_n$ be the free Lie algebra on the abstract generators $g_1,\dots,
g_n$, considered as the Lie algebra of the pro-nilpotent completion
of $L_n$. We make $s \in \SN$ act as the automorphism defined by
$g_i \mapsto g_{s(i)}$ for $1 \leq i \leq n$,
and $t_{ij}$ as the derivation defined by
$$
\begin{array}{l}
t_{ij}.g_k = 0 \mbox{ if } k \not\in \{i,j\} \\
t_{ik}.g_k = [g_k,g_i] \\
\end{array}
$$
This leads to a well-defined
action by automorphisms of $\AN$ on $\mathsf{U} \mathcal{L}_n$, in the sense
that grouplike elements in $\AN$ act as automorphisms and primitive elements in $\AN$ act as derivations, because the relation
$
[t_{ij},t_{ik} + t_{kj}] (g_i) = 0$ means $[[g_i,g_j],g_k] - [[g_i,g_k],g_j]
-[g_i,[g_j,g_k]] = 0$, which is a version of
the Jacobi identity. As a consequence of this,
the semi-direct product $\AN \ltimes \mathsf{U} \mathcal{L}_n$ has a well-defined Hopf algebra structure, and the map $g_k \mapsto t_{k,n+1}$ together
with the natural inclusion $\AN \subset \mathfrak{B}_{n+1}$ induce a
Hopf algebra inclusion $\AN \ltimes \mathsf{U} \mathcal{L}_n \subset \mathfrak{B}_{n+1}$, thus identifying $\mathcal{L}_n$ with $\mathcal{F}_{n+1}$. The natural map $L_n \rtimes B_n \to \Z \times B_n$
corresponds to the following. Let $x$ be an indeterminate, and
make $\AN$ act trivially on the 1-dimensional space $\k x$ (i.e. $s.x=x$ for $s \in \SN$, $t_{ij}.x=0$). Now, $\k x$ can be considered as a (commutative) 1-dimensional Lie algebra and the Lie algebra morphism $\mathcal{L}_n \to
\k x$ defined
by $g_i \mapsto x$ is $\AN$-equivariant, hence leads to the
awaited mapping $\mathsf{U} \mathcal{L}_n \rtimes \AN \to \k[x] \times \AN$.
Here we let $\k[x]$ be the ring of polynomials in $x$ and we identify it with the universal envelopping algebra of $\k x$.

We next apply the associator $\Phi$. The image of $\xi_{i,n+1} \in P_{n+1}$
in $\widehat{\mathfrak{B}_{n+1}}$ under the associated Drinfeld morphism belongs to $\widehat{\AN \ltimes \mathsf{U} \mathcal{F}_{n+1}}$. Indeed, $\xi_{i,n+1}$ is a
conjugate of $\sigma_{n}^2$ by elements of $B_n$, and the image of $\sigma_{n}^2$ is
$\Phi(t_{n,n+1}, Y_{n+1}) e^{2\la t_{n,n+1}} \Phi(Y_{n+1}, t_{n,n+1})$ hence belongs to $\widehat{\U \mathcal{F}_{n+1}}$. As a result, we get the following commutative diagram with all arrows injective.
$$
\xymatrix{
B_{n+1} \ar[r] & \widehat{\mathfrak{B}_{n+1}} \\
B_n \ltimes L_n \ar[r] \ar[u] & \widehat{\AN \ltimes \mathsf{U} \mathcal{L}_n} \ar[u] \\
}
$$

\subsubsection{Long induction}

In \cite{LONG}, D.D. Long generalizes the classical Magnus construction to get (unitary)
1-parameter families of representations of $B_n$ from a (unitary) representation of $F_n \ltimes B_n$. In particular, it is possible to
deduce from given representations of $B_{n+1}$, or even representations of $B_n$,
richer representations of $B_n$. We refer to \cite{LONG} for the geometric
interpretation of this construction, and recall the algebraic
construction.

Let $I$ be the augmentation ideal of $\k L_n$, and $R$ a representation of $L_n \ltimes B_n$ on some $\k$-vector space $V$ of finite dimension $m$. There are then well-defined actions of $B_n$ on $I$ and $V$, and the associated representation $R^+$ of $B_n$ on the $nm$-dimensional $\k$-vector space $I \otimes_{\k L_n} V$ is given by making $B_n$ act simultaneously on
both factors. From the above observations, it follows that every
representation $R$ of $B_{n+1}$ (resp. $B_n$), by restriction to
$B_n \ltimes F_n$ and rescaling (resp. by the morphism $B_n \ltimes F_n \to B_n \times \Z$) leads to
a family of $B_n$-representation $(R^{\alpha})^+$. 
Using geometric means and a
theorem of Deligne and Mostow, it is shown in \cite{LONG} that,
if $\k = \C$ and $R$ is unitary, then for generic $\alpha \in \C$ of modulus 1 the representation $(R^{\alpha})^+$ is unitary with respect to some non-degenerate
hermitian form. This hermitian form is not explicitely given at the algebraic level,
and comes from a Poincar\'e duality pairing.

We now give an infinitesimal analogue of this construction. Apart from linearizing the situation, it has the advantage that the infinitesimal bilinear form is explicitely described.
Let $\rho$ be a representation of $\mathsf{U} \mathcal{L}_n \rtimes \AN$,
and $\rho^{\alpha}$ the representation deduced from $\rho$ by the
$\AN$-invariant twisting $g_i \mapsto g_i + \alpha$. We denote by $V$ the
underlying $\k$-vector space. A linear $\AN$-action $\rho^+$
is defined on $V^n$ by making $s \in \SN$ act as  
$$
s.(v_1,v_2,\dots,v_n) = (s.v_{s^{-1}(1)},s.v_{s^{-1}(2)},\dots,
s.v_{s^{-1}(n)})
$$
and, denoting $x^+ = \rho^+(x)$, $t_{ij}^+v = \tilde{v}$ with
$$
\left\lbrace
\begin{array}{lcl}
\tilde{v}_k & = & t_{ij} v_k \mbox{ if } k \not\in \{i,j\} \\
\tilde{v}_i & = & t_{ij} v_i + g_j v_i - g_j v_j \\
\tilde{v}_j & = & t_{ij} v_j + g_i v_j - g_i v_i \\
\end{array} \right.
$$
where $v = (v_1,\dots,v_n)$, $\tilde{v} = (\tilde{v}_1,\dots,\tilde{v}_n)$.
It is easily checked that $\rho^+$ is a representation of $\AN$.  
From the previous observations, we then get representations
$(\rho^{\alpha})^+$ of $\mathfrak{B}_n$ from any representation $\rho$ of $\mathfrak{B}_{n+1}$
or $\AN$.

\subsubsection{Infinitesimal forms}

We assume that a representation $\rho$ of $\AN \ltimes \U \mathcal{L}_n$
on a $\k$-vector space $V$ is fixed, and that $V$ is endowed with a
non-degenerate bilinear form $( \ | \ )$ such that $\rho(s)$ is isometric for every $s \in \SN$. By abuse of notation, we shall write $tv = \rho(t)v$
(resp. $t^+ \underline{v} = \rho^+(t) \underline{v}$) for every $t \in \AN$. We introduce on $V$ the twisted bilinear forms
$\beta_i(v,v') = (g_i v | v')$. We shall denote by $\underline{v} = (v_1,\dots,v_n)$ and $\underline{v}' = (v'_1,\dots v'_n)$ elements
of $V^n$. We define on $V^n$ the bilinear form
$$
(\underline{v}, \underline{v}' ) = \sum_{i=1}^n (g_i v_i|v'_i)
= \sum_{i=1}^n \beta_i( v_i, v'_i).
$$
Last, we decompose each $t_{ij}^+$ into two endomorphisms,
$t_{ij}^+ = d_{ij} + m_{ij}$, with
$$
\begin{array}{lcl}
d_{ij} \underline{v} & = & (t_{ij} v_1, \dots,t_{ij} v_i, \dots, t_{ij} v_j , \dots, t_{ij} v_n) \\
m_{ij} \underline{v} & = & (0,\dots, g_j(v_i - v_j),\dots,g_i(v_j-v_i),\dots,0) \\
\end{array}
$$
We need a technical lemma.
\begin{lemma}
1) For all $s \in \SN$, $\rho^+(s)$ is isometric. 2) If $\rho(t_{ij})^{\dagger} = \pm \rho(t_{ij})$ then, for all $\underline{v},\underline{v}' \in V^n$,
$
(d_{ij} \underline{v},\underline{v}') = \pm (\underline{v},d_{ij} \underline{v}') + ([g_i,g_j] v_j|v'_j) + ([g_j,g_i]v_i|v'_i).
$ 
3) If $\rho(g_i)^{\dagger} = \pm \rho(g_i)$, then
$
(m_{ij} \underline{v},\underline{v}') = \pm (\underline{v},m_{ij} \underline{v}') - ([g_i,g_j] v_j|v'_j) - ([g_j,g_i]v_i|v'_i). 
$
\end{lemma}
\begin{proof}
For $s \in \SN$, $\rho(s)^{\dagger} \rho(s) = 1$ hence, for all $i \in [1,n]$,
and $v,v' \in V$,
$
\beta_i(sv,sv') = (g_i sv,sv') = (sg_{s^{-1}(i)} v | sv') = \beta_{s^{-1}(i)}(v,v')
$
and
$$
(s^+ \underline{v}, s^+ \underline{v}')
= \sum_{i=1}^n \beta_{s(i)} (sv_i,sv'_i) = \sum_{i=1}^n \beta_i(v,v') = (\underline{v},\underline{v}')
$$
so we proved 1). We let $\Delta_{ij} = ([g_i,g_j] v_j|v'_j) + ([g_j,g_i]v_i|v'_i)$. Since $t_{ij}.g_r = 0$ for $r \not\in \{ i,j \}$
and $t_{ij}.g_j = [g_j,g_i]$ we have, for all $v,v' \in V$,
$$
\begin{array}{lcl}
\beta_r(t_{ij}v,v') & = & (g_r t_{ij} v,v') = (t_{ij} g_r v,v') \\
\beta_i(t_{ij} v,v') &  = & (t_{ij} g_i v|v') + ([g_j,g_i] v| v') \\
\beta_j(t_{ij} v,v') &  = & (t_{ij} g_j v|v') + ([g_i,g_j] v| v') \\
\end{array}
$$
It follows, if $\rho(t_{ij})^{\dagger} = \pm \rho(t_{ij})$,
that $\beta_r(t_{ij}v,v') = \pm \beta_r(v,t_{ij}v')$ for $r \not\in \{i,j\}$, and $(d_{ij} \underline{v},\underline{v}') = (\underline{v},d_{ij}
\underline{v}') + \Delta_{ij}$ so 2) is proved. The last assertions are consequences of the following
easy computations, where $u= v_i, w = v_j$, $u' = v'_i, w' = v'_j$ :
$$
\begin{array}{l}
(m_{ij} \underline{v},\underline{v}') = (g_i g_j u|u') - (g_i g_j w,u') + (g_j g_i w | w') - (g_i g_j u | w') \\
(\underline{v}, m_{ij} \underline{v}') = (g_i u| g_j u') - (g_j w| g_i u') + (g_j w,g_i w') - (g_i u|g_j w') \\
\end{array}
$$
\end{proof}

In addition to this lemma, note that the (skew-)symmetricity of
$(\ , \ )$ depends on the (skew-)symmetricity of $(\ | \ )$ and the (anti-)selfadjointness of the $\rho(g_i)$'s. This immediately leads to the following result, where a representation $\rho$ of $\mathcal{L}_n$ is said to be
orthogonal (resp. symplectic, unitary) with respect to some non-degenerate
bilinear form $\beta$ if, in accordance with proposition \rref{defsym}, $\beta$ is symmetric (resp. skew-symmetric, symmetric) and $\rho(g_i)^{\dagger} = -\rho(g_i)$ (resp. $\rho(g_i)^{\dagger} = -\rho(g_i)$,
$\rho(g_i)^{\dagger} = \rho(g_i)$).
\begin{prop}
If the restrictions of $\rho$ to $\AN$ and $\mathcal{L}_n$ are
symplectic (resp. orthogonal, unitary) then $\rho^+$ is orthogonal
(resp. symplectic, unitary).
\end{prop}

Note that, for generic values of $\alpha$ (more precisely for $-\alpha \not\in Sp(g_i)$), the form $(\ , \ )$ is non-degenerate. In the case where $\k \subset \R$ and the restrictions of $\rho$ are unitary, for a
symmetric form $(\ | \ )$ which is positive definite, $( \ , \ )$ is also
positive definite for $\alpha \in \R$ large enough. We then recover unitary
representations of $B_n$.

\vfill\eject
\usection{Appendix A : Formal and real unitarity}
Let $\k$ be one of the topological fields $\R$, $\Q_p$, or one of their finite extensions. In particular, it is endowed with a natural non-trivial valuation.
A power series $f$ in $\k[[h]]$ is called
convergent if it has a non-zero radius of convergence with respect to
the given valuation. 
We let $\k \{ h \}$ be the ring of convergent power series, and by $\k ( \{ h \} )$ its quotient field -- that is, the set of Laurent series $f \in \k (( h ))$ such that $h^r f$ is a convergent power series for some integer $r$. We let $K = \k((h))$, $K^* = \k ( \{ h \} )$ and, forgetting the given topology on $\k$, we endow these fields with the ordinary $h$-adic topology.

The purpose of this appendix is to prove approximation results in
the framework of representation theory. More precisely, we work in the following scope. Let $R : G \to GL_N(K)$ be a representation of a \emph{finitely generated} group $G$. Because of this finite generation hypothesis, the image of $G$ in fact lies in $GL_N(L)$, for $L$ a subfield of $K$ which is finitely generated over $\k$. Let us assume next that there exists a subfield $L^*$ of $K^*$ which is isomorphic to $L$. The representation $R^*$ deduced from $R$ by this isomorphism is convergent, and the two representations $R$ and $R^*$
are essentially equivalent -- they are conjugated by an element in $\Gal(K/\Q)$. Here we prove more. Assuming that $L$ contains
$\k(h)$, we prove the existence of $\k(h) \subset L^* \subset K^*$ such that $L^*$ is isomorphic to $L$ as an extension of $\k(h)$ -- this means that this isomorphism is the identity on $\k(h)$. We also say that $L$ is isomorphic to $L^*$ over $\k(h)$. Then $R$ and $R^*$ are conjugated by an element in $\Gal(K/\k(h))$. 

We next deal with a more subtle question, involving unitarity. We let now
$\k = \R$, and $\eps$ be the automorphism of $K$ defined by $f(h) \mapsto
f(-h)$. This is the only continuous automorphism of $K$ such that $\eps(h) = -h$. It leaves $K^*$ invariant -- we say that $K^*$ is an $\eps$-invariant
subfield of $K$. Using the notations of section 3.2.2, we denote by $U_N^{\eps}(K)$ the group $\{ x \in GL_N(K) \ \mid \ x^{-1} =
 ^t \eps(x) \}$, and assume that we are given a representation $R : G \to U_N^{\eps}(K)$. Again because $G$ is finitely generated,
$R(G)$ lies in $U_N^{\eps}(L)$ for $L$ some finitely generated
subfield of $K$, which may be assumed to contain $\k(h)$ and to be $\eps$-invariant. We prove that there exists a finitely generated subfield
$L^*$ of $K^*$, containing $\k(h)$ and $\eps$-invariant, such that
$U_N^{\eps}(L) \simeq U_N^{\eps}(L^*)$. We then explain how, after a convenient rescaling of the local parameter $h$, we get from this a convergent representation such that specialization in $h$ leads to representations of $G$ into the ordinary unitary group $U_N$.

\bigskip


We first prove the following result.
\begin{theo}
Let $L$ be a subfield of $K$ which contains $\k(h)$ and is finitely
generated over $\k$. Then there exists a subfield $L^*$ of $K^*$ which contains $\k(h)$ and is isomorphic to $L$ over $\k(h)$.

Moreover, for any finite family $a_1,\dots,a_t \in L$ and
$m \geq 0$, then the isomorphism $L \to L^*$ can be chosen
such that each $a_i$ is congruent to its image modulo $h^m$.
\end{theo}
\begin{proof}
Let $\k(h) \subset L_0 \subset L$ be a maximal purely transcendantal extension of $\k(h)$ contained in $L$. Since $L$ is finitely generated over $\k$, its transcendance degree over $\k(h)$ is finite, whence there exists
$f_1,\dots,f_r \in K$ algebraically independant over $\k(h)$ such that
$L_0 = \k(h)(f_1,\dots,f_r)$. The field $L$ is by definition an algebraic
extension of $L_0$, and a finite extension because $L$ is
finitely generated over $\k$. The field $K^*$ is an extension of $\k(h)$
of infinite transcendance degree -- for instance, if $\k = \R, \C$ or
$\Q_p$, the family $\{ e^{h^d} \}_{d > 0}$ is algebraically independant over $\k(h)$ and its elements belong to $K^*$. On the other hand, if $g_1,\dots,g_r \in K^*$ are algebraically independant over $\k(h)$, so are $g_1+P_1,\dots,g_r+P_r$ for $P_1,\dots,P_r$ arbitrary elements of $\k(h)$. It follows from these
two facts that a family $g_1,\dots,g_r$ of convergent power series, algebraically independant over $\k(h)$, can be chosen as close to $f_1,\dots,f_r$ as we want with respect to the $h$-adic topology. Let us now
introduce the abstract field $\tilde{L}_0 = \k(h)(Y_1,\dots,Y_r)$, where $Y_1,\dots,Y_r$ are indeterminates. Every family $g=(g_1,\dots,g_r)$
as above yield an embedding of $\tilde{L}_0$ into $K^*$, through $Y_i \mapsto g_i$. We denote by $L_0^g$ the image subfield of $K^*$. The extensions $\tilde{L}_0$, $L_0$ and $L_0^g$ of $\k(h)$ are by definition isomorphic, and the above considerations show that there exists a family $g_1,\dots,g_r$ such that $L_0^g \subset K^*$, so this proves the result in case $L$ is a purely transcendantal extension. In order to prove the general case,
we make use of the additional fact that $(g_1,\dots,g_r)$ may be
chosen arbitrarily near to $(f_1,\dots,f_r)$. All families $g_1,\dots,g_r$ chosen below are assumed to be algebraically independant over $\k(h)$.

Since $L$ is finite and separable as an extension of $L_0$, the primitive element theorem yield $\alpha \in K$ such that $L= L_0(\alpha)$. Because $L_0$ contains $\k(h)$, we can and will assume $\alpha \in h\k[[h]]$. Let
$P \in L_0[X]$ be a minimal polynomial of $\alpha$ over $L_0$. By not requiring $P$ to be monic, we may assume $P \in \k[h,f_1,\dots,f_r]$. The field $L$ is then isomorphic to $L_0[X]/(P)$. Since $L_0$ is isomorphic
to $\tilde{L}_0$ over $\k(h)$, there exists $\tilde{P} \in \tilde{L}_0[X]$
such that $P$ (resp. its derivative $P'$) is deduced from $\tilde{P}$
(resp. $\tilde{P}'$) through the substitution $Y_i \mapsto f_i$.
Similarly, the substitution $Y_i \mapsto g_i$ leads to polynomials $P_g$ and $P'_g$ in $L_0^g[X]$, such that 
$$
L \simeq L_0[X]/(P) \simeq \tilde{L}_0[X]/(\tilde{P}) \simeq L_0^g[X]/(P_g)
$$
as extensions of $\k(h)$. In order to map $L_0^g[X]/(P_g)$ into $K^*$, it is then sufficient to find convergent $g_1,\dots,g_r$ such that $P_g$ admit a root in $K^*$. Because of Artin's theorem
(see Artin \cite{ARTIN} th. 1.2), it is sufficient for this that $P_g$
admits a root in $h \k[[h]]$. Since $P$ is a minimal polynomial
for $\alpha$ over $L_0$, and $L_0$ has characteristic 0, $P(\alpha) = 0$
and $P'(\alpha) \neq 0$. Since $P'(\alpha) \in \k[[h]]$, there exists than $s\geq 0$, $\beta \in \k \setminus \{ 0 \}$ such that $P'(\alpha) = h^s \beta$ modulo $h^{s+1}$. Choosing $g_1,\dots,g_r \in \k \{ h \}$ close enough to $f_1,\dots,f_r$ so that $P_g(\alpha) = 0$ modulo $h^{2s}$ and $P_g'(\alpha) = h^s \beta$ modulo $h^{s+1}$, we have $P_g(\alpha) \in P'_g(\alpha)^2 h \k[[h]]$.
Hensel's lemma
then asserts (see \cite{EISENBUD}, theorem 7.3) the existence
of $\gamma \in \k[[h]]$ such that $P_g(\gamma) = 0$ and $\gamma - \alpha
\in P_g'(\alpha) h \k[[h]] \subset h\k[[h]]$, whence we conclude that
$\gamma \in h \k[[h]]$ and $L \simeq L_0^g[X]/(P_g) \simeq L_0^g(\gamma) \subset K^*$ as extensions of $\k(h)$.

Now consider a family $a_1,\dots,a_t$ as in the statement, let
$f_1,\dots,f_r$ and $\alpha$ be as above. There exists
polynomials $Q_1,\dots,Q_t$ with coefficients in $\k(h)$ such that
$a_i = Q_i(f_1,\dots,f_r,\alpha)$. By continuity
of these polynomials, there exists $M_0 \in \N$ such that
$f_i \equiv g_i$ and $\alpha \equiv \beta$ modulo $h^{M_0}$
implies $a_i \equiv Q_i(g_1,\dots,g_r,\beta)$ modulo $h^m$. By
requiring the $g_1,\dots,g_r$ chosen above to be congruent
modulo $h^{M_1}$ for some $M_1 \geq M_0$, one may assume
that $s \geq M_0$ hence $\gamma \equiv \alpha$ modulo
$h^{M_0}$. Then Artin's approximation theorem says that the corresponding
root $\beta$ of $P_g$ in $K^*$ can be chosen arbitrarily close to $\gamma$,
and the conclusion follows.
\end{proof}

Recall that we denote by $\eps$ the automorphism of $K$ defined by $f(h) \mapsto f(-h)$.
We will need the following proposition only for $\k = \R$,
although it can be proved in a more general context.

\begin{prop}
Let $L$ be an $\eps$-invariant finitely generated subfield of $K$ containing $\k(h)$. Then there exists an $\eps$-invariant subfield $L^*$ of $K^*$ containing $\k(h)$ and a field isomorphism $\Omega : L \to L^*$
such that $\eps \circ \Omega  = \Omega \circ \eps$.

Moreover, for any finite family $a_1,\dots,a_t \in L$ and
$m \geq 0$, then this isomorphism $L \to L^*$ can be chosen
such that each $a_i$ is congruent to its image modulo $h^m$.
\end{prop}
\begin{proof}
Let $L^{\eps} = \{ x \in L \ \mid \ \eps(x) = x \}$ be the set
of elements of $L$ which are fixed by $\eps$. It obviously contains $\k(h^2)$ and is contained in $K^{\eps} = \k((h^2))$. Then $L$ is a quadratic extension of $L^{\eps}$, $L = L^{\eps} \oplus h L^{\eps}$
as a $L^{\eps}$-vector space and $L \simeq L^{\eps}[X]/(X^2-h^2)$ as a field. We let $\Phi$ be the canonical field isomorphism $\k((h^2)) \to K$ defined by $f(h^2) \mapsto f(h)$,
and let $\Lambda = \Phi(L^{\eps})$. Note that $\Phi$ and $\Phi^{-1}$
send convergent series to convergent series.

We have $L \simeq \Lambda_+ = \Lambda[X]/(X^2-h)$ and the action of $\eps$
on $\Lambda_+$ is the non-trivial element of the Galois group $\Gal(\Lambda_+/\Lambda)$. The theorem claims that there exists a subfield $\Lambda^*$ of $K^*$, which contains $\k(h)$ and is isomorphic to $\Lambda$ over $\k(h)$. Thus there exists a field isomorphism between the field extensions $\Lambda^+/\Lambda$ and $\Lambda_+^*/\Lambda^*$ ; it sends
non-trivial elements in $\Gal(\Lambda/\Lambda_+)$ to non-trivial elements
in $\Gal(\Lambda^*/\Lambda_+^*)$. We finally let
$L^*_- = \Phi^{-1}(\Lambda^*)$  and $L^* = L_-^*(h) \subset K^*$. The composite of this isomorphisms is a field isomorphism $\Omega : L \to L^*$, and it is easily checked that $\Omega \circ \eps = \eps \circ \Omega$.
We summarize the situation by the following commutative diagram in the category of field extensions. All vertical unlabelled arrows represent the only non-trivial element in the Galois group of the corresponding quadratic extensions.
$$
\xymatrix{
L/\k(h^2) \ar[r] \ar[d]_{\eps} & L/L^{\eps} \ar[r]^{\Phi} \ar[d] & \Lambda_+/\Lambda \ar[r] \ar[d] & \Lambda_+^*/\Lambda^*   {\ar[r]^{\Phi^{-1}}} \ar[d] & L^*/L^*_- \ar[d]_{\eps} \\
L/\k(h^2) \ar[r]  & L/L^{\eps} \ar[r]^{\Phi} & \Lambda_+/\Lambda \ar[r] & \Lambda_+^*/\Lambda^* {\ar[r]^{\Phi^{-1}}}& L^*/L^*_- \\
}
$$
For the last part, $a_i \in L = L^{\eps} \oplus h L^{\eps}$ can be decomposed as $a_i = a_i^+ + h a_i^-$. Then the isomorphism
$\Lambda \to \Lambda^*$ can be chosen such that the $a_i^{\pm}$'s equal their images modulo $h^m$
and the conclusion follows.
\end{proof}

We finally use this proposition in the set-up of the introduction. Let $R : G \to U_N^{\eps}(K)$ be a representation of the finitely generated group $G$ into the formal unitary group $U_N^{\eps}(K)$ with $\k = \R$. Our purpose is to deduce from this non-trivial representations of $G$ into the unitary group $U_N$. Let $L$ be the smallest $\eps$-invariant subfield of $K = \R((h))$ containing $\R(h)$ and the coefficients of $R(g)$ for $g \in G$. Because $G$ is finitely generated as a group, and $\eps$ has finite order,
$L$ is finitely generated over $\R$ and we can apply the proposition to this field. Let $L^* \subset K^* = \R(\{h \})$ and $\Omega : L \to L^*$ be the field and field isomorphism given by the proposition. We extend $\Omega : L \to L^*$ coefficientwise to a group isomorphism $GL_N(L) \to GL_N(L^*)$. Since $L$ contains the coefficients of $R$, $R$ factorizes through
$U_N^{\eps}(L) \subset GL_N(L)$, and we have
$$
\begin{array}{lcl}
U_N^{\eps}(L^*) & = & \{ x \in GL_N(L^*) \ \mid \ x^{-1} = \ ^t \eps(x) \} \\
 & = & \{ \Omega(y) \ \mid \ y \in GL_N(L), \  \Omega(y)^{-1} = \ ^t \eps \circ \Omega(y) \}
 \\
 & = & \{ \Omega(y) \ \mid \ y \in GL_N(L), \  \Omega(y^{-1})= \Omega(\ ^t \eps(y)) \} \\
 & = & \Omega(U_N^{\eps}(L)).
 \end{array}
$$
hence $\Omega \circ \rho : G \to U_N^{\eps}(M)$ is a representation of $G$ into $U_N^{\eps}(\R(\{h \}))$. Let $c$
be the automorphism of $\C((h))$ defined by $f(h) \mapsto f(\ii h)$, and $\eta$ the automorphism $f \mapsto
\bar{f}$ induced by the complex conjugation of the coefficients. These two automorphisms leave $\C(\{ h \})$
invariant, and it is readily checked that $\eps$ and $\eta$ coincide on $c(\R((h))) = \R((\ii h))$. In particular, if
we let $J = c(L^*) \subset \C(\{h \})$, then
$$
c \circ \Omega \circ \rho : G \to U_N^{\eps}(K) = U_N^{\eta}(K) \subset U_N^{\eta}(\C(\{ h \}))
$$
hence, by specialization in $h$ real and close to 0, we get morphisms $G \to U_N$, i.e. \emph{unitary representations of $G$} in the ordinary sense. 

\vfill\eject
\usection{Appendix B : Combinatorial aspects}
For a given $\k$-algebra $\mathcal{A}$, a \emph{chain} is a filtration of algebras $\k = \mathcal{A}_1 \subset \dots \subset \mathcal{A}_n = \mathcal{A}$. A representation $\rho$ of $\mathcal{A}$ is called multiply
semisimple (MSS) with respect to this chain if the restriction of $\rho$ to
each $\mathcal{A}_r$, $1 \leq r \leq n$, is semisimple. For a group $G$,
a chain is a filtration of subgroups $\{ e \} \subset G_1 \subset \dots
\subset G_n = G$, and a representation $\rho$ is said to be (MSS) with respect to
this chain if and only if, as a $\k G$-module, it is (MSS) with respect to 
the corresponding filtration of the group algebra.

Under this assumption one can introduce the Bratteli diagram of a representation $\rho$.
We consider this diagram as an oriented graph $\Gamma$ with two distinguished
vertices, the vertex $\mathbf{O}$ corresponding to the trivial representation
of $\k = \mathcal{A}_1$, and the vertex $\rho$. We choose by convention that no edge ends at $\mathbf{O}$. The \emph{level} of a vertex is by definition $1$ plus its distance to $\mathbf{O}$ in the underlying non-oriented graph. A vertex of level $r$ thus corresponds
to an irreducible representation of $\mathcal{A}_r$. We say that
such a Bratteli diagram is \emph{multiplicity free} if every two vertices are connected by at most one edge.

\subsection*{B.1\ \ \ Infinitesimal representations}
\addcontentsline{toc}{subsection}{B.1\ \ \ Infinitesimal representations}
\paragraph{Bratteli diagrams for essentially pure representations.}

We expose a combinatorial device to deal with an irreducible essentially pure $\AN$-representation $\rho$. We moreover assume that the (MSS) condition holds. Note that this condition is automatically satisfied
if $\rho$ is agregating or unitary. One can then
define its Bratteli diagram with respect to the chain
$$
\k \simeq \mathfrak{B}_1 \subset \mathfrak{B}_2 \subset \dots \subset
\mathfrak{B}_n
$$
or to
the chain
$$
\k \simeq \U \mathcal{T}_1 \subset \U \mathcal{T}_2 \subset \dots
\subset \U \mathcal{T}_n.$$
Under our assumption of essential purity, it is clear that one gets the same Bratteli diagram $\Gamma$, choosing one chain or another. We let $T_r$ be the sum of all $t_{ij}$'s in $\mathcal{T}_r$, so that $Y_r = T_r - T_{r-1}$.

We associate to each vertex of level $r$ in $\Gamma$ the action $T_r(\rho')$
of $T_r$ on the corresponding $\mathfrak{B}_r$-representation $\rho'$ : it
belongs to $\End(\rho')$ as a $\mathfrak{B}_r$-module, since $T_r$ belongs to the center of $\mathfrak{B}_r$. We also associate to each edge $p \to q$ with $q$ of level $r$ corresponding to the representation $\rho'$
the action $Y_r(p \to q ) \in \End_{\C}(\rho')$ on the components of (the representation corresponding to) $q$ which are isomorphic to $p$. We then get a \emph{natural coloring} of the sets of vertices and edges of the Bratteli diagram.

\paragraph{Multiplicity free diagrams.}

We now assume that the Bratteli diagram $\Gamma$ of $\rho$ is multiplicity free.
Let $\rho_1,\dots,\rho_r$ be the
irreducible components of the restriction of $\rho$ to $\mathfrak{B}_{n-1}$,
and $Z_n(\rho)$ (resp. $Z_{n-1}(\rho_i)$) the action of $Z_n = 2T_n/n(n-1)$, 
(resp. $Z_{n-1} = 2T_{n-1}/(n-1)(n-2)$) on $\rho$ (resp. $\rho_i$).
If $Z_n(\rho)$ and $Z_{n-1}(\rho_i)$ are \emph{scalars}, then the
following relations hold
$$
\begin{array}{l}
(\dim \rho) Z_n(\rho) = tr(Z_n) = \sum_{i=1}^r tr(Z_{n-1}(\rho_i))
= \sum_{i=1}^r (\dim \rho_i) Z_{n-1}(\rho_i) \\
\dim \rho = \sum_{i=1}^r \dim \rho_i \\
\end{array}
$$
i.e. $(Z_n(\rho), \dim \rho)$ is obtained as the barycentre of the massive
points $(Z_{n-1}(\rho_i),\dim \rho_i)$. It follows that from a Bratteli diagram
and a set of colors on the vertices of level 2 (given by the natural coloring) one can build
a \emph{formal coloring} of the sets of vertices and edges by scalars, defining
$Z_r(\rho)$ by the above formula, $T_r(\rho) = r(r-1)Z_r(\rho)$,
and $Y_r(p \to q) = T_{r}(q) -T_{r-1}(p)$ for $p$ of level $r<n$. It is clear that the formal and natural colorings are the same as soon as
every $T_r$ acts as a scalar on each vertex of level $r$ --- in particular,
if $\k$ is considered as a $\Q$-affine space, the value of $Z_n(\rho)$ then lies inside the convex hull of the spectrum
of $\rho(t_{12})$. We 
say that a coloring of such a graph is \emph{injective} if the
map which to a path from $\mathbf{O}$ to $\rho$ associate the corresponding $n$-tuple of colors of the vertices (or, equivalently, the $(n-1)$-tuple
of colors of the edges) is injective. We then have

\begin{prop}
Let us suppose $\rho$ as above. If the formal coloring of
its Bratteli diagram is injective, then $\rho$ is agregating
and the natural coloring is the same as the formal coloring.
\end{prop}
\begin{proof}
We prove this by induction on $n$. The case $n=2$ is obvious. Let us suppose
$n \geq 3$, and denote by $\rho_1,\dots,\rho_r$ the irreducible components
of the restriction of $\rho$ to $\mathfrak{B}_{n-1}$. 
The coloring of the $\rho_i$'s is injective, then they are agregating and the natural coloring of their diagram is the same as their formal coloring. In particular their natural coloring is a coloring by scalars and is injective. It easily follows (see \cite{IRRED}) that $\rho$ is agregating, hence is absolutely irreducible as stated in proposition \rref{agregtens}.
Then $T_n$ acts as a scalar hence the natural and formal coloring are the same.
We conclude by induction.
\end{proof}

As a corollary, the diagram $\Gamma$ together with given colors on the
level 2 vertices is a sufficient data in this case to determine whether the corresponding monodromy representation $R$ is agregating, or simplicial. It also enables one to check whether $R$ factors through $Hurw_n$, $\Z$ or $\widetilde{\mathfrak{S}}_n$. Indeed, the action of $Y_n$ is readily determined from the coloring, and $R$ factors through $\widetilde{\mathfrak{S}}_n$ iff it factors through $\Z$ iff $\Gamma$ is a segment --- because we assumed $\rho$ to be essentially pure. Note that the
eigenvalues of the $\rho(t_{ij})'s$, the dimension of $\rho$
and the Bratteli diagram of $R$ with respect to $B_1 \subset \dots \subset B_n$ are also determined by these datas.

\subsection*{B.2\ \ \ Braid representations over an algebraic closed field}
\addcontentsline{toc}{subsection}{B.2\ \ \ Braid representations over an algebraic closed field}
We denote by $K$ an algebraically closed field of characteristic zero. Let
$R$ be a (MSS) representation of $B_n$ with respect to the chain $B_1 \subset \dots \subset B_n$, such that the corresponding Bratteli diagram $\Gamma$ is multiplicity free. In the same vein as before, we get a ``natural'' coloring of the vertices (resp. of the edges) from the action of the $\gamma_r$'s (resp. the $\delta_r$'s) --- since $K$ is algebraically closed, $R(\gamma_n)$ is a scalar by Schur's lemma, and $R(\delta_n)=
R(\gamma_n)/R(\gamma_{n-1})$ is split.

We want to know how much information one can get from $\Gamma$ and the natural coloring of the vertices of level 2. First notice that, if no two level-2 vertices have the same natural coloring, it implies that $R(P_n) = R(B_n)$. Indeed, since $\delta_2 = \sigma_1^2$, it means that $R(\sigma_1)$ is a polynomial in $R(\sigma_1^2)$, hence every $R(\sigma_i)$ is a polynomial in $R(\sigma_i^2)$. We now try to recover the natural coloring from these datas. For this purpose, let us introduce $z_n(R) = R(\gamma_n)^{\frac{1}{n(n-1)}} \in K^{\times}/\mu_{\infty}(K)$. Since $\gamma_n = (\sigma_1 \dots \sigma_{n-1})^n$, 
$$\det R(\gamma_n) = (\det R(\sigma_1) )^{n(n-1)} = R(\gamma_n)^{\dim R}
$$
and $z_n(R)^{\dim R}$ equals $\det R(\sigma_1)$ modulo $\mu_{\infty}(K)$. Let us denote by $R_1,\dots,R_r$ the vertices of level $n-1$ in $\Gamma$. Since $\det R(\sigma_1)$ is the product of the
$\det R_i(\sigma_1)$ for $1 \leq i \leq r$, it follows that
$$
z_n(R) = \prod_{i=1}^r z_{n-1}(R_i)^{\frac{\dim R_i}{\dim R}}
$$
and the natural coloring can be recovered up to roots of unity.

Now assume that the natural coloring of $\Gamma$ is known, as well as the representations $R_1,\dots,R_r$. Assuming this coloring known, we
do not need the algebraic closeness assumption anymore. For a complete
description of $R$, the only missing piece is $R(\sigma_{n-1})$.
Since $\sigma_{n-1}$ commutes to $B_{n-2}$, it suffices to determine the action of $\sigma_{n-1}$ on each vector space $\Hom_{B_{n-2}}(\check{R},R)$, for $\check{R}$ a vertex of level $n-2$ in $\Gamma$. The equation $\delta_n = \sigma_{n-1} \delta_{n-1} \sigma_{n-1}$ can be rewritten as $(\sigma_{n-1} \delta_{n-1})^2 = (\delta_{n-1} \sigma_{n-1})^2 = \gamma_n$. Depending on the complexity of $\Gamma$, this yields to a sometimes very tractable
set of equations on each space $\Hom_{B_{n-2}}(\check{R},R)$.
Applications of this method can be found in \cite{THESE}. Its interest heavily depends on the knowledge we have on $\Gamma$, the roots of unity involved in its coloring, and the spectrum of $\sigma_1$.

As we have seen before, these datas are known if $R = \Phi(\rho)$ with $\Phi \in \Ass_{\la}(\k)$, $\la \neq 0$ and $\rho$ an essentially pure irreducible representation of $\AN$. This combinatorial approach has thus the advantage of avoiding the (intricated) calculation of $\Phi(x,y)$ for $x = \rho(t_{i,i+1})$, $y = \rho(Y_i)$ and nevertheless getting an explicit (matrix) description of $R$.

\vfill\eject

\end{document}